\documentclass[5p,twocolumn,number,times]{elsarticle}
\usepackage{amssymb,amsmath,amsthm,epsfig,verbatim,pstricks,url}
\newtheorem{thm}{Theorem}
\newtheorem{lem}[thm]{Lemma}
\newdefinition{rem}[thm]{Remark}
\newproof{pf}{Proof}

\newcommand{\R}{{\mathbb R}}
\newcommand{\spa}{\operatorname{span}}

\newcommand{\diam}{\operatorname{diam}}
\newcommand{\vol}{\mathcal{L}^d} %{\operatorname{vol}}

\newcommand{\dH}[1]{\;{\rm d}{\cal H}^{#1}} % Hausdorff measure
\newcommand{\dL}[1]{\;{\rm d}{\cal L}^{#1}} % Lebesgue measure
\newcommand{\bigchi}{\ensuremath{\mathrm{\mathcal{X}}}}
\newcommand{\charfcn}[1]{\bigchi_{#1}} % characteristic function

\newcommand{\Vh}{\underline{V}(\Gamma^m)}
\newcommand{\Wh}{W(\Gamma^m)}
\newcommand{\Vht}{\underline{V}(\Gamma^h(t))}
\newcommand{\Wht}{W(\Gamma^h(t))}
\newcommand{\uspace}{\mathbb{U}}
\newcommand{\pspace}{\mathbb{P}}
\newcommand{\kspace}{\mathbb{K}}
\newcommand{\xspace}{\mathbb{X}}
\newcommand{\sigmaO}{o}
\newcommand{\nabs}{\nabla_{\!s}}
\newcommand{\Id}{I\!d}

\newcommand{\ddt}{\frac{\rm d}{{\rm d}t}}
\newcommand{\errorXx}{\|\vec{X} - \vec{x}\|_{L^\infty}}

\newcommand{\errorUu}[1]{\|\vec U - \vec I^h_{#1}\,\vec u\|_{L^\infty}}

\newcommand{\LerrorPp}{\|P - p\|_{L^2}}
\newcommand{\LerrorPpc}{\|P_c - p_c\|_{L^2}}
\newcommand{\errorLl}{\|\lambda^h - \lambda\|_{L^\infty}}
\newcommand{\unitn}{\vec{\rm n}}
\newcommand{\XFEMGAMMA}{XFEM$_\Gamma$}
\def\epsilon{\varepsilon} 
\newcommand{\mat}[1]{\underline{\underline{#1}}\rule{0pt}{0pt}}

\begin{document}
\title{Eliminating spurious velocities with a stable approximation \\ of 
viscous incompressible two-phase Stokes flow}
\author[ic]{John W. Barrett}
\ead{j.barrett@imperial.ac.uk}
\author[reg]{Harald Garcke}
\ead{harald.garcke@mathematik.uni-regensburg.de}
\author[ic]{Robert N\"urnberg\corref{rn}}
\ead{robert.nurnberg@imperial.ac.uk}

\address[ic]{Department of Mathematics, 
Imperial College London, London, SW7 2AZ, UK}
\address[reg]{Fakult{\"a}t f{\"u}r Mathematik, Universit{\"a}t Regensburg, 
93040 Regensburg, Germany}
\cortext[rn]{Corresponding author. Tel.: +44 2075948572}

\begin{abstract}
We present a parametric finite element approximation of two-phase flow. This
free boundary problem is given by the Stokes equations in the two
phases, which are coupled via jump conditions across the interface. Using
a novel variational formulation for the interface evolution gives rise to a
natural discretization of the mean curvature of the interface. In addition, 
the mesh quality of the parametric approximation of the interface
does not deteriorate, in general, over time; and an equidistribution property 
can be shown
for a semidiscrete continuous-in-time variant of our scheme in two space
dimensions.
Moreover, 
on using a simple XFEM pressure space enrichment, we obtain exact
volume conservation for the two phase regions.
Furthermore, 
our fully discrete
finite element approximation can be shown to be unconditionally stable. 
We demonstrate the applicability of our method with some numerical
results which, in particular, demonstrate that spurious velocities can
be avoided in the classical test cases.
\end{abstract} 


\maketitle

\thispagestyle{plain}

\section{Introduction} \label{sec:1}

It is well-known that non-physical velocities can appear in 
the numerical approximation of two-phase
incompressible flows. These so-called spurious currents appear in
different representations of the interface, %and they appear also 
with and without surface tension. If surface tension effects are taken
into account, a jump discontinuity in the pressure results, and this
poses serious challenges for the numerical method. As the pressure
jump is balanced by terms involving the curvature of the unknown
interface, it is necessary to accurately approximate the interface, its
curvature and the pressure. In this paper we introduce a new stable
parametric finite element method with good mesh properties, which leads to
an approximation of the interface, and the conditions on it,
with the property that undesired spurious velocities are either small or vanish completely. 
Although B{\"a}nsch, \cite{Bansch01}, proved a stability result for
a (highly nonlinear) 
parametric discretization of the Navier--Stokes equations with a free 
capillary surface, to our knowledge, our stability result, for a linear scheme,
is the first in the 
literature for a parametric discretization of two-phase flow. 

When approximating two-phase flows, one has to decide first of all on
how to represent the interface. The most direct choice is an explicit
tracking of the interface. In these tracking methods 
the interface is either triangulated 
or represented by a connected set of particles, which carry
forces. The interface is
then transported using the flow velocity. Variants of these approaches 
are called front tracking methods, immersed interface methods
or immersed boundary methods, see
e.g.\ \cite{UnverdiT92,LeVequeL97,Bansch01,Peskin02,GanesanT08} %,LiYLSJK13} 
for details.

In a second completely different approach the interface is captured
implicitly by defining a function on the whole domain. In the volume of
fluid (VOF) method 
the characteristic function of one of the fluid phases is used
in order to evolve the interface, see
e.g.\ \cite{HirtN81,RenardyR02,Popinet09}. In the level set method,
instead of a characteristic function, the interface is represented 
as the zero level set of a smooth function, see e.g.\ 
\cite{SussmanSO94,Sethian,OsherF03,GrossR07} for this approach. Finally, in
the phase field approach, instead of a sharp interface description, the
interface is considered to be diffuse with a small interfacial layer
in which a phase field variable rapidly changes from 
the different constant values in the two phases, see e.g.\ 
\cite{AndersonMW98,Jacqmin99,Feng06}. 

Spurious velocities have already been observed in the numerical 
approximation of one-phase
incompressible fluid flow with external forces, see \cite{GerbeauLB97},
where a projection method to deal with this phenomenon is also proposed.
In two-phase flow with surface tension effects it is well-known that the 
imbalance between the discrete computation of the
curvature of the interface and the pressure jump at the interface can
create spurious velocity fields near the interface, even in
situations where the exact solution has zero velocity. %is stationary. 
Several methods to
compute the discrete curvature and different choices of enriching
the pressure space have been introduced to reduce spurious velocities,
see e.g.\ \cite{PopinetZ99,RenardyR02,JametTB02,FrancoisCDKSW06,TongW07,%
GanesanMT07,GrossR11,ZahediKK12,AusasBI12}.

\begin{figure*}
\center
\includegraphics[angle=-0,width=0.33\textwidth]{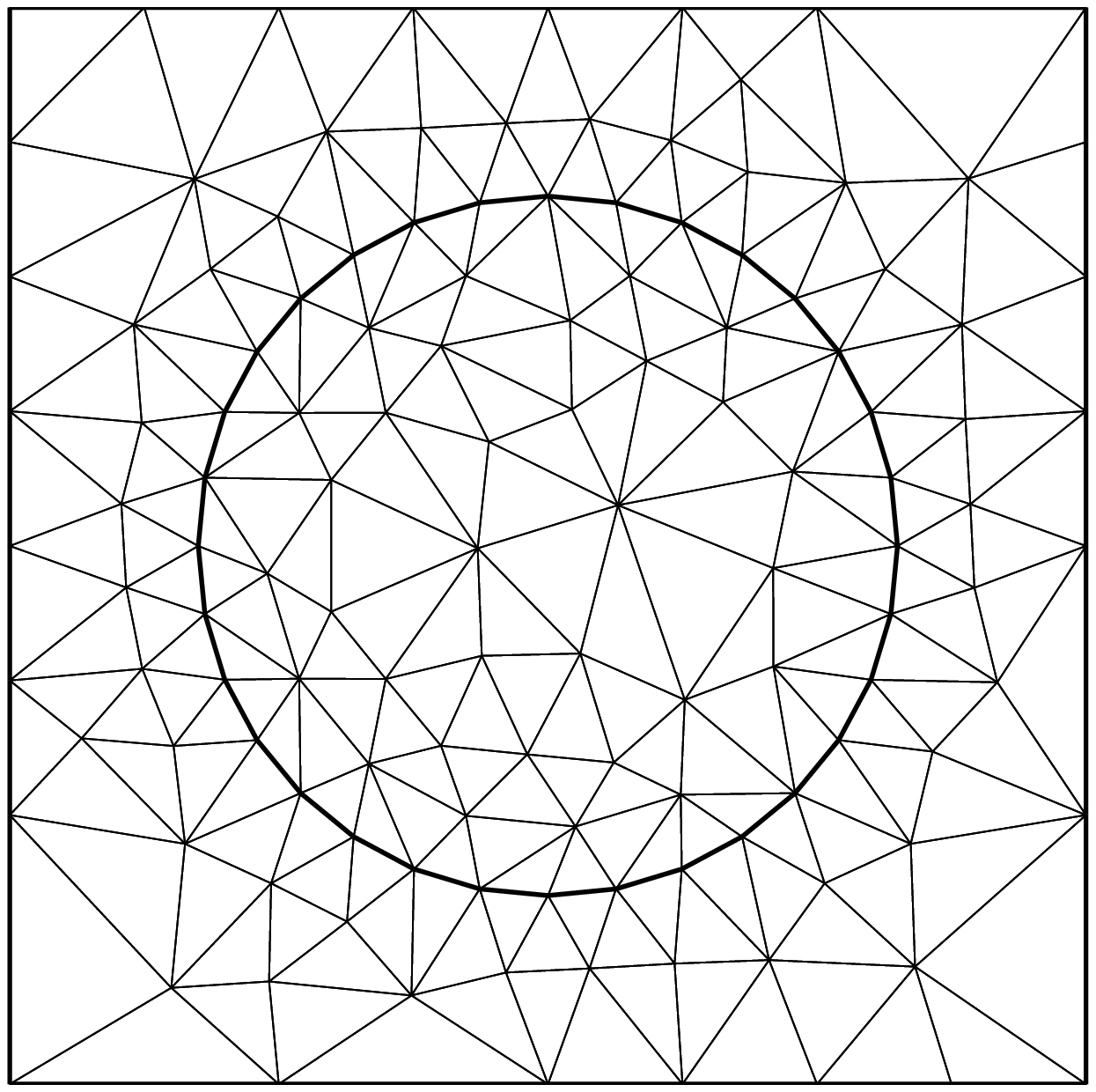}\qquad\quad
\includegraphics[angle=-0,width=0.33\textwidth]{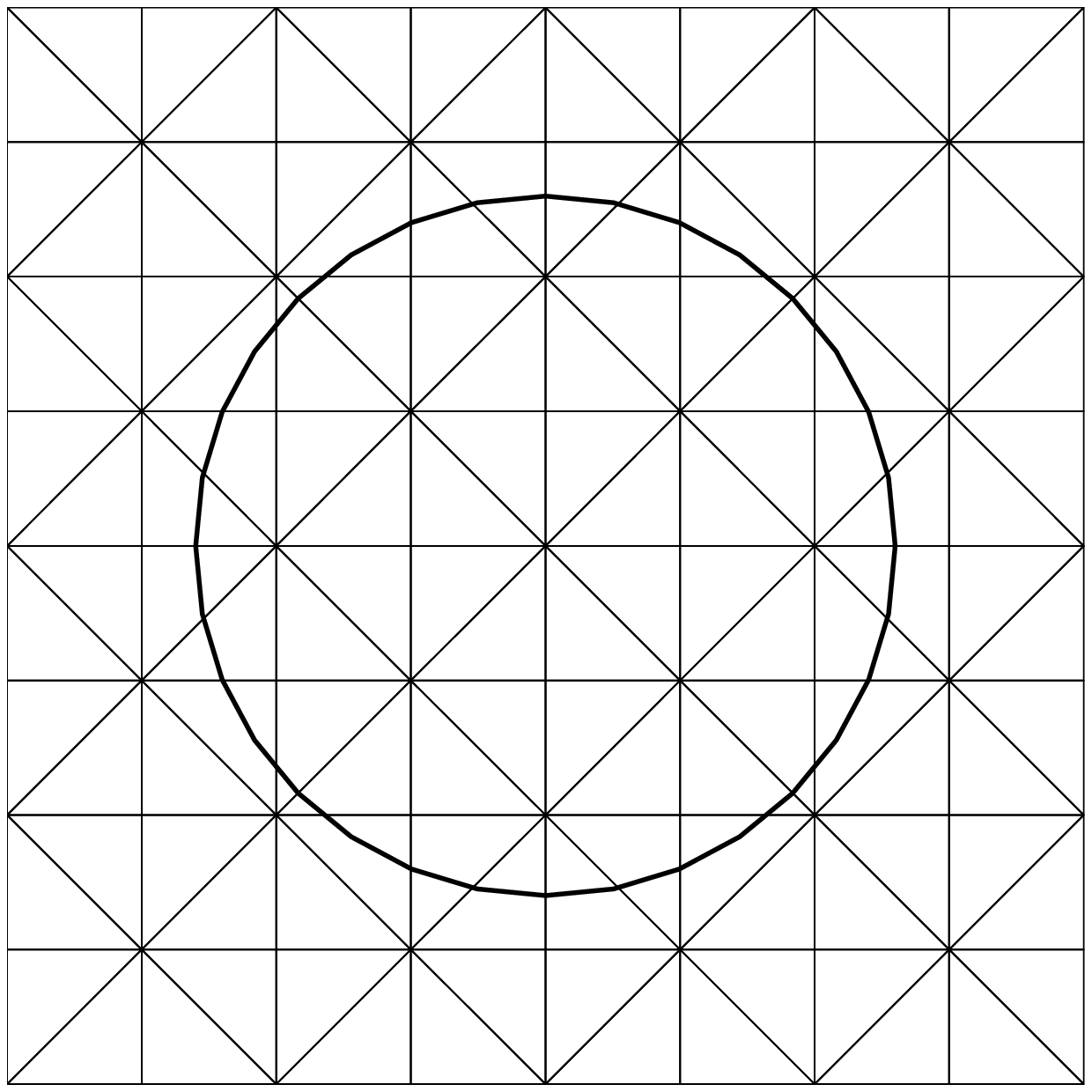}
\caption{Fitted and unfitted bulk finite element meshes for a circular 
interface.}
\label{fig:fitted}
\end{figure*}

In this paper we propose a numerical method for two-phase
incompressible Stokes flow based on a parametric representation of
the interface. We use finite elements to approximate the velocity and the 
pressure, and the interface is approximated using a lower dimensional mesh. 
Here we allow both for a fitted approach,
where the bulk mesh is adapted to the interface, and an unfitted approach,
where the bulk and interface meshes are totally independent. 
Typical meshes for
both approaches are shown in Figure~\ref{fig:fitted}.
In this paper we use the unfitted approach for our numerical results, which
means that we can avoid
remeshings of the bulk mesh at every time step.
As discussed for example in \cite{GanesanMT07}, unfitted bulk %non-adapted
meshes for two-phase flow with pressure jumps, due to surface tension
effects, lead to a poor approximation of the pressure. We avoid
this by using locally refined bulk meshes at the interface in practice.
But we stress that all our theoretical results presented in this paper also
hold for %the principal our approach can be modified also for 
fitted meshes. 
Another strategy in the context of the unfitted approach, 
which can be combined with our proposed method, is to
enrich the pressure finite element space by functions providing
additional degrees of freedom close to the interface; and, in particular,
allow for pressure jumps in the elements cut by the interface. Such an
enrichment technique is an example of the extended finite element method (XFEM),
and has been used for two-phase problems
in the context of level set methods. A major drawback of this approach
is that the resulting linear algebraic system is typically very ill-conditioned,
because the linear independence of the finite element basis
deteriorates as enrichment functions with very small support arise. 
Although strategies  have been developed to reduce the problem of
ill-conditioning, see e.g.\ \cite{GrossR11,AusasBI12,SauerlandF12}, 
the computational effort
still remains large  %updating the mesh 
due to the reconstruction of the XFEM basis as the interface
moves. An XFEM approach is also possible within the context of our
method, and will be discussed later. 

A common problem in approaches which directly parameterize an evolving
interface is that typically the mesh deteriorates, and often
computations cannot be continued without remeshing the interface. 
Situations often occur in which distances between
some interface mesh points or some angles in the interface
triangulation become very small. In earlier work, the present authors
introduced a new methodology to approximate curvature driven curve and
surface evolution, see \cite{triplej,triplejMC,gflows3d}. The method
has the important feature that interface mesh properties remain good during the
evolution. In fact, for curves semidiscrete versions of the approach
lead to polygonal approximations where the vertices are equally
spaced throughout the evolution. The approach has been successfully
used for various geometric evolution equations, including surface diffusion, 
\cite{triplej}, and grain boundary motion, \cite{ejam3d}, and was
recently applied to crystal growth phenomena, see
\cite{dendritic,jcg}. In this paper we generalize the approach to two-phase 
incompressible Stokes flow. %governed by Stokes flow. 
The generalization to the Navier--Stokes case will be considered in the
forthcoming article \cite{fluidfbp}.

Studies of other groups reveal that the main source of spurious
velocities in two-phase flow with surface tension is the fact that
discontinuous functions allowing for jumps at the interface are not in
the pressure space for an unfitted bulk mesh. 
In addition, it has been observed that the size of the
spurious velocities depend both on the calculation of the interface curvature 
in the surface tension term and on the approximation properties of the
pressure space, see e.g.\ \cite{GanesanMT07}. We address
the %first 
issue of obtaining a good approximation of the curvature term
by using the tangential degrees of freedom in the interface representation 
appropriately.
This approach is based on our earlier work on geometric
evolution equations, where we observed good approximation properties
for the curvature also for surfaces in three spatial dimensions, see
\cite{gflows3d,willmore}. 
The %second 
issue concerning the approximation properties of the pressure space
is taken into account by an adaptive refinement of the bulk mesh close to the
interface. In addition, we extend the pressure space by one degree of
freedom by the addition of the characteristic function of one of the
phases. It turns out that this enriched space in a semidiscrete
version leads to exact volume 
conservation (area conservation in 2d) for the two
phase regions. This new pressure space also ensures that stationary spherical
states can be computed %almost 
exactly. %with respect to the fluid velocity.  
In particular, it
turns out that the approach eliminates spurious velocities in the
standard test case of a spherical drop in equilibrium, see e.g.\ 
\cite{GanesanMT07} for the difficulties other approaches have with this
simple test case.
Moreover, in more general situations spurious velocities either
do not appear at all or are small. We also observe that the conditioning
of the resulting linear systems is not so badly effected, in contrast to other
XFEM approaches which involve far more degrees of freedom.

Let us state a few properties of our scheme.

\begin{itemize}
\item The semidiscrete continuous-in-time version of our
scheme is stable in the sense that for a
  closed system the total interface energy decreases in time at
  a rate given by the energy dissipated.
  Similarly, for the fully discrete scheme the rate of decrease
  is at least that given by the energy dissipated. 
  
\item If no outer forces act, then we can easily show that any discrete 
stationary solution must have zero velocity, i.e.\ we
can prove that there are no stationary solutions with spurious
velocities. In addition, we can prove the existence of such stationary
solutions for our scheme.

\item For the semidiscrete %, i.e.\ continuous in time and discrete in space,
  version of our scheme we obtain in two spatial dimensions that the
  interface mesh points are equally spaced. In three space dimensions
  the semidiscrete solutions are conformal polyhedral surfaces, see
  \cite{gflows3d}, which are known to have good mesh properties. In
  practice we also observe for the fully discrete scheme that the
  computed discrete interfaces %polyhedral surfaces 
  have good mesh properties. In particular, no remeshings of the 
  discrete interface are necessary.
   
\item In two spatial dimensions polygonal curves with equidistributed
  vertices on a circle have constant discrete mean curvature and lead
  to discrete solutions with zero velocity and a constant pressure jump across
  the interface. In three spatial dimensions we numerically compute
  stationary polyhedral approximations of a sphere with constant discrete mean
  curvature. These polyhedral surfaces lead to zero velocity solutions
  with a constant pressure jump. 
  These solutions are the natural
  discrete analogues of the stationary circle/sphere, which are the only
  stationary solutions of this incompressible two-phase flow in the
  case when no outer forces act.

\item Our scheme can be applied both in the fitted interface mesh approach, as
well as in the unfitted approach. The latter avoids the repeated 
remeshing of the bulk mesh to ensure that it remains fitted 
to the interface mesh, while the former naturally captures pressure jumps at
the interface with standard pressure finite element spaces.

\end{itemize}

The remainder of the paper is organized as follows. In Section~\ref{sec:2} we
describe the mathematical model of two-phase flow and introduce a suitable 
weak formulation. In Section~\ref{sec:3} we propose our discretization and
establish existence, stability and other theoretical results. Finally, in
Section~\ref{sec:6} we present some numerical results in two and three space
dimensions, including some convergence experiments for stationary and expanding
bubble problems.

\section{Mathematical setting}\label{sec:2}

\begin{figure}
\begin{center}
\unitlength15mm
\psset{unit=\unitlength,linewidth=1pt}
\begin{picture}(4,3.5)(0,0)
\psline[linestyle=solid]{-}(0,0)(4,0)(4,3.5)(0,3.5)(0,0)
\psccurve[showpoints=false,linestyle=solid] 
 (1,1.3)(1.5,1.6)(2.7,1.0)(2.5,2.6)(1.6,2.5)(1.1,2.7)
\psline[linestyle=solid]{->}(2.92,1.6)(3.4,1.6)
\put(3.25,1.7){{\black $\vec\nu$}}
\put(2.9,1.0){{$\Gamma(t)$}}
\put(2,2.1){{$\Omega_-(t)$}}
\put(0.5,0.5){{$\Omega_+(t)$}}
\end{picture}
\end{center}
\caption{The domain $\Omega$ in the case $d=2$.}
\label{fig:sketch}
\end{figure}
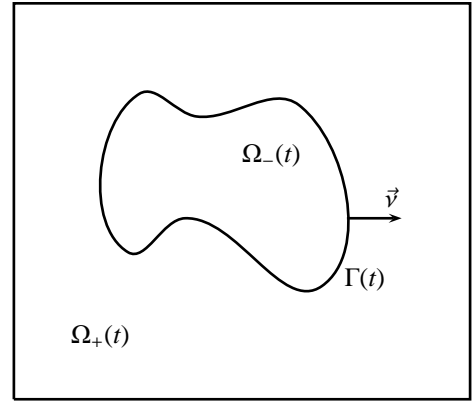

We consider the governing equations for the motion of unsteady,
viscous, incompressible, immiscible two fluid systems. For low
Reynolds numbers one can neglect the inertia terms, and so the equations
for the velocity $\vec u$ and the pressure $p$
are given by
\begin{equation*}
-\mu_{\pm}\, \Delta\,\vec{u} +\nabla\, p =\vec{f} 
\,,\qquad
\nabla\,.\,\vec{u} =0 \qquad \mbox{ in } \Omega_\pm (t)\,;
\end{equation*}
where $\Omega_+(t)$ and $\Omega_-(t)$ are the time dependent regions
occupied by the two fluid phases. For a given domain
$\Omega\subset\mathbb{R}^d$, where $d=2$ or $d=3$, we now seek a time
dependent interface $(\Gamma(t))_{t\in[0,T]}$ such that $\Gamma(t)$ is
completely contained in $\Omega$ and %such that $\Gamma(t)$ 
separates it
into the two domains $\Omega_+(t)$ and $\Omega_-(t)$. Here
the phases could represent two different liquids, or a liquid and a
gas. Common examples are oil/water or water/air interfaces, 
see Figure~\ref{fig:sketch} for an illustration.
For later use, we assume that
$(\Gamma(t))_{t\in [0,T]}$ 
is a sufficiently smooth evolving
hypersurface parameterized by $\vec{x}(\cdot,t):\Upsilon\to\R^d$,
where $\Upsilon\subset \R^d$ is a given reference manifold, i.e.\
$\Gamma(t) = \vec{x}(\Upsilon,t)$. Then
$\mathcal{V} := \vec{x}_t \,.\,\vec{\nu}$ is
the normal velocity of the evolving hypersurface $\Gamma$,
where $\vec\nu$ is the unit normal on $\Gamma(t)$ pointing into $\Omega_+(t)$.
We consider cases in which the viscosity of the two fluids can be
different and introduce 
$\mu(t) = \mu_+\,\charfcn{\Omega_+(t)} + \mu_-\,\charfcn{\Omega_-(t)}$,
with $\mu_\pm \in \R_{>0}$ denoting the fluid viscosities; where, here and
throughout, $\charfcn{\mathcal{A}}$ denotes the characteristic function for a
set $\mathcal{A}$.

In order to define the conditions that have to hold on the interface
$\Gamma(t)$, we introduce the stress tensor
\begin{equation} \label{eq:sigma}
\mat\sigma = \mu \,(\nabla\,\vec u + (\nabla\,\vec u)^T) - p\,\mat\Id
= 2\,\mu\, \mat D(\vec u)-p\,\mat\Id\,,
\end{equation}
where $\mat\Id \in \R^{d \times d}$ denotes the identity matrix and 
$\mat D(\vec u):=\frac12\, (\nabla\vec u+(\nabla\vec u)^T)$ is the
rate-of-deformation tensor. Now on the free surface $\Gamma(t)$, the
following conditions have to hold:
\begin{equation*}
[\vec u]_-^+ = \vec 0\,, \qquad \quad
[\mat\sigma\,\vec \nu]_-^+ = -\gamma\,\varkappa\,\vec\nu\,, 
\qquad \quad
\mathcal{V} = \vec u\,.\,\vec \nu\,, 
\end{equation*}
where $\gamma>0$ is a positive constant and
$\varkappa$ denotes the mean curvature of $\Gamma(t)$, i.e.\ the sum of
the principal curvatures of $\Gamma(t)$. Here we have adopted the sign
convention that $\varkappa$ is negative where $\Omega_-(t)$ is locally convex.
Moreover, as usual, $[\vec u]_-^+ := \vec u_+ - \vec u_-$ and
$[\mat\sigma\,\vec\nu]_-^+ := \mat\sigma_+\,\vec\nu - \mat\sigma_-\,\vec\nu$
denote the jumps in the velocity and the normal component of stress across the interface
$\Gamma(t)$. Here and throughout we employ the shorthand notation
$\vec g_\pm := \vec g\!\mid_{\Omega_\pm(t)}$ for a function 
$\vec g : \Omega \times [0,T] \to \R^d$; and similarly for scalar and
matrix-valued functions.
To close the system we prescribe the initial data 
$\Gamma(0) = \Gamma_0$ 
and the boundary condition 
$\vec u = \vec 0$ %\qquad\mbox{ 
on $\partial \Omega$.
Using the fact that the velocity is divergence free, we can rewrite the
total system as
\begin{subequations}
\begin{alignat}{2}
- \mu\,\nabla\,.\,(\nabla\,\vec u + (\nabla\,\vec u)^T)
+ \nabla\,p & = \vec f \quad &&\mbox{in } 
\Omega_\pm(t)\,, \label{eq:2a} \\
\nabla\,.\,\vec u & = 0 \quad &&\mbox{in } \Omega_\pm(t)\,, \label{eq:2b} \\
\vec u & = \vec 0 \qquad &&\mbox{on } \partial\Omega\,, \label{eq:2c} \\
[\vec u]_-^+ & = \vec 0 \quad &&\mbox{on } \Gamma(t)\,, \label{eq:2d} \\ 
[\mu \,(\nabla\,\vec u + (\nabla\,\vec u)^T)\,\vec\nu - p\,\vec \nu]_-^+ 
& = -\gamma\,\varkappa\,\vec\nu \quad &&\mbox{on } \Gamma(t)\,, \label{eq:2e} \\ 
\mathcal{V} &= \vec u\,.\,\vec \nu \quad &&\mbox{on } \Gamma(t)\,, 
\label{eq:2f}  \\
\Gamma(0) & = \Gamma_0 \,,
\label{eq:2g} 
\end{alignat}
\end{subequations}
which is appropriate for the weak formulation considered in this
paper. 

For the system (\ref{eq:2a}--g) an a priori energy bound holds, and our
goal is to introduce a discretization that satisfies a discrete
analogue. First, on noting (\ref{eq:sigma}) and (\ref{eq:2e}), we have that 
\begin{align}
&
\int_{\Omega_+(t)\cup\Omega_-(t)} (\nabla\,.\,\mat\sigma)\,.\, \vec \xi \dL{d} 
\nonumber \\ 
& \quad = - \int_{\Omega} \mat\sigma : \nabla\,\vec \xi \dL{d}
- \int_{\Gamma(t)} [\mat\sigma\,\vec\nu]_-^+\, . \,\vec \xi  \dH{d-1}  
\nonumber \\ 
& \quad = \int_{\Omega} \left( p\,\nabla\,.\,\vec \xi
-2\,\mu\,\mat D(\vec u) : \mat D(\vec \xi) \right) \dL{d} 
\nonumber \\ & \qquad \qquad
+ \gamma \int_{\Gamma(t)} \varkappa\,\vec \nu \,.\, \vec \xi \dH{d-1}
\quad \forall \ \vec \xi \in H^1_0(\Omega, \R^d)\,,
\label{eq:sigmaibp}
\end{align}
where, here and throughout, $\mathcal{L}^d$ 
and $\mathcal{H}^{d-1}$
denote the Lebesgue measure in $\R^d$ and the $(d-1)$-dimensional Hausdorff measure, respectively.
Then using the identity 
\begin{equation}
\ddt\,\mathcal{H}^{d-1}(\Gamma(t)) 
= -\int_{\Gamma(t)} \varkappa\,\mathcal{V} \dH{d-1} \,,
\label{eq:dtarea}
\end{equation}
see e.g.\ \cite[Lemma~2.1]{DeckelnickDE05},
we obtain from (\ref{eq:2a}--c,f) and (\ref{eq:sigmaibp}) that
\begin{align}
& \gamma\,\ddt\,\mathcal{H}^{d-1}(\Gamma(t))
= -\gamma \int_{\Gamma(t)} %[\mat \sigma\,\vec \nu]^+_-
\varkappa\,\vec \nu \,.\,\vec u \dH{d-1} \nonumber\\
& \qquad =-2\int_\Omega\mu\, \mat D(\vec u):\mat D(\vec u) \dL{d}
+\int_\Omega \vec f\,.\, \vec u \dL{d}\,.\label{eq:testD}
\end{align}
Due to the incompressibility condition, the volume of the two fluids is
preserved. We have from (\ref{eq:2b},f) that
\begin{align}
\frac{\rm d}{{\rm d}t} \vol(\Omega_-(t)) & = 
\int_{\Gamma(t)} \mathcal{V}\dH{d-1}=\int_{\Gamma(t)}\vec
u\,.\,\vec\nu \dH{d-1} \nonumber \\ & = 
\int_{\Omega_-(t)} \nabla\,.\,\vec u \dL{d} 
=0\,. \label{eq:conserved}
\end{align}
It is a further aim that our discretization also preserves the fluid volumes.

In order to compute a discrete version of the mean curvature for
polyhedral surfaces, we note the identity
\begin{equation} \label{eq:LBop}
\Delta_s\, \vec x = \varkappa\, \vec\nu
\qquad \mbox{on $\Gamma(t)$}\,,
\end{equation}
where $\Delta_s = \nabs\,.\,\nabs$ is the Laplace--Beltrami operator on 
$\Gamma(t)$, with $\nabs\,.\,$ and $\nabs$ denoting surface divergence and
surface gradient on $\Gamma(t)$, respectively. We note that the sign convention
in (\ref{eq:LBop}) is such that $\varkappa < 0$ where $\Omega_-(t)$ is locally
convex.
It is possible to use a weak formulation of this identity, which was
first suggested by Dziuk, \cite{Dziuk91}, and was used in \cite{Bansch01} for
the Navier--Stokes equations with a free capillary surface. 
A variant of Dziuk's approach, which leads to good mesh properties, has been
introduced in \cite{triplej} for $d=2$ and in \cite{gflows3d} for
$d=3$, and this will be the basis of our novel weak formulation. The main
novelty, which inherently leads %naturally 
to good mesh properties, 
is that in \cite{triplej,gflows3d} the mean curvature is
treated as a scalar and is discretized separately from the normal $\vec\nu$, 
in contrast to discretizing $\vec\varkappa:=\varkappa\,\vec\nu$ directly 
as in \cite{Dziuk91,Bansch01,GanesanMT07}, see also \S\ref{sec:36} below.

Before introducing our finite element approximation 
of (\ref{eq:2a}--g), we will state an appropriate weak formulation. 
With this in mind, 
we introduce the function spaces
\begin{align*}
\uspace &:= H^1_0(\Omega, \R^d)\,,\qquad \pspace := L^2(\Omega)\,, \\
\widehat\pspace & := \{\eta \in \pspace : \int_\Omega\eta \dL{d}=0 \}\,, \\
\xspace &:= H^1(\Upsilon,\R^d) \qquad\mbox{and}\qquad
\kspace := L^2(\Upsilon,\R)\,, % \label{eq:VW}
\end{align*}
where we recall that $\Upsilon$ is a given reference manifold.  
Let $(\cdot,\cdot)$ and $\langle \cdot, \cdot \rangle_{\Gamma(t)}$
denote the $L^2$--inner products on $\Omega$ and $\Gamma(t)$, respectively.
On recalling (\ref{eq:sigmaibp}), a possible weak formulation of 
(\ref{eq:2a}--g) is then given as follows.
Find time dependent functions $\vec u$, $p$, $\vec x$ and $\varkappa$ such that
$\vec u(\cdot,t) \in \uspace$, $p(\cdot,t) \in \widehat\pspace$,
$\vec{x}(\cdot,t)\in \xspace$, $\varkappa(\cdot,t) \in \kspace$ and
\begin{subequations}
\begin{align}
& 2\left(\mu\,\mat D(\vec u), \mat D(\vec \xi)\right)
- \left(p, \nabla\,.\,\vec \xi\right)
 - \gamma\,\left\langle \varkappa\,\vec\nu, \vec\xi\right\rangle_{\Gamma(t)}
= \left(\vec f, \vec \xi\right) \nonumber \\ & \hspace{6.5cm}
\forall\ \vec\xi \in \uspace \,, \label{eq:weaka}\\
& \left(\nabla\,.\,\vec u, \varphi\right) = 0 
\quad \forall\ \varphi \in \widehat\pspace\,,
\label{eq:weakb} \\
&  \left\langle \vec x_t - \vec u, \chi\,\vec\nu \right\rangle_{\Gamma(t)}
= 0
 \quad \forall\ \chi \in \kspace\,,
\label{eq:weakc} \\
& \left\langle \varkappa\,\vec\nu, \vec\eta \right\rangle_{\Gamma(t)}
+ \left\langle \nabs\,\vec x, \nabs\,\vec \eta \right\rangle_{\Gamma(t)}
= 0  \quad \forall\ \vec\eta \in \xspace\,
\label{eq:weakd}
\end{align}
\end{subequations}
holds for almost all times $t \in (0,T]$, as well as the initial condition (\ref{eq:2g}). 
Here we have observed that if $p \in \pspace$ is part of a
solution to (\ref{eq:2a}--g), then so is $p + c$ for an arbitrary $c
\in \R$.
Moreover, we note that for convenience 
we have adopted a slight abuse of notation
in (\ref{eq:weaka}--d). In particular, 
in this paper we will identify functions defined
on the reference manifold $\Upsilon$ with functions defined on $\Gamma(t)$.
That is, we identify $v \in \kspace$ with $v \circ \vec{x}^{-1}$ on 
$\Gamma(t)$, where we recall that $\Gamma(t) = \vec{x}(\Upsilon,t)$, and we
denote both functions simply as $v$.
For example, $\vec{x} \equiv \vec{\rm id}$ is also 
the identity function on $\Gamma(t)$.

\section{Discretization} \label{sec:3}

We consider the partitioning  $0= t_0 < t_1 < \ldots < t_{M-1} < t_M = T$ 
of $[0,T]$ into possibly variable time steps $\tau_m := t_{m+1}-t_m$, 
$m=0 ,\ldots, M-1$. %We then set $\tau := \max_{m=0,\ldots,  M-1}\tau_m$. 
Let us now introduce the space discretization. 

\subsection{Finite element spaces} \label{sec:31}
For simplicity we consider $\Omega$ to be a polyhedral domain. Then for
all $m\ge 0$, let  ${\cal T}^m$ 
be a regular partitioning of $\Omega$ into disjoint open simplices
$\sigmaO^m_j$, $j = 1 ,\ldots, J^m_\Omega$. 
We set $h^m:= \max_{j=1 ,\ldots, J^m_\Omega}
\mbox{diam}( \sigmaO^m_j)$.
Associated with ${\cal T}^m$ are the finite element spaces
\begin{equation*} % \label{eq:Sh}
 S^m_k := \{\chi \in C(\overline{\Omega}) : \chi\!\mid_{\sigmaO^m} \in
\mathcal{P}_k(\sigmaO^m) \ \forall\ \sigmaO^m \in {\cal T}^m\} 
\subset H^1(\Omega)
\end{equation*}
for $k \in \mathbb{N}$,
where $\mathcal{P}_k(\sigmaO^m)$ denotes the space of polynomials of degree $k$
on $\sigmaO^m$. We also introduce $S^m_0$, the space of  
piecewise constant functions on ${\cal T}^m$.
Then our approximation to the velocity and pressure on ${\cal T}^m$
will be finite element spaces
$\uspace^m\subset\uspace$ and $\pspace^m\subset\pspace$.
We require also the space $\widehat\pspace^m:= \pspace^m \cap \widehat\pspace$. 
The velocity/pressure finite element spaces $(\uspace^m,\pspace^m)$ 
satisfy the LBB inf-sup condition if there exists a
$C_0 \in \R_{>0}$, independent of $h^m$, such that
\begin{equation} \label{eq:LBB}
\inf_{\varphi \in \widehat\pspace^m} \sup_{\vec \xi \in \uspace^m}
\frac{( \varphi, \nabla \,.\,\vec \xi)}
{\|\varphi\|_0\,\|\vec \xi\|_1} \geq C_0 > 0\,,
\end{equation}
where $\|\cdot\|_1 := \|\cdot\|_0 + \|\nabla\,\cdot\|_0$ denotes the 
$H^1$--norm on $\Omega$; %, and
see e.g.\ \cite[p.~114]{GiraultR86}. Here we take the reduced pressure space
$\widehat\pspace^m$ in (\ref{eq:LBB}) because
\begin{equation*} % \label{eq:LBBconst}
\left(1, \nabla\,.\,\vec \xi\right) = \int_{\partial\Omega} \vec\xi\,.\,\unitn
\dH{d-1} = 0 \quad \forall\ \vec \xi \in \uspace^m\,,
\end{equation*}
where $\unitn$ denotes the outer unit normal to $\Omega$.
For example, we may choose the lowest order Taylor--Hood element 
P2--P1, the P2--P0 element or the P2--(P1+P0) element on setting $\uspace^m=[S^m_2]^d\cap\uspace$, 
and $\pspace^m = S^m_1,\,S^m_0$ or $S^m_1+S^m_0$, respectively.
It is
well-known that these choices satisfy the LBB condition (\ref{eq:LBB}). 
We only remark that results for P2--P1 and P2--(P1+P0) need
the weak constraint that all simplices have a vertex in $\Omega$, see
\cite{BoffiCGG12} for the P2--(P1+P0) element.

The parametric finite element spaces in order to approximate $\vec{x}$
and $\varkappa$ in (\ref{eq:weaka}--d) are defined as follows.
Similarly to \cite{gflows3d}, we introduce the following discrete
spaces, based on the work of Dziuk, \cite{Dziuk91}. Let
$\Gamma^{m}\subset\R^d$ be a $(d-1)$-dimensional {\em polyhedral
  surface}, i.e.\ a union of non-degenerate $(d-1)$-simplices with no
hanging vertices (see \cite[p.~164]{DeckelnickDE05} for $d=3$),
approximating the closed surface $\Gamma(t_m)$, $m=0 ,\ldots, M$.  In
particular, let $\Gamma^m=\bigcup_{j=1}^{J^m_\Gamma}
\overline{\sigma^m_j}$, where $\{\sigma^m_j\}_{j=1}^{J^m_\Gamma}$ is a
family of mutually disjoint open $(d-1)$-simplices with vertices
$\{\vec{q}^m_k\}_{k=1}^{K^m_\Gamma}$.
Then for $m =0 ,\ldots, M-1$, let
$\Vh := \{\vec\chi \in C(\Gamma^m,\R^d):\vec\chi\!\mid_{\sigma^m_j}
\in \mathcal{P}_1(\sigma^m_j), j=1,\ldots, J^m_\Gamma\}
=: [\Wh]^d \subset H^1(\Gamma^m,\R^d)$,
where $\Wh \subset H^1(\Gamma^m,\R)$ is the space of scalar continuous
piecewise linear functions on $\Gamma^m$, with 
$\{\chi^m_k\}_{k=1}^{K^m_\Gamma}$ denoting the standard basis of $\Wh$.
For later purposes, we also introduce 
$\pi^m: C(\Gamma^m,\R)\to \Wh$, the standard interpolation operator
at the nodes $\{\vec{q}_k^m\}_{k=1}^{K^m_\Gamma}$,
and similarly $\vec\pi^m: C(\Gamma^m,\R^d)\to \Vh$.
Throughout this paper, we will parameterize the new closed surface 
$\Gamma^{m+1}$ over $\Gamma^m$, with the help of a parameterization
$\vec{X}^{m+1} \in \Vh$, i.e.\ $\Gamma^{m+1} = \vec{X}^{m+1}(\Gamma^m)$.
Moreover, for $m \geq 0$, we will often identify $\vec{X}^m$ with 
$\vec{\rm id} \in \Vh$, the identity function on $\Gamma^m$. 

For scalar and vector %and matrix
functions $v,w$ on $\Gamma^m$ 
we introduce the $L^2$--inner 
product $\langle\cdot,\cdot\rangle_{\Gamma^m}$ over
the current polyhedral surface $\Gamma^m$ %, which is described by the vector
as follows
\begin{equation*} % \label{def:ipm}
\left\langle v, w \right\rangle_{\Gamma^m} := 
\int_{\Gamma^m} v\,.\,w \dH{d-1}\,.
\end{equation*}
If $v,w$ are piecewise continuous, with possible jumps
across the edges of $\{\sigma_j^m\}_{j=1}^{J^m_\Gamma}$,
we introduce the mass lumped inner product
$\langle\cdot,\cdot\rangle_{\Gamma^m}^h$ as
\begin{equation} \label{eq:masslump}
\left\langle v, w \right\rangle^h_{\Gamma^m}  :=
\tfrac1d \sum_{j=1}^{J^m_\Gamma} \mathcal{H}^{d-1}(\sigma^m_j)\,\sum_{k=1}^{d} 
(v\,.\,w)((\vec{q}^m_{j_k})^-),
\end{equation}
where $\{\vec{q}^m_{j_k}\}_{k=1}^{d}$ are the vertices of $\sigma^m_j$,
and where
we define $v((\vec{q}^m_{j_k})^-):=
\underset{\sigma^m_j\ni \vec{p}\to \vec{q}^m_{j_k}}{\lim}\, v(\vec{p})$.

Given $\Gamma^m$, we 
let $\Omega^m_+$ denote the exterior of $\Gamma^m$ and let
$\Omega^m_-$ denote the interior of $\Gamma^m$, so that
$\Gamma^m = \partial \Omega^m_- = \overline\Omega^m_- \cap 
\overline\Omega^m_+$. 
We then partition the elements of the bulk mesh 
$\mathcal{T}^m$ into interior, exterior and interfacial elements as follows.
Let
\begin{align}
\mathcal{T}^m_- & := \{ o^m \in \mathcal{T}^m : o^m \subset
\Omega^m_- \} \,, \nonumber \\
\mathcal{T}^m_+ & := \{ o^m \in \mathcal{T}^m : o^m \subset
\Omega^m_+ \} \,, \nonumber \\
\mathcal{T}^m_{\Gamma^m} & := \{ o^m \in \mathcal{T}^m : o^m \cap
\Gamma^m \not = \emptyset \} \,. \label{eq:partT}
\end{align}
Clearly $\mathcal{T}^m = \mathcal{T}^m_- \cup \mathcal{T}^m_+ \cup
\mathcal{T}^m_{\Gamma^m}$ is a disjoint partition, which in practice
can easily be found e.g.\ with the Algorithm~4.1 in \cite{crystal}. Here we
assume that $\Gamma^m$ has no self intersections, and for the numerical
experiments in this paper this was always the case.
In addition, we define the piecewise constant unit normal 
$\vec{\nu}^m$ to $\Gamma^m$ such that $\vec\nu^m$ points into
$\Omega^m_+$.
Of course, in the case of a fitted bulk mesh it holds that
$\mathcal{T}^m_{\Gamma^m} = \emptyset$.

\subsection{Finite element approximation} \label{sec:32}

Following a similar approach used by the authors in \cite{crystal}
for crystal growth phenomena, 
we consider an unfitted finite element approximation of (\ref{eq:weaka}--d).
On recalling (\ref{eq:partT}), we introduce the discrete viscosity 
$\mu^m \in S^m_0$, for $m\geq 0$, as 
\begin{equation} \label{eq:rhoma}
\mu^m\!\mid_{o^m} = \begin{cases}
\mu_- & o^m \in \mathcal{T}^m_-\,, \\
\mu_+ & o^m \in \mathcal{T}^m_+\,, \\
\tfrac12\,(\mu_- + \mu_+) & o^m \in \mathcal{T}^m_{\Gamma^m}\,.
\end{cases}
\end{equation}
Clearly, for a fitted bulk mesh $\mathcal{T}^m$, (\ref{eq:rhoma}) reduces to
$\mu^m = \mu_+\,\charfcn{\Omega^m_+} + \mu_-\,\charfcn{\Omega^m_-}\in S^m_0$.

Our finite element approximation is then given as follows.
Let $\Gamma^0$, an approximation to $\Gamma(0)$, be given.
For $m=0,\ldots, M-1$, find $\vec U^{m+1} \in \uspace^m$, 
$P^{m+1} \in \widehat\pspace^m$, $\vec{X}^{m+1}\in\Vh$ and 
$\kappa^{m+1} \in \Wh$ such that
\begin{subequations}
\begin{align}
&
 2\left(\mu^m\,\mat D(\vec U^{m+1}), \mat D(\vec \xi) \right)
- \left(P^{m+1}, \nabla\,.\,\vec \xi\right)
 - \gamma\,\left\langle \kappa^{m+1}\,\vec\nu^m,
   \vec\xi\right\rangle_{\Gamma^m}
\nonumber \\ & \hspace{4cm}
= \left(\vec f^{m+1}, \vec \xi\right)
\quad \forall\ \vec\xi \in \uspace^m \,, \label{eq:HGa}\\
& \left(\nabla\,.\,\vec U^{m+1}, \varphi\right)  = 0 
\quad \forall\ \varphi \in \widehat\pspace^m\,,
\label{eq:HGb} \\
&  \left\langle \frac{\vec X^{m+1} - \vec X^m}{\tau_m} ,
\chi\,\vec\nu^m \right\rangle_{\Gamma^m}^h
- \left\langle \vec U^{m+1}, 
\chi\,\vec\nu^m \right\rangle_{\Gamma^m}  = 0
\quad \forall\ \chi \in \Wh\,,
\label{eq:HGc} \\
& \left\langle \kappa^{m+1}\,\vec\nu^m, \vec\eta \right\rangle_{\Gamma^m}^h
+ \left\langle \nabs\,\vec X^{m+1}, \nabs\,\vec \eta \right\rangle_{\Gamma^m}
 = 0
\quad \forall\ \vec\eta \in \Vh
\label{eq:HGd}
\end{align}
\end{subequations}
and set $\Gamma^{m+1} = \vec{X}^{m+1}(\Gamma^m)$.
Here we have defined $\vec f^{m+1}(\cdot) := \vec I^m_2\,\vec
f(\cdot,t_{m+1})$, where $\vec I^m_2$ is the standard interpolation operator
onto $[S^m_2]^d$. Note that here $\nabs$ denotes the surface gradient on
$\Gamma^m$, and so it depends on $m$.
We observe that (\ref{eq:HGa}--d) is a linear scheme in that
it leads to a linear system of equations for the unknowns 
$(\vec U^{m+1}, P^{m+1}, \vec{X}^{m+1}, \kappa^{m+1})$ at each time level.

Furthermore, 
we note that in line with our earlier work in e.g.\ \cite{triplej,gflows3d} 
we employ mass lumping on the first terms in (\ref{eq:HGc},d).
The latter is necessary in order to be able to
prove good mesh properties, while the former is then
enforced to yield a stable scheme. On the other hand, it will become clear in
\S\ref{sec:34} below that for conservation of the phase volumes 
we must use true integration
for the remaining term in (\ref{eq:HGc}), which in turn enforces true
integration for the third term in (\ref{eq:HGa}) for stability reasons.

\subsection{Existence and stability results} \label{sec:33}

Before we can prove existence and uniqueness results for a solution to
our approximation to (\ref{eq:weaka}--d), we have to
introduce the notion of a vertex normal on $\Gamma^m$, which is given
as a weighted sum of the neighbouring normals. We will combine this
definition with a mild assumption that is needed in the existence proof.

\begin{itemize}
\item[$({\cal A})$]
We assume for $m=0,\ldots, M-1$ that $\mathcal{H}^{d-1}(\sigma^m_j) > 0$ 
for all $j=1,\ldots, J^m_\Gamma$,
and that $\Gamma^m \subset \overline\Omega$.
For $k= 1 ,\ldots, K^m_\Gamma$, let
$\Xi_k^m:= \{\sigma^m_j : \vec{q}^m_k \in \overline{\sigma^m_j}\}$
and set
\begin{align*}
\Lambda_k^m & := \bigcup_{\sigma^m_j \in \Xi_k^m} \overline{\sigma^m_j}
 \qquad \mbox{and} \qquad \\
\vec\omega^m_k & := \frac{1}{\mathcal{H}^{d-1}(\Lambda^m_k)}
\sum_{\sigma^m_j\in \Xi_k^m} \mathcal{H}^{d-1}(\sigma^m_j)
\;\vec{\nu}^m_j\,. % \label{omegamj}
\end{align*}
Then we further assume that 
$\dim \spa\{\vec{\omega}^m_k\}_{k=1}^{K^m_\Gamma} = d$,
$m=0,\ldots, M-1$.
\end{itemize}
We stress that $({\cal A})$ is a very mild assumption that is only violated
in very rare occasions. For example, it always holds for surfaces without 
self-intersections; see \cite{triplej,gflows3d} for more details.
Given the above definitions, we introduce the piecewise linear 
vertex normal function
\begin{equation} \label{eq:omega}
\vec\omega^m := \sum_{k=1}^{K^m_\Gamma} \chi^m_k\,\vec\omega^m_k \in \Vh \,,
\end{equation}
and, on recalling (\ref{eq:masslump}), note that 
\begin{subequations}
\begin{equation} 
\left\langle \vec{v}, w\,\vec\nu^m\right\rangle_{\Gamma^m}^h =
\left\langle \vec{v}, w\,\vec\omega^m\right\rangle_{\Gamma^m}^h 
\quad \forall\ \vec{v} \in \Vh\,,\ w \in \Wh \,.
\label{eq:NI}
\end{equation}
and
\begin{equation} 
\left\langle \vec{v}, \vec\omega^m\right\rangle_{\Gamma^m}^h =
\left\langle \vec{v}, \vec\nu^m\right\rangle_{\Gamma^m}^h =
\left\langle \vec{v}, \vec\nu^m\right\rangle_{\Gamma^m}
\quad \forall\ \vec{v} \in \Vh\,. \label{eq:NI1}
\end{equation}
\end{subequations}

The next theorem establishes the existence of a unique solution to 
(\ref{eq:HGa}--d) under some natural conditions. For later developments we also
consider the reduced system {\rm (\ref{eq:HGa},c,d)} with
$\uspace^m$ replaced by 
$$\uspace^m_0 :=\{ \vec U \in \uspace^m : (\nabla\,.\,\vec U, \varphi) = 0 
\quad \forall\ \varphi \in \widehat\pspace^m \} \,.$$
Clearly, any solution $(\vec U^{m+1}, P^{m+1}, \vec{X}^{m+1}, \kappa^{m+1}) 
\in \uspace^m\times \widehat\pspace^m \times \Vh \times \Wh$ to 
{\rm (\ref{eq:HGa}--d)} is such that $(\vec U^{m+1}, 
\vec{X}^{m+1}, \kappa^{m+1}) \in \uspace^m_0 \times \Vh \times \Wh$ solves this
reduced system.

\begin{thm} \label{thm:stab}
Let the assumption ($\mathcal{A}$) hold,
and let $(\uspace^m,\widehat \pspace^m)$ satisfy the LBB condition 
{\rm (\ref{eq:LBB})}, $m=0 ,\ldots, M-1$. 
Then for $m=0 ,\ldots, M-1$ there exists a unique solution
$(\vec U^{m+1}, P^{m+1}, \vec{X}^{m+1}, \kappa^{m+1}) 
\in \uspace^m\times \widehat\pspace^m \times \Vh \times \Wh$ to 
{\rm (\ref{eq:HGa}--d)}. 
In the absence of an LBB condition for $(\uspace^m,\widehat \pspace^m)$ there
exists a unique solution $(\vec U^{m+1}, \vec{X}^{m+1}, \kappa^{m+1}) 
\in \uspace^m_0\times \Vh \times \Wh$ to {\rm (\ref{eq:HGa},c,d)} with
$\uspace^m$ replaced by 
$\uspace^m_0$.
\end{thm}
\begin{pf}
As the system (\ref{eq:HGa}--d) is linear, existence follows from uniqueness.
In order to establish the latter, we consider the system:
Find $(\vec U, P, \vec{X}, \kappa) \in \uspace^m\times\widehat\pspace^m 
\times \Vh \times \Wh$ such that
\begin{subequations}
\begin{align}
&
 2\left(\mu^m\,\mat D(\vec U), \mat D(\vec \xi) \right)
- \left(P, \nabla\,.\,\vec \xi\right)
 - \gamma\,\left\langle \kappa\,\vec\nu^m, \vec\xi\right\rangle_{\Gamma^m}
= 0 %\nonumber \\ & \hspace{6cm}
\quad \forall\ \vec\xi \in \uspace^m \,, \label{eq:proofa}\\
& \left(\nabla\,.\,\vec U, \varphi\right)  = 0 
\quad \forall\ \varphi \in \widehat\pspace^m\,,
\label{eq:proofb} \\
&  \left\langle \frac{\vec X}{\tau_m} ,
\chi\,\vec\nu^m \right\rangle_{\Gamma^m}^h
- \left\langle \vec U, 
\chi\,\vec\nu^m \right\rangle_{\Gamma^m} = 0
 \quad\forall\ \chi \in \Wh\,,
\label{eq:proofc} \\
& \left\langle \kappa\,\vec\nu^m, \vec\eta \right\rangle_{\Gamma^m}^h
+ \left\langle \nabs\,\vec X, \nabs\,\vec \eta \right\rangle_{\Gamma^m}
 = 0  \quad\forall\ \vec\eta \in \Vh\,.
\label{eq:proofd}
\end{align}
\end{subequations}
Choosing $\vec\xi=\vec U$ in (\ref{eq:proofa}), 
$\varphi =  P$ in (\ref{eq:proofb}), 
$\chi = \gamma\,\kappa$ in (\ref{eq:proofc}) 
and $\vec\eta=\gamma\,\vec{X}$ in (\ref{eq:proofd})
yields %, on noting (\ref{eq:NI}), 
that
\begin{equation}
2\,\tau_m\left(\mu^m\,\mat D(\vec U), \mat D(\vec U) \right)
+ \gamma\,\left\langle \nabs\,\vec{X}, \nabs\,\vec{X} \right\rangle_{\Gamma^m} 
=0\,. \label{eq:proof2}
\end{equation}
It immediately follows from (\ref{eq:proof2}) and Korn's inequality
that $\vec U = \vec 0$.
In addition, it holds that $\vec{X} = \vec{X}_c \in \R^d$. 
Together with (\ref{eq:proofc}) for $\vec U=\vec 0$, 
(\ref{eq:NI}) and the assumption $(\mathcal{A})$ this
immediately yields that $\vec{X} = \vec0$, while
(\ref{eq:proofd}) with $\vec\eta=\vec\pi^m[\kappa\,\vec\omega^m]$, 
recall (\ref{eq:NI}), implies that $\kappa = 0$.
Finally, it now follows from (\ref{eq:proofa}) with $\vec U = \vec 0$ and
$\kappa = 0$, on recalling (\ref{eq:LBB}), that $P = 0$. %\in \widehat\pspace^m$.
Hence there exists a unique solution
$(\vec U^{m+1}, P^{m+1}, \vec{X}^{m+1}, \kappa^{m+1}) \in \uspace^m\times
\widehat\pspace^m \times \Vh \times \Wh$ to {\rm (\ref{eq:HGa}--d)}.

The proof for the reduced system {\rm (\ref{eq:HGa},c,d)}, with
$\uspace^m$ replaced by $\uspace^m_0$ is analogous.
\end{pf}

We now demonstrate that the discretization (\ref{eq:HGa}-d) satisfies
an energy estimate, which corresponds to the computation 
(\ref{eq:testD}) in the continuous case. In particular, we obtain
unconditional stability for our scheme.

\begin{thm} \label{thm:stabstab}
Let $0\leq m \leq M-1$ and let
$(\vec U^{m+1},\vec{X}^{m+1}, \kappa^{m+1}) 
\in \uspace^m\times \widehat\pspace^m \times \Vh \times \Wh$ 
be the unique solution to {\rm (\ref{eq:HGa},c,d)},
with $\uspace^m$ replaced by $\uspace^m_0$.
Then
\begin{align}
& \gamma\, \mathcal{H}^{d-1}(\Gamma^{m+1})
+ 2\,\tau_m\left(\mu^m\,\mat D(\vec U^{m+1}), \mat D(\vec U^{m+1}) \right)
\nonumber \\ & \qquad \qquad
\leq \gamma\, \mathcal{H}^{d-1}(\Gamma^m) 
+ \tau_m\left( \vec f^{m+1}, \vec U^{m+1} \right)\,.
\label{eq:stab}
\end{align}
In addition, let $\{t_k\}_{k=0}^M$ be an
arbitrary partitioning of $[0,T]$. Then it holds that
\begin{align}
& \gamma\,\mathcal{H}^{d-1}(\Gamma^{m+1}) 
+ 2 \sum_{k=0}^m  \tau_k\left(\mu^k\,\mat D(\vec U^{k+1}), \mat D(\vec U^{k+1})
\right)%\nonumber \\ & \hspace{2cm}
\nonumber \\ & \qquad \qquad
\leq \gamma\,\mathcal{H}^{d-1}(\Gamma^0)
+ \sum_{k=0}^m \tau_k\left(\vec f^{k+1}, \vec U^{k+1} \right)
\label{eq:stabstab}
\end{align}
for $m=0,\ldots, M-1$.
\end{thm}
\begin{pf}
Choosing $\vec\xi = \vec U^{m+1} \in \uspace^m_0$ in (\ref{eq:HGa}), 
$\chi = \gamma\,\kappa^{m+1}$ in (\ref{eq:HGc}) and
$\vec\eta=\gamma\,({\vec{X}^{m+1}-\vec{X}^m})$ in (\ref{eq:HGd}) yields that
\begin{align*}
&
2\,\tau_m\left(\mu^m\,\mat D(\vec U^{m+1}), \mat D(\vec U^{m+1}) \right)
+ \gamma\,\left\langle \nabs\,\vec{X}^{m+1}, \nabs\,(\vec{X}^{m+1} - \vec{X}^m) 
\right\rangle_{\Gamma^m}
\nonumber \\ & \qquad 
= \tau_m\left(\vec f^{m+1}, \vec U^{m+1} \right)\,.
\end{align*}
Hence (\ref{eq:stab}) follows immediately, where
we have used the result that
$\left\langle \nabs\,\vec{X}^{m+1}, \nabs\,(\vec{X}^{m+1} - \vec{X}^m) 
\right\rangle_{\Gamma^m}
\geq \mathcal{H}^{d-1}(\Gamma^{m+1}) - \mathcal{H}^{d-1}(\Gamma^{m})$,
see e.g.\ \cite[Proof of Theorem~2.3]{triplej} and 
\cite[Proof of Theorem~2.2]{gflows3d} 
for the proofs for $d=2$ and $d=3$, respectively.
The desired result (\ref{eq:stabstab}) immediately follows from 
(\ref{eq:stab}).
\end{pf}

\begin{rem} \label{rem:cons}
It is worthwhile to consider a continuous-in-time semidiscrete version of our 
scheme {\rm (\ref{eq:HGa}--d)}. 
For $t\in [0,T]$, let $\mathcal{T}^h(t)$ be a %n arbitrarily fixed
regular partitioning of $\Omega$ into disjoint open simplices
and define the finite element spaces $S^h_k(t)$, $\uspace^h(t)$ and
$\pspace^h(t)$ similarly to $S^m_k$, $\uspace^m$ and $\pspace^m$, with the
corresponding interpolation operators $I^h_k$ and 
discrete approximations $\mu^h(t) \in S^h_0(t)$, which
will depend on $\Gamma^h(t)$ via the analogue of {\rm (\ref{eq:rhoma})}.
Then, given $\Gamma^h(0)$, %and $\vec U^h(0) \in \uspace^h(0)$, 
for $t\in (0,T]$ find
$\vec U^h(t) \in \uspace^h(t)$, $P^h(t) \in \widehat\pspace^h(t)$,
$\vec{X}^h(t)\in \Vht$ and $\kappa^h(t) \in \Wht$ such that
\begin{subequations}
\begin{align}
&2\left(\mu^h\,\mat D(\vec U^h), \mat D(\vec \xi) \right)
- \left(P^h, \nabla\,.\,\vec \xi\right)
 - \gamma\,\left\langle \kappa^h\,\vec\nu^h,
   \vec\xi\right\rangle_{\Gamma^h(t)}
= \left(\vec f^h, \vec \xi\right)
\nonumber \\ & \hspace{6cm}
\forall\ \vec\xi \in \uspace^h(t) \,, \label{eq:sda}\\
& \left(\nabla\,.\,\vec U^h, \varphi\right)  = 0 
\quad \forall\ \varphi \in \widehat\pspace^h(t)\,,
\label{eq:sdb} \\
&  \left\langle \vec X^h_t ,
\chi\,\vec\nu^h \right\rangle_{\Gamma^h(t)}^h
- \left\langle \vec U^h, 
\chi\,\vec\nu^h \right\rangle_{\Gamma^h(t)} = 0
\quad \forall\ \chi \in \Wht\,,
\label{eq:sdc} \\
&\left\langle \kappa^h\,\vec\nu^h, \vec\eta \right\rangle_{\Gamma^h(t)}^h
+ \left\langle \nabs\,\vec X^h, \nabs\,\vec \eta \right\rangle_{\Gamma^h(t)}
 = 0  
\quad
\forall\ \vec\eta \in \Vht\,,
\label{eq:sdd}
\end{align}
\end{subequations}
where $\vec f^h := \vec I^h_2\,\vec f(t)$.
In {\rm (\ref{eq:sda}--d)} 
we always integrate over the current surface $\Gamma^h(t)$, with normal
$\vec\nu^h(t)$, described by the identity function %parameterization 
$\vec{X}^h(t) \in \Vht$.
Moreover, $\langle \cdot,\cdot\rangle^{h}_{\Gamma^h(t)}$
is the same as $\langle \cdot,\cdot \rangle_{\Gamma^m}^{h}$ with 
$\Gamma^m$ and $\vec{X}^m$ replaced by $\Gamma^h(t)$ and $\vec{X}^h(t)$, 
respectively;
and similarly for $\langle \cdot,\cdot\rangle_{\Gamma^h(t)}$.

Using the results from \cite{gflows3d} it is straightforward to show that
\begin{equation*}
\ddt \mathcal{H}^{d-1}(\Gamma^h(t)) =  
\left\langle \nabs\,\vec{X}^h, \nabs\,\vec{X}^h_t \right\rangle_{\Gamma^h(t)} 
\,,
\end{equation*}
which is the discrete analogue of {\rm (\ref{eq:dtarea})} on noting 
{\rm (\ref{eq:sdd})}.
It is then not difficult to derive the following
energy bound for the solution $(\vec U^h, P^h, \vec{X}^h, \kappa^h)$ of the
semidiscrete scheme {\rm (\ref{eq:sda}--d)}{\rm :}
\begin{align}
& \gamma\,\ddt\,\mathcal{H}^{d-1}(\Gamma^h(t)) 
+ 2\,\|[\mu^h]^\frac12\,\mat D(\vec U^h)\|^2_{0}
= \left(\vec f^h, \vec U^h\right) \,.
\label{eq:stabsd}
\end{align}
Clearly, {\rm (\ref{eq:stabsd})} is the natural discrete analogue of
{\rm (\ref{eq:testD})}.
In addition, it is possible to prove that the vertices of $\Gamma^h(t)$ are
well distributed. As this follows already from the equations 
{\rm (\ref{eq:sdd})}, we
refer to our earlier work in \cite{triplej,gflows3d} for further details. In
particular, we observe that in the case $d=2$, i.e.\ for the planar two-phase
problem, an equidistribution property for the vertices of $\Gamma^h(t)$ can be
shown, while in the case $d=3$ it can be shown that $\Gamma^h(t)$ is a 
conformal polyhedral surface; see also {\rm (\ref{eq:conformal})} below.
\end{rem}

\subsection{\XFEMGAMMA\ for conservation of the phase volumes} \label{sec:34}
In general, the fully discrete approximation (\ref{eq:HGa}--d) will not
conserve mass, which means in particular 
that the volume $\mathcal{L}^d(\Omega^m_-)$, 
enclosed by $\Gamma^m$ will in general not be preserved.
Clearly, given that volume conservation holds on the continuous level, recall
(\ref{eq:conserved}), it would be desirable to %have an analogous property also
preserve the volume of the two phases also on the discrete level.

For the semidiscrete approximation {\rm (\ref{eq:sda}--d)} 
from Remark~\ref{rem:cons} we can show 
conservation of the two phase volumes in the case that 
\begin{equation} \label{eq:XFEM1}
\charfcn{\Omega_-^h(t)} \in \pspace^h(t)\,.
\end{equation}
Choosing $\chi=1$ in {\rm (\ref{eq:sdc})} and
$\varphi=(\charfcn{\Omega_-^h(t)} -
\frac{\mathcal{L}^d(\Omega_-^h(t))}{\mathcal{L}^d(\Omega)})
\in \widehat\pspace^h(t)$ in {\rm (\ref{eq:sdb})}, we
then obtain, on recalling (\ref{eq:NI1}) and as $\vec U^h \in \uspace$, that
\begin{align}
& \frac{\rm d}{{\rm d}t} \vol(\Omega_-^h(t)) = 
\left\langle \vec{X}^h_t , \vec\nu^h \right\rangle_{\Gamma^h(t)}
= \left\langle \vec{X}^h_t , \vec\nu^h \right\rangle^h_{\Gamma^h(t)}
 \nonumber \\ & \
= \left\langle \vec U^h, \vec\nu^h \right\rangle_{\Gamma^h(t)}
= \int_{\Omega_-^h(t)} \nabla\,.\,\vec U^h \dL{d}
= \left( \nabla\,.\,\vec U^h, \charfcn{\Omega_-^h(t)} \right) 
 \nonumber \\ & \
= \left( \nabla\,.\,\vec U^h, \charfcn{\Omega_-^h(t)} -
\frac{\mathcal{L}^d(\Omega_-^h(t))}{\mathcal{L}^d(\Omega)} \right)
=0\,; \label{eq:cons}
\end{align}
which is the discrete analogue of {\rm (\ref{eq:conserved})}. 
Clearly, for fitted bulk meshes $\mathcal{T}^h(t)$ with
$S^h_0(t) \subset \pspace^h(t)$ the condition {\rm (\ref{eq:XFEM1})} trivially
holds. In the case of the unfitted approach, on the other hand, 
discrete pressure spaces $\pspace^h(t)$ 
based on piecewise polynomials, such as
$S^h_0$, $S^h_1$ or $S^h_0 + S^h_1$, will in general not satisfy
the condition {\rm (\ref{eq:XFEM1})}. 
However, the assumption {\rm (\ref{eq:XFEM1})} 
can now be satisfied with the help of the extended finite element method 
(XFEM), see e.g.\ \cite[\S7.9.2]{GrossR11}. 
Here the pressure spaces $\pspace^m$ need to
be suitably extended, so that they satisfy the time-discrete analogue of
(\ref{eq:XFEM1}), i.e.\ 
\begin{equation*} % \label{eq:XFEM1m}
\charfcn{\Omega_-^m} \in \pspace^m\,,
\end{equation*}
which means that then (\ref{eq:HGb}) implies
$\langle \vec U^{m+1}, \vec\nu^m \rangle_{\Gamma^m} = 0$, %\,,
which together with (\ref{eq:HGc}) and (\ref{eq:NI1}) 
then yields that
\begin{equation} \label{eq:consm}
\left\langle \vec X^{m+1} - \vec X^m, \vec \nu^m \right\rangle_{\Gamma^m} =
0\,.
\end{equation}
Hence the obvious strategy to guarantee (\ref{eq:consm}) 
in the context of unfitted bulk meshes is to
add only a single new basis function to the basis of $\pspace^m$, 
namely $\charfcn{\Omega_-^m}$. We remark that in practice
(\ref{eq:consm}) leads to excellent phase volume %mass
conservation properties for the fully discrete scheme (\ref{eq:HGa}--d).
Moreover, we note that the contributions to (\ref{eq:HGa},b)
coming from $\charfcn{\Omega_-^m} -
\frac{\mathcal{L}^d(\Omega_-^m)}{\mathcal{L}^d(\Omega)} \in \widehat\pspace^m$
can be written in terms of integrals over $\Gamma^m$, since
$(\nabla\,.\,\vec \xi, 1 ) = 0$ and
\begin{equation} \label{eq:XFEM1eq}
\left( \nabla\,.\,\vec \xi, \charfcn{\Omega_-^m} \right) =
\int_{\Omega_-^m} \nabla\,.\,\vec \xi \dL{d} = 
\left\langle \vec \xi, \vec\nu^m \right\rangle_{\Gamma^m}
\end{equation}
for all $\vec \xi \in \uspace^m$.
We will call this particular enrichment procedure the \XFEMGAMMA\ approach.
For example, $\pspace^m$ may be given by one of 
\begin{align} \label{eq:PXFEM}
\pspace^m & = S^m_0 + \spa\{\charfcn{\Omega^m_-}\}\,,
\quad
\pspace^m = S^m_1 + \spa\{\charfcn{\Omega^m_-}\} \nonumber \\
\quad\text{or}\quad
\pspace^m & = S^m_0 +S^m_1 + \spa\{\charfcn{\Omega^m_-}\}\,,
\end{align}
with $\mathcal{T}^m$ being independent of $\Gamma^m$.

We note that the above XFEM approach is different to that in e.g.\ 
\cite{GrossR11,AusasBI12,SauerlandF12},
where in the level set context
a standard finite element pressure space is enriched with
numerous discontinuous basis functions in the vicinity of the interface
in order to improve its approximation of the pressure jump
across the interface. Some of these basis functions
may have support on a small fraction of a bulk element 
cut by the interface, and this can lead to ill-conditioning of the
associated linear system. 
Of course, the support of the single additional basis function in the
\XFEMGAMMA\ approach is very non-local, and is given by
$\overline{\Omega^m_-}$. However, its contributions are easily calculated with
the help of surface integrals; recall (\ref{eq:XFEM1eq}). 

Similarly to the other XFEM approaches, we are unable to prove that the 
elements $(\uspace^m, \widehat\pspace^m)$, where $\pspace^m$ is given by one of 
(\ref{eq:PXFEM}), satisfy the LBB condition (\ref{eq:LBB}). This means 
that we cannot easily prove existence and 
uniqueness of the discrete pressure $P^{m+1} \in \widehat\pspace^m$
for the system {\rm (\ref{eq:HGa}--d)} with $\pspace^m$ given as in
(\ref{eq:PXFEM}). However, the relevant result in Theorem~\ref{thm:stab}, 
in the absence of an LBB condition, and Theorem~\ref{thm:stabstab}
still hold.

\subsection{Properties of discrete stationary solutions} \label{sec:35}

We now consider stationary states, 
$\Gamma^{m+1} = \Gamma^m$, of the fully discrete system.
It follows from {\rm (\ref{eq:NI})} that 
a stationary solution to {\rm (\ref{eq:HGa}--d)}
satisfies
\begin{align}\label{eq:conformal}
& \left\langle
\nabs\, \vec X^m,\nabs\,\vec\eta\right\rangle_{\Gamma^m}=0
\nonumber \\ & \quad \forall
\,\vec\eta\in \Vh \ \text{with}\
\vec\eta(\vec{q}_k^m)\,.\,\vec\omega_k^{m} = 0,\, k = 1 ,\ldots, 
K^m_\Gamma \,,
\end{align}
where we note {\rm (\ref{eq:omega})}. We recall from \cite{triplej} that
(\ref{eq:conformal}) in the case $d=2$ implies that $\Gamma^m$ is
equidistributed, with the possible exception of elements $\sigma^m_j$ that are
locally parallel to each other. 
Moreover, we recall from \cite{gflows3d} that surfaces in $\R^3$ that satisfy
(\ref{eq:conformal}) are called conformal polyhedral surfaces.

Next we consider discrete stationary states % in the case where 
when no outer forces act, i.e.\ when $\vec f = \vec 0$. 
Here it turns out that our stability results from \S\ref{sec:33} have an
immediate consequence.

\begin{lem}\label{lem:stat1}
Let
$(\vec U^{m+1},P^{m+1},\vec X^{m+1}, \kappa^{m+1}) 
\in
\uspace^m\times \widehat \pspace^m\times \Vh\times\Wh$ 
be a solution to {\rm (\ref{eq:HGa}--d)} with 
$\vec f^{m+1} = \vec 0$.
If $\vec X^{m+1} = \vec X^m$, then $\vec U^{m+1} = \vec 0$. 
\end{lem}
\begin{pf}
On recalling Theorem~\ref{thm:stabstab}, %and \ref{thm:stabXFEM1},
the solution $(\vec U^{m+1}, \vec X^{m+1})$ fulfills
(\ref{eq:stab}) with $\Gamma^{m+1}$ replaced by $\Gamma^m$ and 
$\vec f^{m+1} = \vec 0$. 
Hence we obtain $(\mu^m\,\mat D (\vec U^{m+1}), \mat D(\vec U^{m+1})) = 0$, 
and so Korn's inequality implies $\vec U^{m+1} = \vec 0$. 
\end{pf}

The above lemma guarantees, independently of the choice of $\mu_\pm$,
that no spurious velocities appear for
discrete stationary solutions, $\Gamma^{m+1} = \Gamma^m$. 
For the \XFEMGAMMA\ approach we
can even show that polyhedral surfaces with constant discrete mean
curvature and zero velocity are stationary solutions. 

\begin{lem} \label{lem:stat2}
Let $\charfcn{\Omega_-^m} \in \pspace^m$ and
let $\Gamma^m$ be a polyhedral surface with constant discrete mean
curvature, i.e.\ there exists a constant
$\overline{\kappa}\in\mathbb{R}$ such that
\begin{equation}\label{eq:constcurv}
\overline{\kappa} \left\langle\vec\nu^m,\vec\eta\right\rangle_{\Gamma^m}
+\left\langle
\nabs\, \vec X^m,\nabs\,\vec\eta\right\rangle_{\Gamma^m}=0\quad \forall
\,\vec\eta\in\Vh\,.
\end{equation}
Then $\Gamma^m$ satisfies {\rm (\ref{eq:conformal})} and %Then 
$(\vec U^{m+1}, \vec X^{m+1}, \kappa^{m+1}) = 
(\vec 0, \vec X^m, \overline{\kappa})\in\uspace^m_0\times\Vh\times \Wh$ 
is the unique solution to the reduced system {\rm (\ref{eq:HGa},c,d)}, 
with $\uspace^m$ replaced by $\uspace^m_0$,
with $\vec f^{m+1} = \vec 0$. 
\end{lem}
\begin{pf} %Assertion (i) follows as in Lemma~\ref{lem:stat1}. 
It immediately follows from (\ref{eq:NI}) that (\ref{eq:conformal}) holds.
Theorem~\ref{thm:stab} implies that in order to establish the 
remaining %desired
result, %(ii), 
we only need to show that $(\vec U^{m+1}, \vec X^{m+1}, \kappa^{m+1}) = 
(\vec 0, \vec X^m, \overline{\kappa})$ is a solution to 
{\rm (\ref{eq:HGa},c,d)}, with $\uspace^m$ replaced by $\uspace^m_0$,
with $\vec f^{m+1} = \vec 0$. But this follows immediately from
$\overline{\kappa}\,\langle\vec\nu^m,\vec\eta\rangle_{\Gamma^m}
= \langle\overline{\kappa}\,\vec\nu^m,\vec\eta\rangle_{\Gamma^m}^h$
for all $\vec\eta \in \Vh$, 
and
\begin{equation*}
\left\langle \vec\nu^m, \vec\xi \right\rangle_{\Gamma^m} = 
\left(\nabla\,.\,\vec\xi,\charfcn{\Omega_-^m}\right) =0
\quad \forall \ \vec\xi \in \uspace^m_0\,,
\end{equation*}
where we have recalled (\ref{eq:XFEM1eq}).
\end{pf}

\begin{rem} \label{rem:stat}
A stationary solution to the continuous problem with $\vec f = \vec 0$ is
a circle $(d=2)$ or a sphere $(d=3)$ with zero velocity
and a piecewise constant pressure with a discontinuity across the interface,
see {\rm (\ref{eq:radialr},b)} below. 

For $d=2$, one can choose $\Gamma^m$ with equidistributed points on a circle 
as an approximation of this circle, i.e.\ a closed regular polygon.
Such a $\Gamma^m$ has constant discrete curvature, i.e.\ 
there exists a $\overline{\kappa} \in \mathbb{R}$
such that {\rm (\ref{eq:constcurv})} is satisfied. 
Hence Lemma~{\rm \ref{lem:stat2}} yields that in this situation
$(\vec U^{m+1}, \vec X^{m+1}, \kappa^{m+1}) = (\vec 0, \vec X^m, 
\overline{\kappa})$ %\in\uspace^0\times \Vhz\times \Whz$ 
is the unique solution to the reduced system
with $\vec f^{m+1} =\vec 0$.   

For $d=3$, 
we observe in practice that conformal approximations of the sphere,
i.e.\ spherical $\Gamma^m$ satisfying {\rm (\ref{eq:conformal})}, also
satisfy {\rm (\ref{eq:constcurv})}; see 
\cite[Fig.\ 11]{gflows3d} and \S{\rm \ref{numexpt3d}} below.
\end{rem}

\subsection{Alternative curvature treatment} \label{sec:36}
As mentioned in Section~\ref{sec:2}, there is an alternative way to 
approximate the curvature vector $\varkappa\,\vec\nu$ 
in (\ref{eq:LBop}). In contrast to
the strategy employed in (\ref{eq:HGa}--d), where $\varkappa$ and $\vec\nu$ are
discretized separately, it is also possible to discretize 
$\vec\varkappa:= \varkappa\,\vec\nu$ directly, as proposed in the seminal paper
\cite{Dziuk91}. We then obtain the following linear finite element 
approximation.
For $m=0,\ldots, M-1$, find $\vec U^{m+1} \in \uspace^m$, 
$P^{m+1} \in \widehat\pspace^m$, $\vec{X}^{m+1}\in\Vh$ and 
$\vec\kappa^{m+1} \in \Vh$ such that
\begin{subequations}
\begin{align}
& 2\left(\mu^m\,\mat D(\vec U^{m+1}), \mat D(\vec \xi) \right)
- \left(P^{m+1}, \nabla\,.\,\vec \xi\right)
 - \gamma\,\left\langle \vec\kappa^{m+1} ,
   \vec\xi\right\rangle_{\Gamma^m}
\nonumber \\ & \hspace{4cm}
= \left(\vec f^{m+1}, \vec \xi\right)
\quad \forall\ \vec\xi \in \uspace^m \,, \label{eq:GDa}\\
& \left(\nabla\,.\,\vec U^{m+1}, \varphi\right) = 0 
\quad \forall\ \varphi \in \widehat\pspace^m\,,
\label{eq:GDb} \\
& \left\langle \frac{\vec X^{m+1} - \vec X^m}{\tau_m} ,
\vec\chi \right\rangle_{\Gamma^m}
- \left\langle \vec U^{m+1}, \vec\chi \right\rangle_{\Gamma^m}  = 0
\quad
 \forall\ \vec\chi \in \Vh\,,
\label{eq:GDc} \\
& \left\langle \vec\kappa^{m+1} , \vec\eta \right\rangle_{\Gamma^m}
+ \left\langle \nabs\,\vec X^{m+1}, \nabs\,\vec \eta \right\rangle_{\Gamma^m}
 = 0  
\quad 
\forall\ \vec\eta \in \Vh
\label{eq:GDd}
\end{align}
\end{subequations}
and set $\Gamma^{m+1} = \vec{X}^{m+1}(\Gamma^m)$. A 
discretization based on (\ref{eq:GDa}--d)
has first been proposed by B{\"a}nsch in \cite{Bansch01} 
for one-phase flow with a free capillary surface in the very
special situation that
\begin{equation} \label{eq:Bansch}
\vec\xi\!\mid_{\Gamma^m} \in \Vh \quad \forall\ \vec\xi\in \uspace^m\,.
\end{equation}
Clearly, (\ref{eq:Bansch}) requires the fitted approach and in that case can be
satisfied e.g.\ for the lowest order Taylor--Hood element, P2--P1, and a
piecewise quadratic variant of $\Vh$, see \cite{Bansch01}, 
or for the MINI element, P1$^{bubble}$--P1, 
with the piecewise linear $\Vh$. % from (\ref{eq:Vh}).
We note that if (\ref{eq:Bansch}) holds, 
then (\ref{eq:GDa}--d) can be equivalently rewritten as
\begin{align} \label{eq:Banscha}
& 2\left(\mu^m\,\mat D(\vec U^{m+1}), \mat D(\vec \xi) \right)
- \left(P^{m+1}, \nabla\,.\,\vec \xi\right)
 + \gamma\,\left\langle \nabs\,\vec X^{m+1}, \nabs\,
   \vec\xi\right\rangle_{\Gamma^m}
\nonumber \\ & \hspace{4cm}
= \left(\vec f^{m+1}, \vec \xi\right)
\quad \forall\ \vec\xi \in \uspace^m \,, 
\end{align}
together with (\ref{eq:GDb},c). For a
nonlinear variant of this scheme involving space-time finite elements,
in the context of the Navier--Stokes equations with a free capillary surface,
B{\"a}nsch proved existence, 
uniqueness and stability of discrete solutions, see \cite{Bansch01}. 

It is not difficult to prove existence, uniqueness and stability results also 
for the linear scheme (\ref{eq:GDa}--d) for two-phase Stokes flow. 
In particular,
one can show that there exists a unique solution to \mbox{(\ref{eq:GDa}--d)} 
that also satisfies the stability bounds (\ref{eq:stab})
and (\ref{eq:stabstab}). Similarly, the analogues of Lemma~\ref{lem:stat1}
and, if (\ref{eq:Bansch}) is satisfied, of Lemma~\ref{lem:stat2} hold. In the
latter case we observe that the curvature part of the discrete stationary
solution is given by the unique $\vec\kappa^{m+1} \in \Vh$ such that
$\langle \vec\kappa^{m+1}, \vec\eta \rangle_{\Gamma^m} =
\overline\kappa\,\langle \vec\nu^{m}, \vec\eta \rangle_{\Gamma^m}$ for all
$\vec\eta\in\Vh$.
We stress that if (\ref{eq:Bansch}) does not hold, then it does not appear 
possible to show the analogue of Lemma~\ref{lem:stat2}, which means that it 
is not possible to prove the existence of
discrete stationary solutions for (\ref{eq:GDa}--d).

However, the crucial difference between (\ref{eq:GDa}--d) and 
(\ref{eq:HGa}--d) is that in (\ref{eq:GDc})
the tangential velocity of the discrete interface is fixed by $\vec U^{m+1}$, 
and this has two consequences. Firstly, there is no
guarantee that the mesh quality of $\Gamma^m$ will be preserved. In fact, as
mentioned in the Introduction, typically the mesh will deteriorate over time.
And secondly, even for the case that $\charfcn{\Omega_-^m} \in \pspace^m$, 
it is not possible to prove (\ref{eq:consm}) for \mbox{(\ref{eq:GDa}--d)}, 
as $\vec\chi = \vec\nu^m$ is not a valid test function in (\ref{eq:GDc}),
and so true 
volume conservation in the semidiscrete setting, recall (\ref{eq:cons}), 
cannot be shown.
It is for these reasons that we prefer to use (\ref{eq:HGa}--d). 

We now return to the equivalent rewrite (\ref{eq:Banscha}), (\ref{eq:GDb},c)  
of (\ref{eq:GDa}--d). This has the advantage that the explicit computation of
the curvature vector $\vec\kappa^{m+1}$ is avoided, although in practice the
gain in computational efficiency is negligible because 
the main computational task is to solve the Stokes equations in the bulk.
We recall that for the equivalence between (\ref{eq:Banscha}), (\ref{eq:GDb},c) and the original (\ref{eq:GDa}--d) it was crucial to enforce the strong 
assumption (\ref{eq:Bansch}), which relies on a fitted bulk mesh. A variant of
this rewrite can be obtained for (\ref{eq:GDa}--d),
with $\langle \cdot, \cdot \rangle_{\Gamma^m}$ replaced by
$\langle \cdot, \cdot \rangle_{\Gamma^m}^h$, also in the absence of the
assumption (\ref{eq:Bansch}). This then leads to (\ref{eq:Banscha}) with
$\nabs\,\vec\xi$ replaced by $\nabs\,(\vec\pi^m\,\vec\xi)$. Being equivalent,
this rewrite inherits all the theoretical properties of (\ref{eq:GDa}--d),
with $\langle \cdot, \cdot \rangle_{\Gamma^m}$ replaced by
$\langle \cdot, \cdot \rangle_{\Gamma^m}^h$, namely the stability results
(\ref{eq:stab}) and (\ref{eq:stabstab}) and the analogue of
Lemma~\ref{lem:stat1}.
 
It is important to note that unconditional stability for (\ref{eq:Banscha}), 
(\ref{eq:GDb},c), even in the presence of (\ref{eq:Bansch}), 
can no longer be shown if the third term in (\ref{eq:Banscha}) is changed to
\begin{equation} \label{eq:alternatives}
\gamma\left\langle \nabs\,\vec X^{m}, \nabs\,\vec
\xi\right\rangle_{\Gamma^m}\,,
\end{equation}
where we note that $\langle \nabs\,\vec X^{m}, \nabs\,\vec \xi
\rangle_{\Gamma^m} = \langle 1, \nabs\,.\,\vec \xi\rangle_{\Gamma^m}$.
Proving existence and uniqueness for
this simpler variant is trivial, but it is no longer possible to establish
stability, or the analogues of Lemmas~\ref{lem:stat1} and \ref{lem:stat2}.
Alternatively, the term in (\ref{eq:alternatives}), in view of 
(\ref{eq:GDc}), can also be replaced by 
\begin{equation} \label{eq:GT08}
\gamma\left\langle \nabs\,\vec X^{m}, \nabs\,\vec
\xi\right\rangle_{\Gamma^m} + \gamma\,\tau_m\,
\left\langle \nabs\,\vec U^{m+1}, \nabs\,\vec \xi\right\rangle_{\Gamma^m}\,,
\end{equation}
which collapses to (\ref{eq:Banscha}) if the condition (\ref{eq:Bansch}) holds.
The formulation (\ref{eq:Banscha}), (\ref{eq:GDb},c) with the third term in 
(\ref{eq:Banscha}) replaced by (\ref{eq:GT08}) has been exploited in 
\cite{GanesanT08} for an axisymmetric formulation of two-phase Navier--Stokes 
flow.

Finally, we mention that the ideas presented in this section on how to
discretize the curvature term arising from (\ref{eq:2e}) can also be applied
in the level set approach. Here (\ref{eq:GDc}) is replaced with an 
approximation of the level set transport equation
\begin{equation} \label{eq:levelset}
\phi_t + \vec u \,.\,\nabla\,\phi = 0 \qquad\text{in $\Omega$}\,,
\end{equation}
and this is then combined with (\ref{eq:GDa},b), where the third term in 
(\ref{eq:GDa}) is replaced by (\ref{eq:alternatives})
with $\Gamma^m$ now being a suitable reconstruction of the discrete interface
arising from the zero level set of
a discretization of the level set function $\phi$ in (\ref{eq:levelset}). 
See e.g.\ \cite{GrossR07,GrossR11,AusasBI12} for some examples.
We stress that the level set method is a convenient computational tool for the
interface motion, but that it does not appear possible to derive 
stability results in the spirit of (\ref{eq:stab}) and (\ref{eq:stabstab}) 
for the level set method.

\subsection{Solution methods} \label{sec:37}
As is standard practice for the solution of linear systems arising from
discretizations of Stokes and Navier--Stokes equations, we avoid the
complications of the constrained pressure space $\widehat\pspace^m$ in practice
by considering an overdetermined linear system with $\pspace^m$ instead. 
With a view towards some numerical test cases in Section~\ref{sec:6},
we also allow for 
an inhomogeneous Dirichlet boundary condition $\vec g$ on $\partial\Omega$
and for ease of exposition
consider only piecewise quadratic velocity approximations.
Then we reformulate \mbox{(\ref{eq:HGa}--d)} as follows.
Find $\vec U^{m+1} \in \uspace^m(\vec g) :=
\{\vec U \in [S^m_2]^d : \vec U = \vec I^m_2\,\vec g\ 
\text{on}\ \partial\Omega\}$,
$P^{m+1} \in \pspace^m$, $\vec{X}^{m+1}\in\Vh$ and 
$\kappa^{m+1} \in \Wh$ such that (\ref{eq:HGa},c,d)
with $\uspace^m = [S^m_2]^d \cap \uspace$ hold together with
\begin{equation} \label{eq:LAb}
 \left(\nabla\,.\,\vec U^{m+1}, \varphi\right) = 
\frac{\left(\varphi, 1\right)}{\mathcal{L}^d(\Omega)}\, 
\int_{\partial\Omega} (\vec I^m_2\,\vec g) \,.\, \unitn \dH{d-1} 
\quad \forall\ \varphi \in \pspace^m\,.
\end{equation}
If $(\uspace^m, \pspace^m)$ satisfy the LBB condition (\ref{eq:LBB}), then
the existence and uniqueness proof for a solution to (\ref{eq:HGa},c,d), 
(\ref{eq:LAb}) is as before. In the absence of (\ref{eq:LBB}), the existence
and uniqueness of a solution to the reduced system that is analogous to 
(\ref{eq:HGa},c,d), with $\uspace^m$ replaced by $\uspace^m_0$,
hinges on the nonemptiness of the set
$\uspace^m_0(\vec g) := 
\{ \vec U \in \uspace^m(\vec g) : (\nabla\,.\,\vec U, \varphi) = 0 \quad
\forall\ \varphi \in \widehat\pspace^m \}$.
The linear system (\ref{eq:HGa},c,d), (\ref{eq:LAb}) 
can be solved with the help of a Schur complement
approach, which reduces the system to a standard saddle point problem
arising from discretizations of Stokes problems. 

\section{Numerical results}  \label{sec:6}

For details on the assembly of the linear system arising at each time step, 
as well as
details on the adaptive mesh refinement algorithm and the solution procedure,
we refer to the forthcoming article \cite{fluidfbp}. In particular, we recall
that the scheme in general uses an adaptive bulk mesh that has a fine mesh size
$h_f$ around $\Gamma^m$ and a coarse mesh size $h_c$ 
further away from it. 
The special case $h_f =
h_c$ leads to a uniform bulk mesh which will be sufficient for some of the
simple test problems considered in this section. For all the numerical results
presented in this paper no refinement or remeshing procedure was applied to the
discrete interface approximations $\Gamma^m$. 
We remark that we implemented our scheme with the help of 
the finite element toolbox ALBERTA, see \cite{Alberta}.

In order to test our finite element approximation (\ref{eq:HGa}--d), we
consider the trivial true solution of a stationary circle/sphere,
as it has been considered in e.g.\ \cite{GanesanMT07,GrossR07}.
In particular, $\Gamma(t) := \{ \vec z \in \R^d : |\vec z| = r(t)\}$, where
\begin{subequations}
\begin{equation} \label{eq:radialr}
r(t) = r(0)\,,
\end{equation}
together with 
\begin{equation} \label{eq:radialup}
\vec u(\vec z, t) = \vec 0 \,,\quad
p(\vec z, t) = \lambda(t)\left[\charfcn{\Omega_-(0)}
-\frac{\mathcal{L}^d(\Omega_-(0))}{\mathcal{L}^d(\Omega)} \right] ,
\end{equation}
\end{subequations}
where $\lambda(t) = \lambda(0) = \gamma\,(d-1)\,[r(0)]^{-1}$, 
is an exact solution to the problem (\ref{eq:2a}--g).

A nontrivial divergence free and radially symmetric solution $\vec u$ can be
constructed on a domain that does not contain the origin. To this end, consider
e.g.\
$\Omega = (-H,H)^d \setminus [-H_0, H_0]^d$, with $0 < H_0 < H$. Then
$\Gamma(t) := \{ \vec z \in \R^d : |\vec z| = r(t)\}$, where
\begin{subequations}
\begin{equation} \label{eq:radialr2}
r(t) = ([r(0)]^d + \alpha\,t\,d)^\frac1d \,,
\end{equation}
together with 
\begin{equation} \label{eq:radialup2}
\vec u(\vec z, t) = \alpha\,|\vec z|^{-d}\,\vec z \,, \quad
p(\vec z, t) = \lambda(t)\left[ \charfcn{\Omega_-(t)} - 
\frac{\mathcal{L}^d(\Omega_-(t))}{\mathcal{L}^d(\Omega)}\right],
\end{equation}
\end{subequations}
where $\lambda(t) = \gamma\,(d-1)\,[r(t)]^{-1} + 
2\,\alpha\,(d-1)\,(\mu_+ - \mu_-)\,[r(t)]^{-d}$,
is an exact solution to the problem (\ref{eq:2a}--g) with the homogeneous right
hand side in (\ref{eq:2c}) replaced by $\vec g$, where 
$\vec g(\vec z) = \alpha\,|\vec z|^{-d}\,\vec z$.

{From} now on we fix $\Gamma(0) = \{ \vec z \in \R^d : |\vec z| = \frac12 \}$.
Throughout this section we use uniform time steps
$\tau_m=\tau$, $m=0,\ldots, M-1$.
For later use, we define $h^m_\Gamma := \max_{j = 1 ,\ldots, J^m_\Gamma} 
\diam(\sigma^m_j)$.
We also define the errors 
$$\errorXx :=
\max_{m=1,\ldots, M} \|\vec{X}^m - \vec{x}(\cdot,t_m)\|_{L^\infty}\,,$$ 
where $\|\vec{X}(t_m) - \vec{x}(\cdot,t_m)\|_{L^\infty}
:=\max_{k=1,\ldots, K^m_\Gamma}
\left\{\min_{\vec{y}\in \Upsilon} |\vec{X}^m(\vec{q}^m_k) - \right.$ $\left.
\vec{x}(\vec{y},t_m)|\right\}$ and
$$\errorUu2 :=
\max_{m=1,\ldots, M}\|U^m - \vec I^m_2\,u(\cdot,t_m)\|_{L^\infty(\Omega)}\,.$$
In order to evaluate the errors in the pressure, we define
$$
\LerrorPp := \left[\tau\,\sum_{m=1}^M \|P^m - p(\cdot,t_m)\|_{L^2(\Omega)}^2 
\right]^\frac12\,.
$$
When we use \XFEMGAMMA, we also evaluate the following errors for the pressure, $\LerrorPpc$
and $\errorLl := \max_{m=1,\ldots, M} |\lambda^m - \lambda(t_m)|$. 
Here
$p_c(\cdot,t_m):= p(\cdot,t_m)-\lambda(t_m)\,\charfcn{\Omega_-(t_m)} \in \R$
for the test problems (\ref{eq:radialr},b) and (\ref{eq:radialr2},b),
and $P^{m}_c := P^{m} - \lambda^{m}\,\charfcn{\Omega^{m-1}_-}$
is piecewise polynomial on $\mathcal{T}^{m-1}$.

\subsection{Numerical results in 2d} \label{sec:61}
For our first set of experiments we fix $\Omega = (-1,1)^2$ and use the true
solution (\ref{eq:radialr},b) for the parameters
$$
\mu = \gamma = 1\,.
$$
This means that the true solution reduces to
$r(t) = \frac12$, $\vec u(\cdot, t) = \vec 0$ and $p(t) =
2\,[\charfcn{\Omega_-(0)} - \tfrac14\,\mathcal{L}^2(\Omega_-(0)) ]$ 
for all $t\geq0$.
Some errors for our approximation (\ref{eq:HGa}--d) for the P2--P1 element 
can be seen in Table~\ref{tab:sb0p2p1}. The same convergence
test for the pressure spaces P0 and P1+P0 
are shown in Tables~\ref{tab:sb0p2p0} and \ref{tab:sb0p2p10}, respectively. 
Here we always choose uniform spatial discretizations such that 
$h_c = h_f = h$ and $h^m_\Gamma \approx h / 8$. We varied the time
discretization parameter $\tau$ from $10^{-2}$ to $10^{-4}$, but as the errors
were nearly indistinguishable for these runs, we only report on the
computations with $\tau = 10^{-2}$.
It appears that the three errors $\errorXx$, $\errorUu2$ and $\LerrorPp$
in Tables~\ref{tab:sb0p2p1}--\ref{tab:sb0p2p10} converge with
$\mathcal{O}(h)$, $\mathcal{O}(h)$ and $\mathcal{O}(h^\frac12)$, respectively.
The latter is to be expected for the approximation of the discontinuous
pressure on uniform meshes.
\begin{table*}
\center
\begin{tabular}{lllll}
\hline
 $2^\frac12/h$ & $\tau$ & $\errorXx$ & $\errorUu2$ & $\LerrorPp$ \\
 \hline % the following numbers are for N_f = N_c
   4  & $10^{-2}$ & 1.7401e-02 & 3.4406e-02 & 5.8656e-01 \\ % & 6.6637e-01 \\
   8  & $10^{-2}$ & 7.9853e-03 & 1.7896e-02 & 4.0791e-01 \\ % & 2.7465e-00 \\
   16 & $10^{-2}$ & 3.5541e-03 & 8.9120e-03 & 2.9411e-01 \\ % & 3.2791e-00 \\
\hline
\end{tabular}
\caption{($\mu=\gamma=1$)
Stationary bubble problem on $(-1,1)^2$ over the time interval $[0,1]$
for the P2--P1 element without \XFEMGAMMA.}
\label{tab:sb0p2p1}
\end{table*}%
\begin{table*}
\center
\begin{tabular}{lllll}
\hline
 $2^\frac12/h$ & $\tau$ & $\errorXx$ & $\errorUu2$ & $\LerrorPp$ \\
 \hline % the following numbers are for N_f = N_c
   4  & $10^{-2}$ & 2.0198e-02 & 3.0165e-02 & 5.9176e-01 \\ % & 5.6213e-01 \\
   8  & $10^{-2}$ & 9.7242e-03 & 1.4778e-02 & 4.5001e-01 \\ % & 1.4163e-00 \\
   16 & $10^{-2}$ & 4.6101e-03 & 8.0770e-03 & 3.1948e-01 \\ % & 3.2247e-00 \\
\hline
\end{tabular}
\caption{($\mu=\gamma=1$)
Stationary bubble problem on $(-1,1)^2$ over the time interval $[0,1]$
for the P2--P0 element without \XFEMGAMMA.}
\label{tab:sb0p2p0}
\end{table*}%
\begin{table*}
\center
\begin{tabular}{lllll}
\hline
 $2^\frac12/h$ & $\tau$ & $\errorXx$ & $\errorUu2$ & $\LerrorPp$ \\
 \hline % the following numbers are for N_f = N_c
   4  & $10^{-2}$ & 5.8351e-03 & 1.1813e-02 & 4.4080e-01 \\ % & 4.5340e-01 \\
   8  & $10^{-2}$ & 2.1014e-03 & 5.6510e-03 & 3.2709e-01 \\ % & 1.0941e-00 \\
   16 & $10^{-2}$ & 5.9531e-04 & 3.3472e-03 & 2.3255e-01 \\ % & 2.9170e-00 \\
\hline
\end{tabular}
\caption{($\mu=\gamma=1$)
Stationary bubble problem on $(-1,1)^2$ over the time interval $[0,1]$
for the P2--(P1+P0) element without \XFEMGAMMA.}
\label{tab:sb0p2p10}
\end{table*}%

\begin{figure*}
\center
\mbox{
\includegraphics[angle=-0,width=0.33\textwidth]{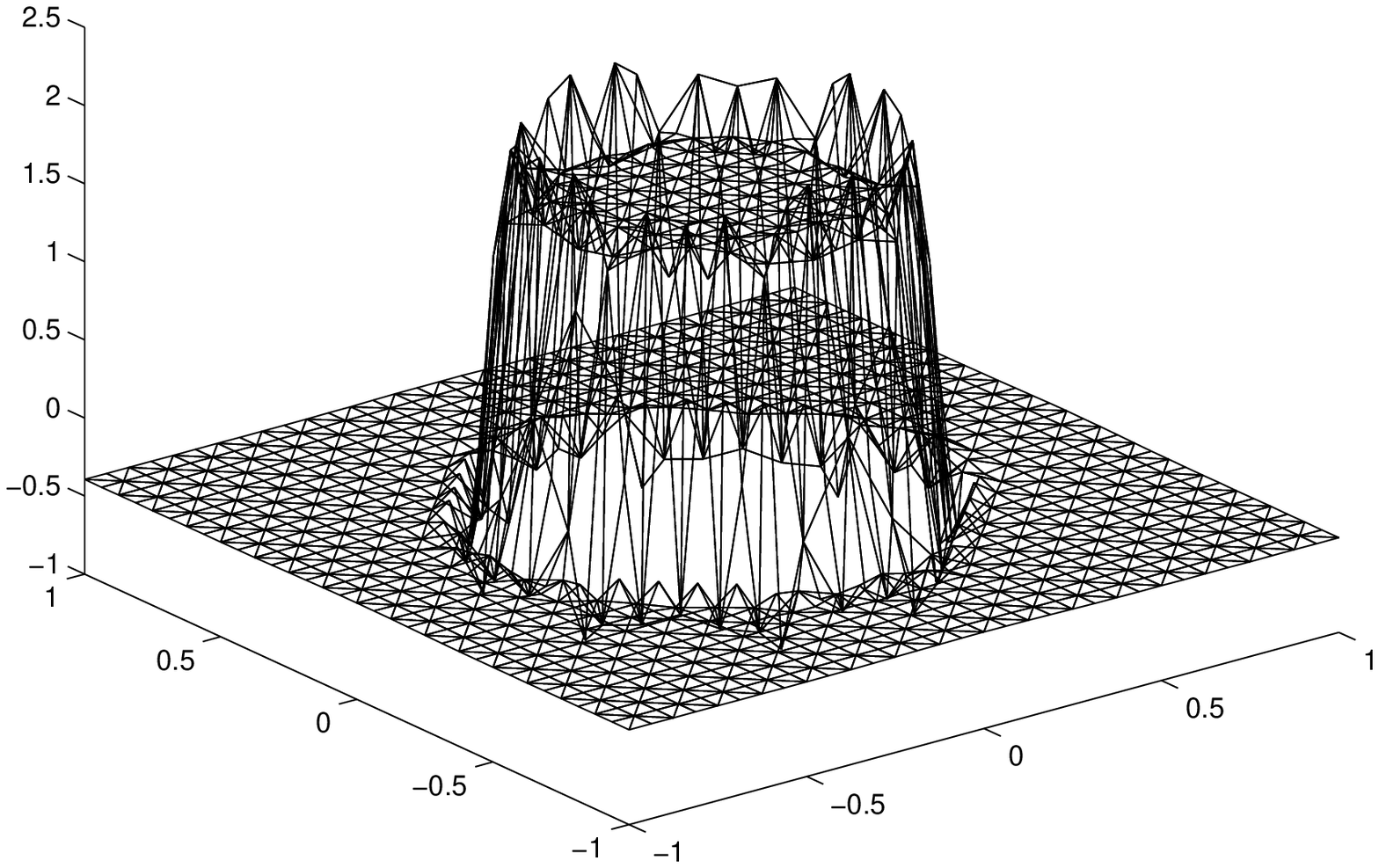}
\includegraphics[angle=-0,width=0.33\textwidth]{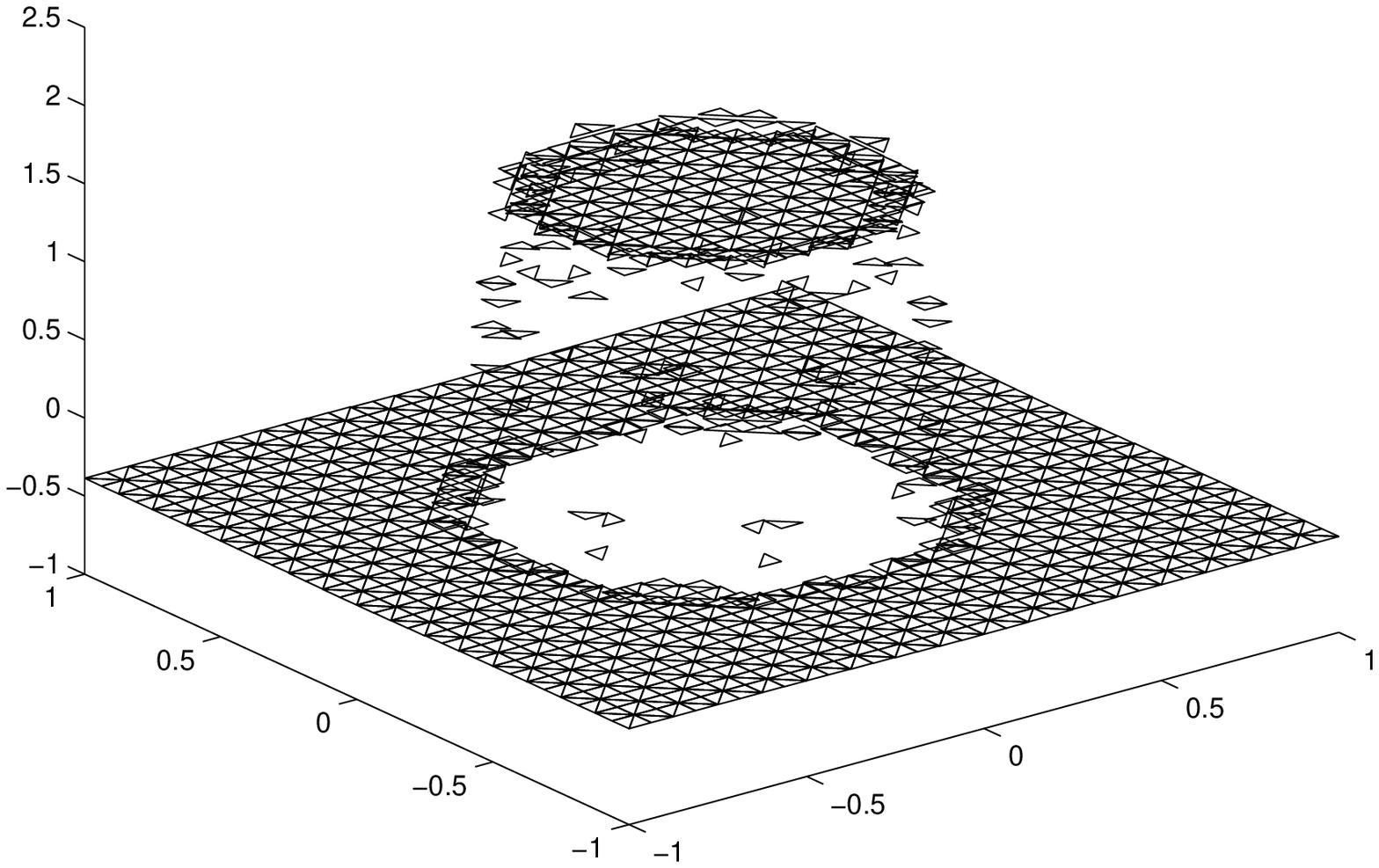}
\includegraphics[angle=-0,width=0.33\textwidth]{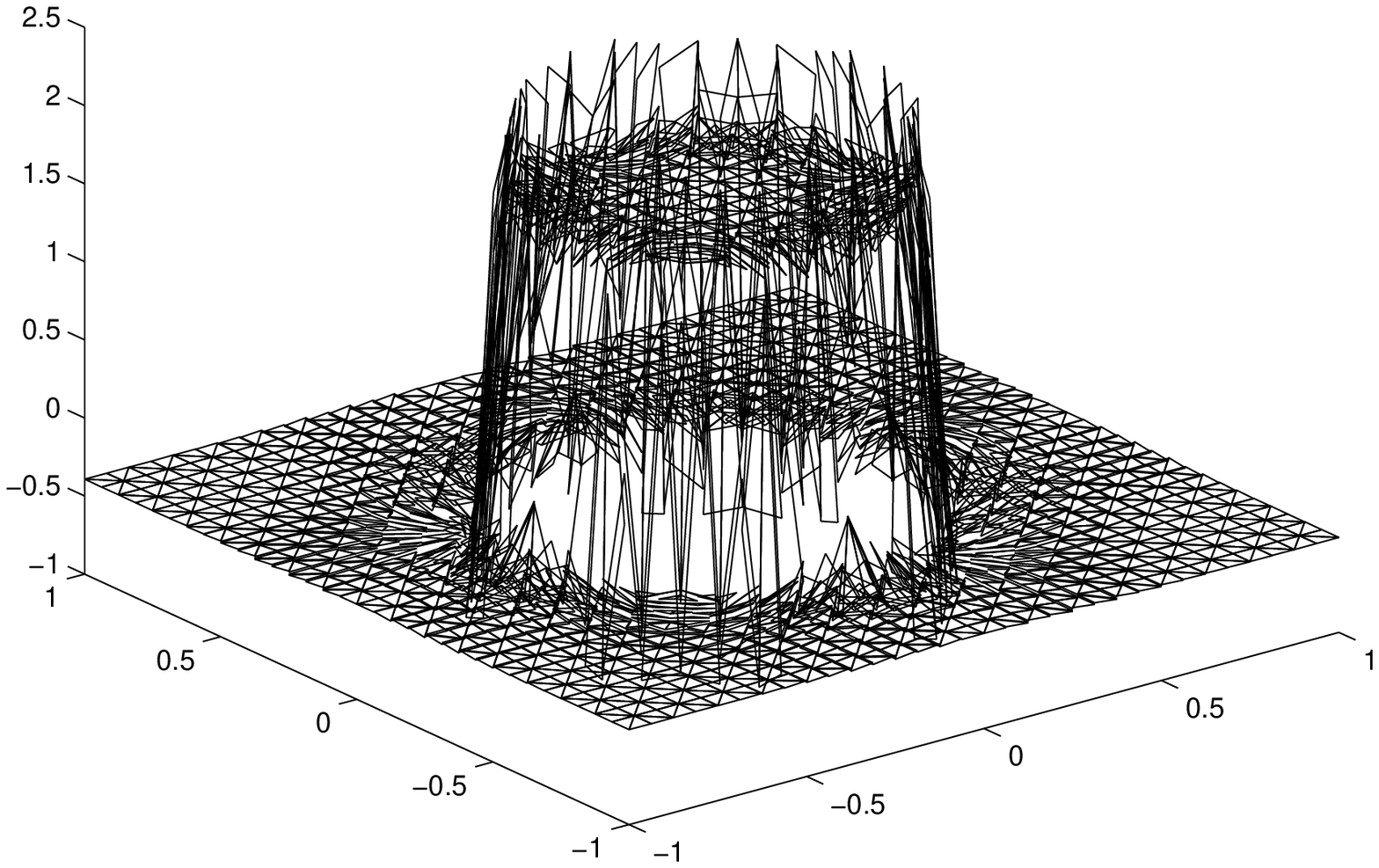}}
\mbox{
\includegraphics[angle=-0,width=0.33\textwidth]{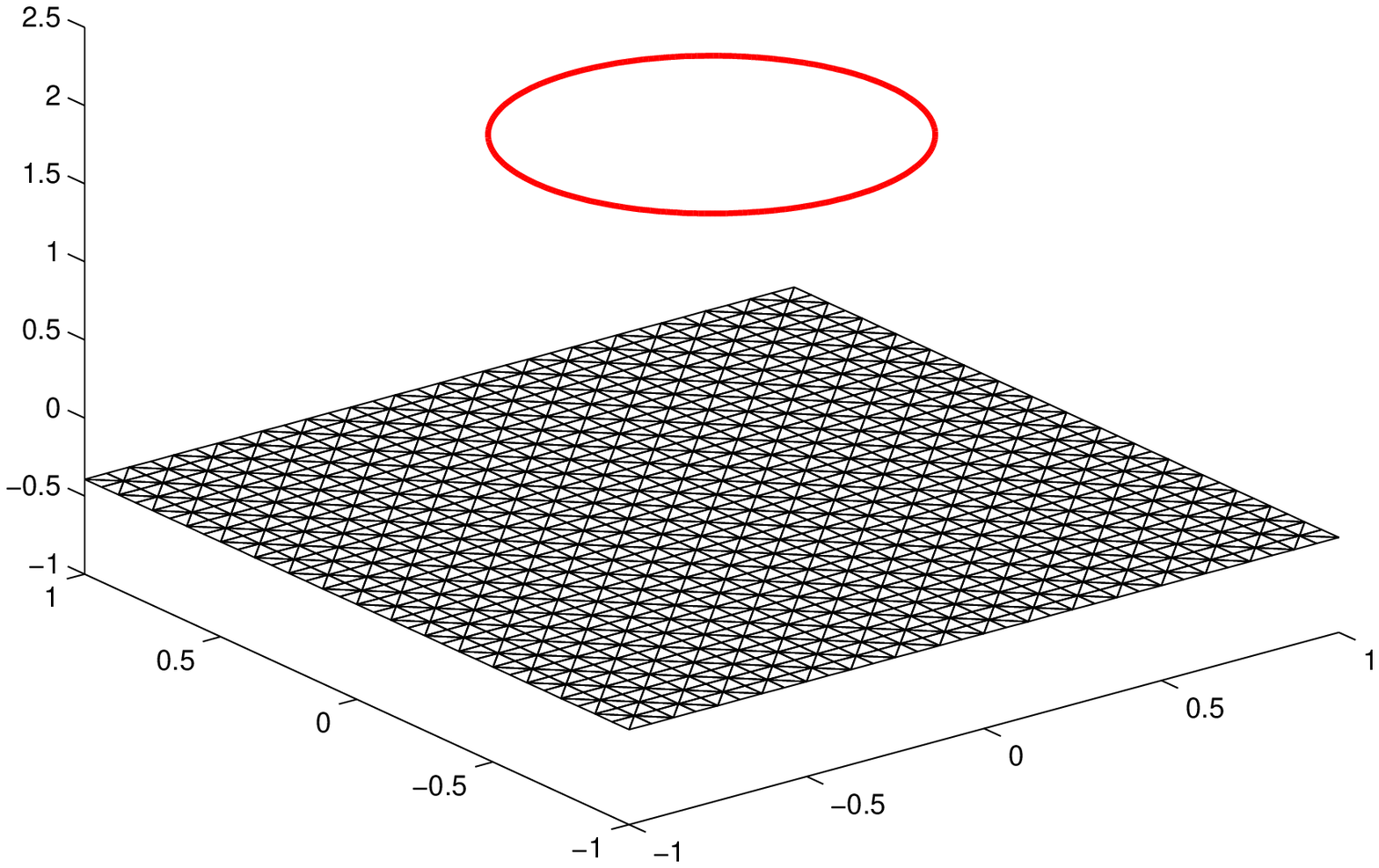}
\includegraphics[angle=-0,width=0.33\textwidth]{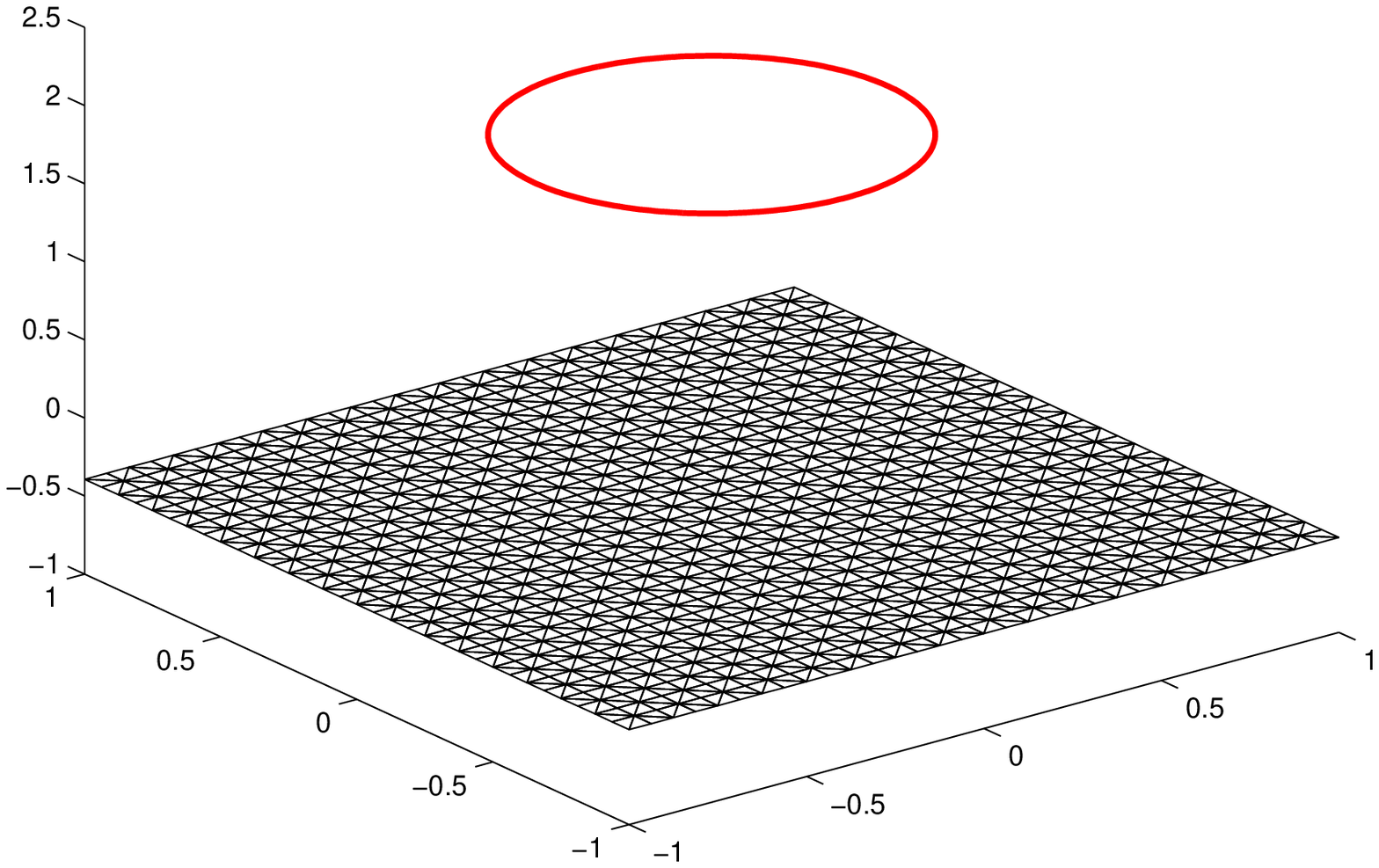}
\includegraphics[angle=-0,width=0.33\textwidth]{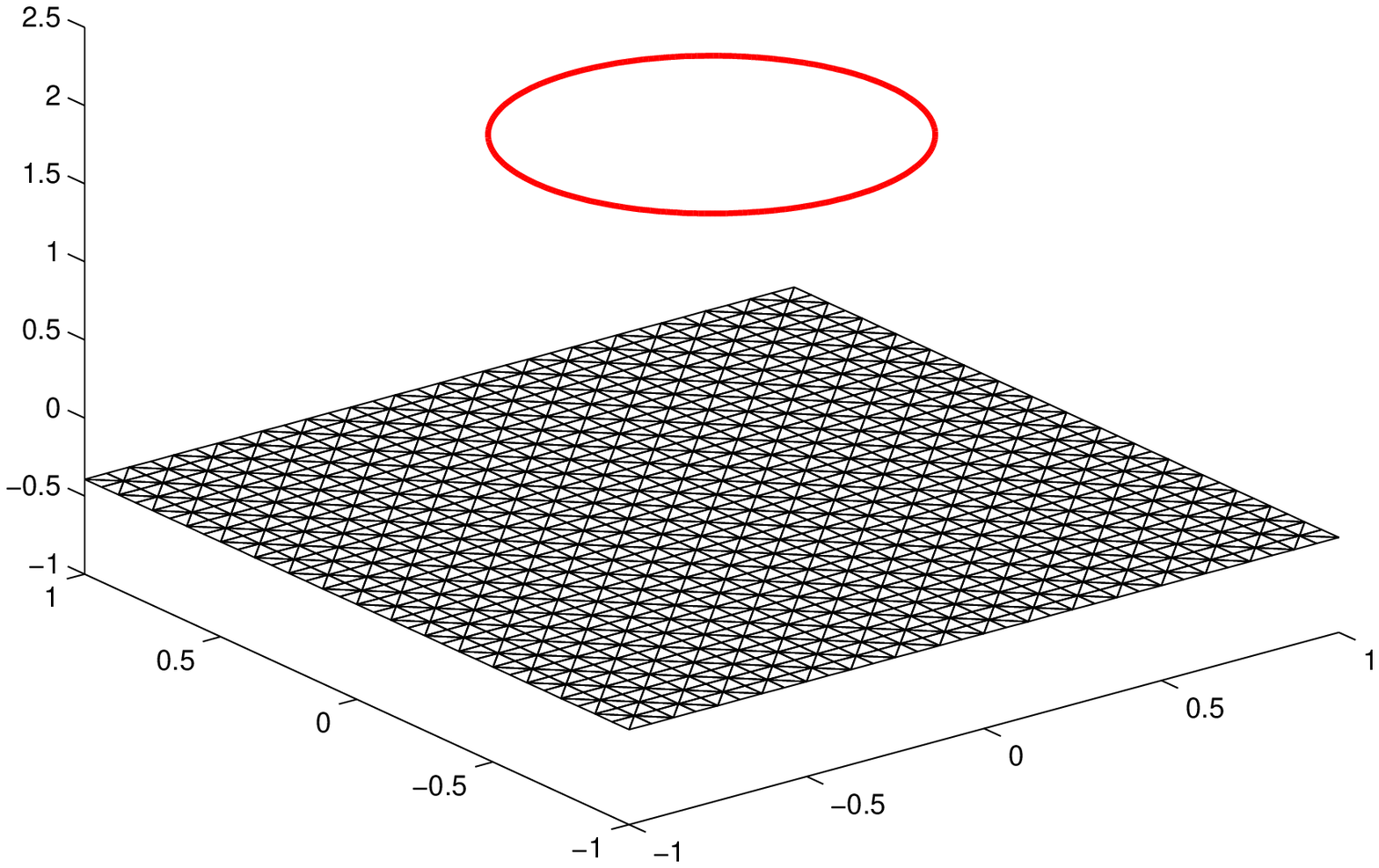}}
\caption{($\mu=\gamma=1$)
Pressure plots at time $T=1$ for the stationary bubble problem. The
pressure spaces are P1, P0 and P1+P0 without and with \XFEMGAMMA.}
\label{fig:pressures}
\end{figure*}

We note that for the experiments in Tables~\ref{tab:sb0p2p1}--\ref{tab:sb0p2p10}
we choose an equidistributed approximation $\Gamma^0$ of the circle 
$\Gamma(0)$. In this special case our approximation with \XFEMGAMMA, i.e.\
(\ref{eq:HGa}--d) with $\charfcn{\Omega^m_-} \in \pspace^m$,
yields the exact solution $\vec U^{m+1} = \vec 0$. 
In particular, 
on recalling Lemma~\ref{lem:stat2}, we have that the unique solution to
the reduced system (\ref{eq:HGa},c,d),
with $\uspace^m$ replaced by $\uspace^m_0$,
for $m=0$ is given by $\vec U^{m+1} = \vec 0$, 
$\vec X^{m+1} = \vec X^0$ and $\kappa^{m+1} = - \lambda^0 \in \R$, where
$- \lambda^0 \approx - 2$ approximates the curvature of $\Gamma(0)$, 
and by induction for all $m=0,\ldots, M-1$. This implies
that $\vec U^{m+1} = \vec 0$, 
$P^{m+1} = \lambda^0\,[\charfcn{\Omega_-^0} 
- \tfrac14\,\mathcal{L}^2(\Omega_-^0) ]$,  
$\vec X^{m+1} = \vec X^0$ and $\kappa^{m+1} = -\lambda^0$ is a solution
to (\ref{eq:HGa}--d), % with $\widehat\pspace^m$ replaced by 
and it is this solution that is
found by our solution method in practice, see 
Table~\ref{tab:sb0p2p1XFEM}. %--\ref{tab:sb0p2p10XFEM}.
We remark that these results remain unchanged for nonconstant $\mu$, e.g.\
when choosing $\mu_+ = 10\,\mu_- = 1$.
\begin{table*}
\center
\begin{tabular}{llccll}
\hline
 $2^\frac12/h$ & $\tau$ & $\errorXx$ & $\errorUu2$ & $\LerrorPpc$ 
 & $\errorLl$ \\ % & $\errorKk$ \\ 
 \hline % the following numbers are for N_f = N_c
   4  & $10^{-2}$ & 0 & 0 & 3.1537e-04 & 2.4120e-03 \\ % & 2.4120e-03\\
   8  & $10^{-2}$ & 0 & 0 & 7.8851e-05 & 6.0254e-04 \\ % & 6.0255e-04\\
   16 & $10^{-2}$ & 0 & 0 & 1.9713e-05 & 1.5061e-04 \\ % & 1.5061e-04\\
\hline
\end{tabular}
\caption{($\mu=\gamma=1$)
Stationary bubble problem on $(-1,1)^2$ over the time interval $[0,1]$
for the P2--P1 element with \XFEMGAMMA. The same numbers are obtained for
the  P2--P0 and the P2--(P1+P0) element, respectively.}
\label{tab:sb0p2p1XFEM}
\end{table*}%
We visualize the final pressures for the finest runs in 
Tables~\ref{tab:sb0p2p1}--\ref{tab:sb0p2p1XFEM} in Figure~\ref{fig:pressures}.
Here in the case of the enrichment \XFEMGAMMA\ being used, we plot
$P^M_c$, which is almost identically equal to a constant, 
and $\lambda^M\,\charfcn{\Omega_-^{M-1}}$ separately.

\begin{figure*}
\center
\mbox{
\includegraphics[angle=-0,width=0.33\textwidth]{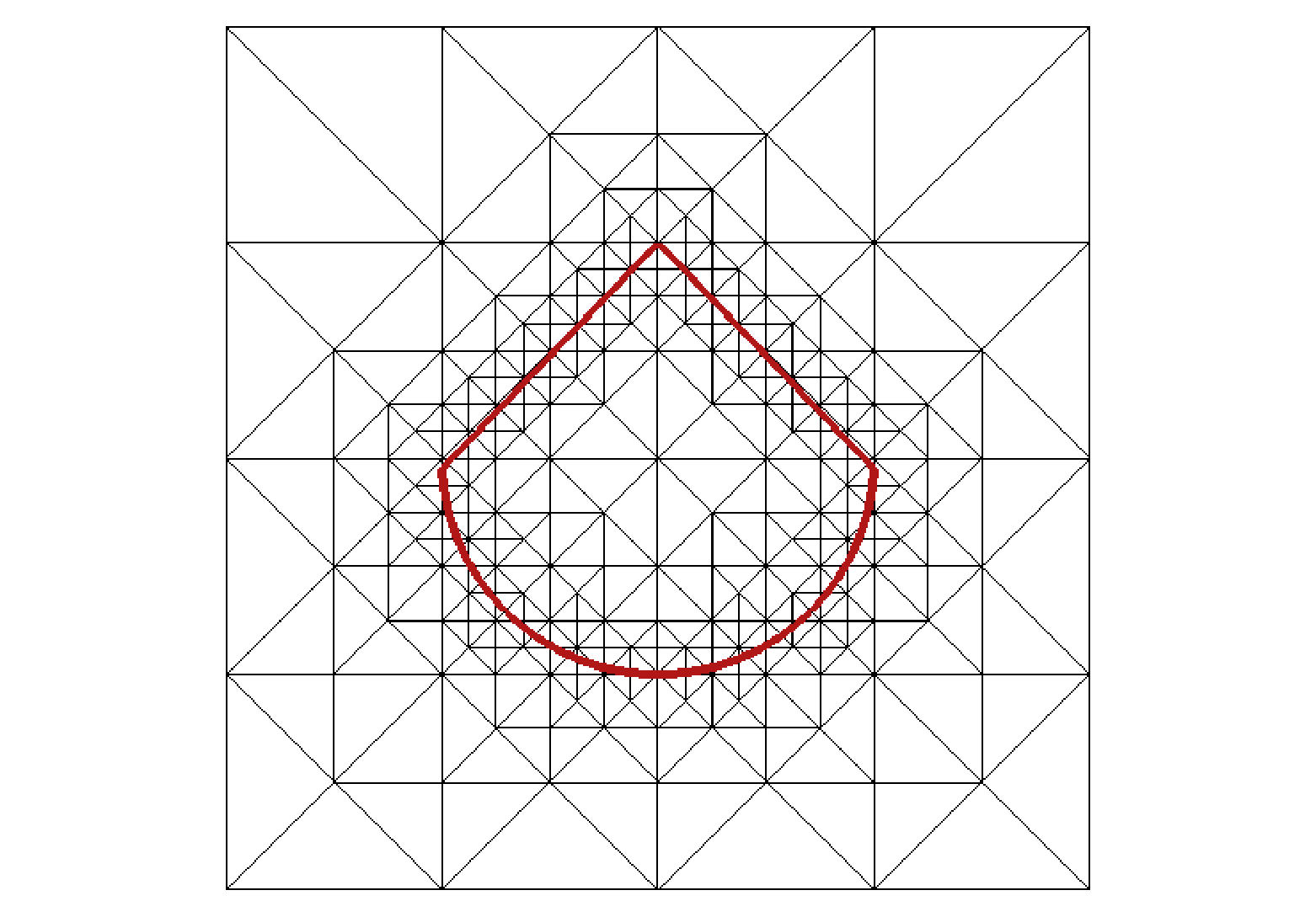}
\includegraphics[angle=-0,width=0.33\textwidth]{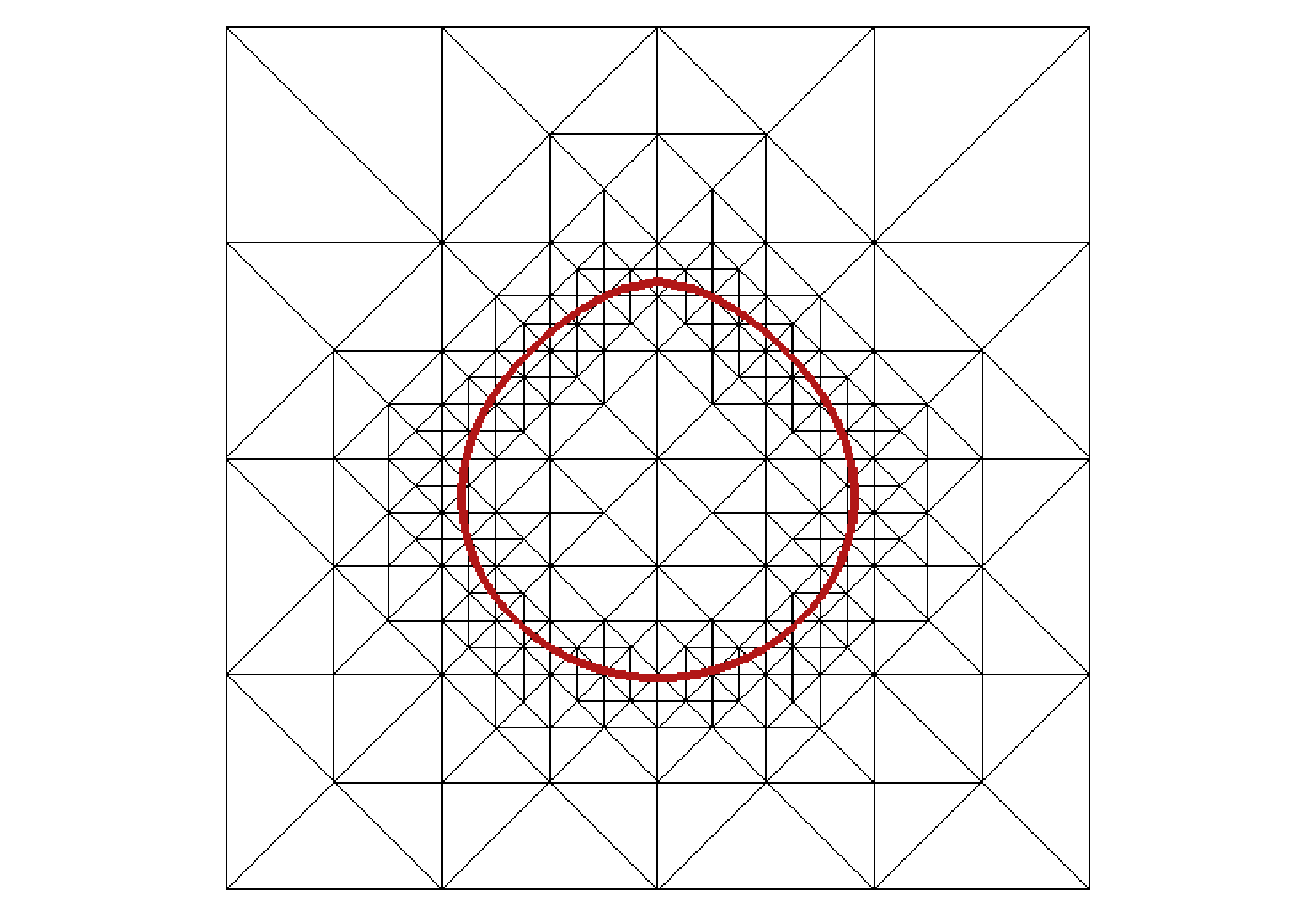}
\includegraphics[angle=-0,width=0.33\textwidth]{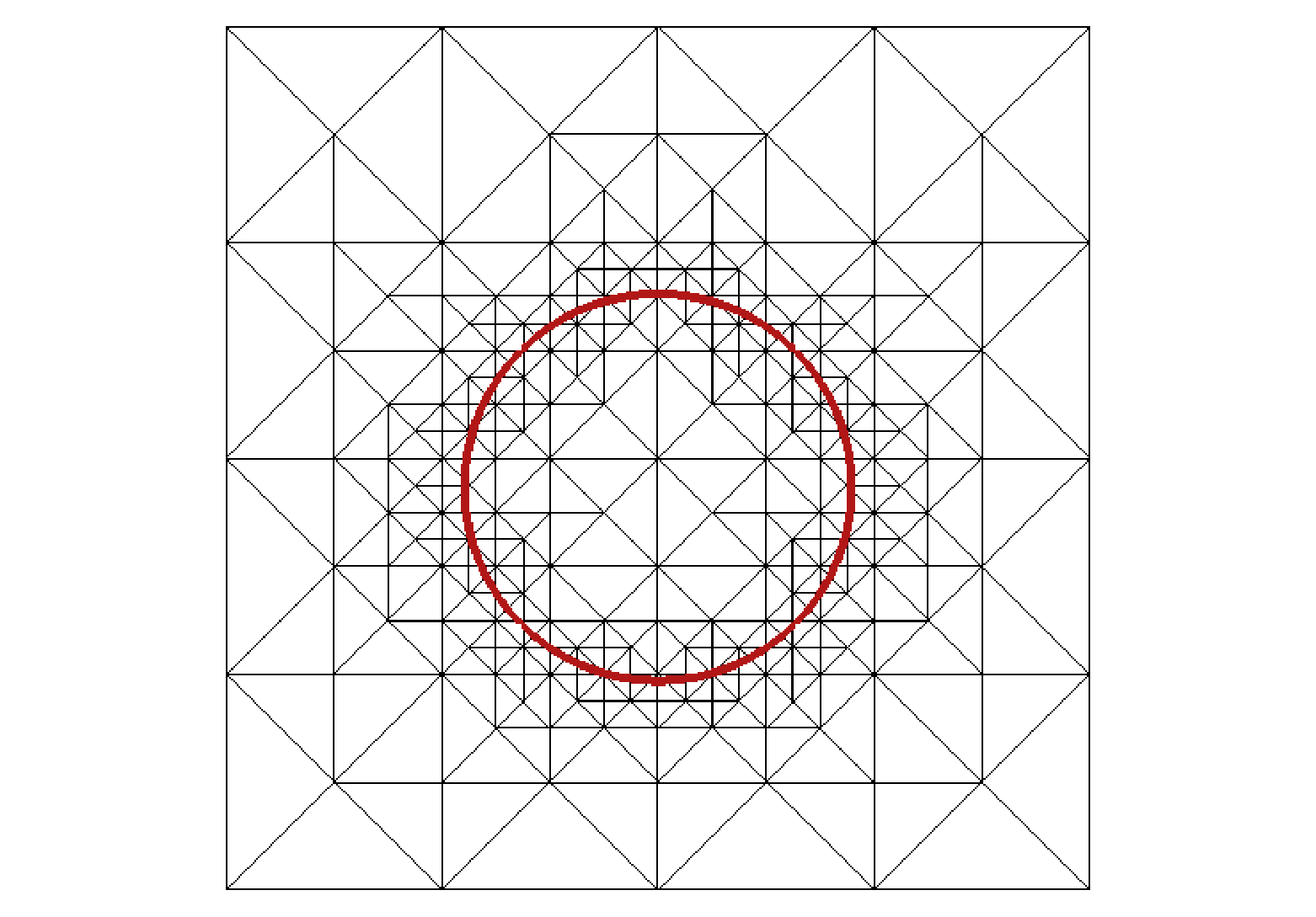}}
\mbox{
\includegraphics[angle=-90,width=0.33\textwidth]{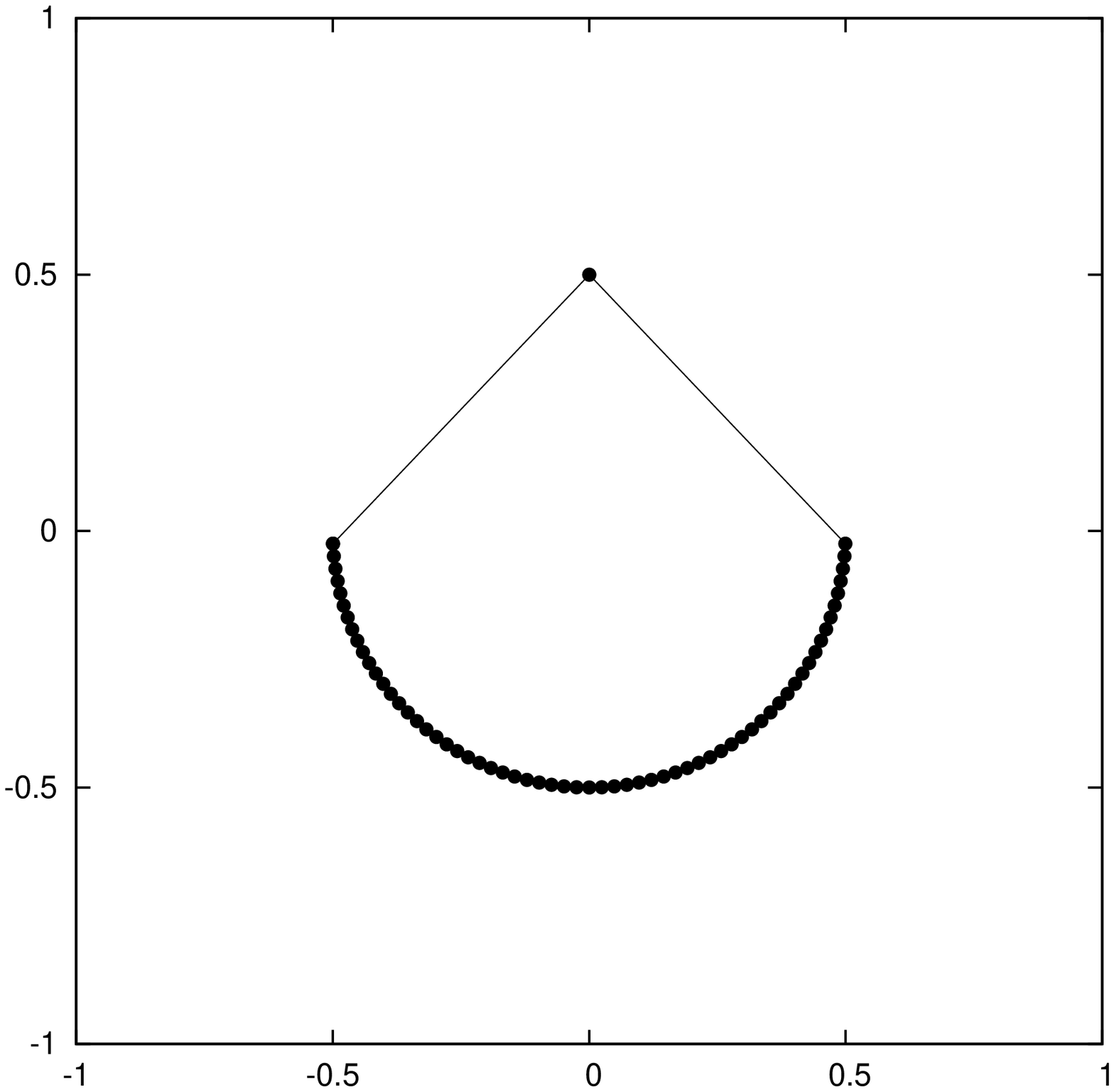}
\includegraphics[angle=-90,width=0.33\textwidth]{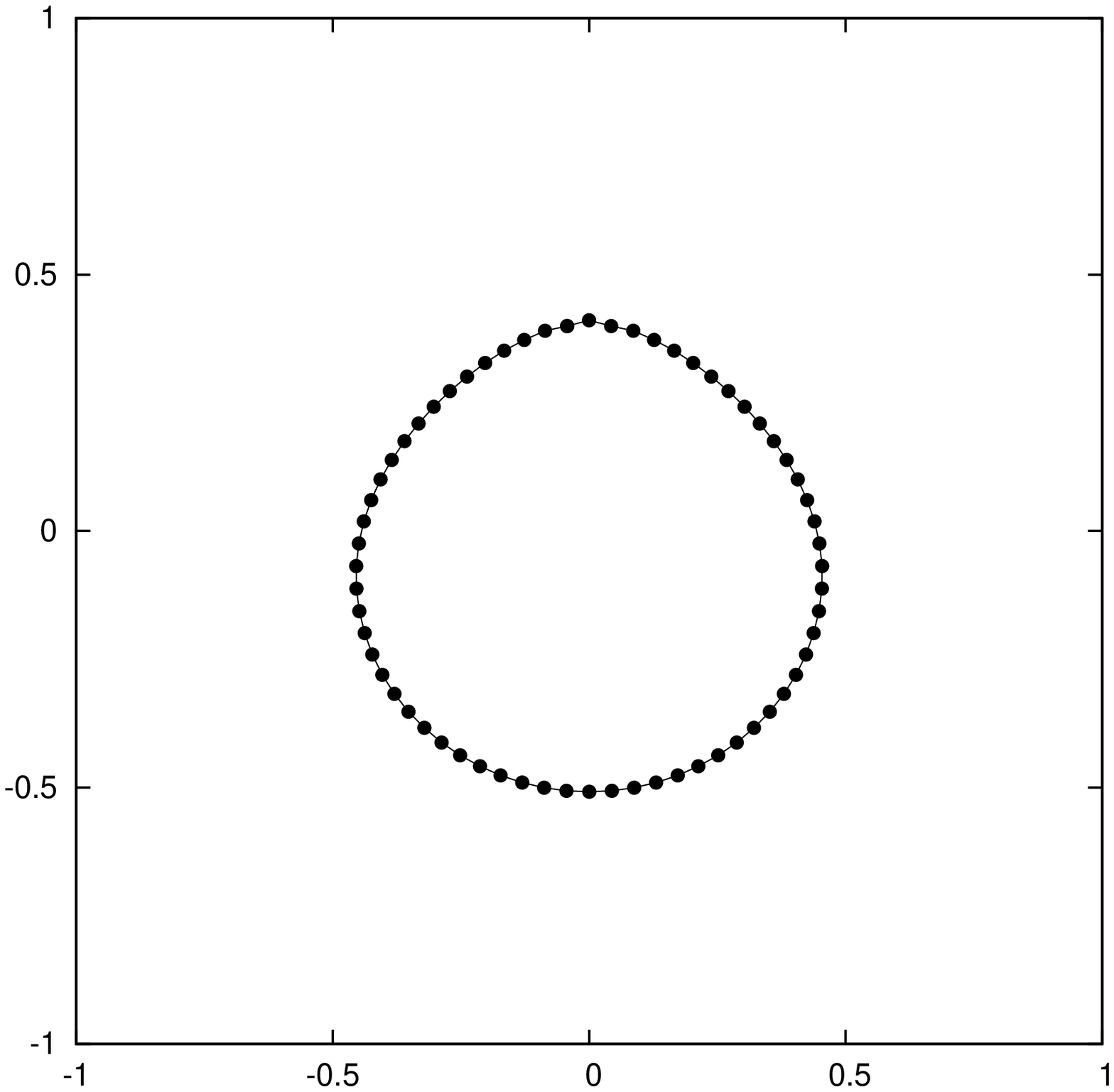}
\includegraphics[angle=-90,width=0.33\textwidth]{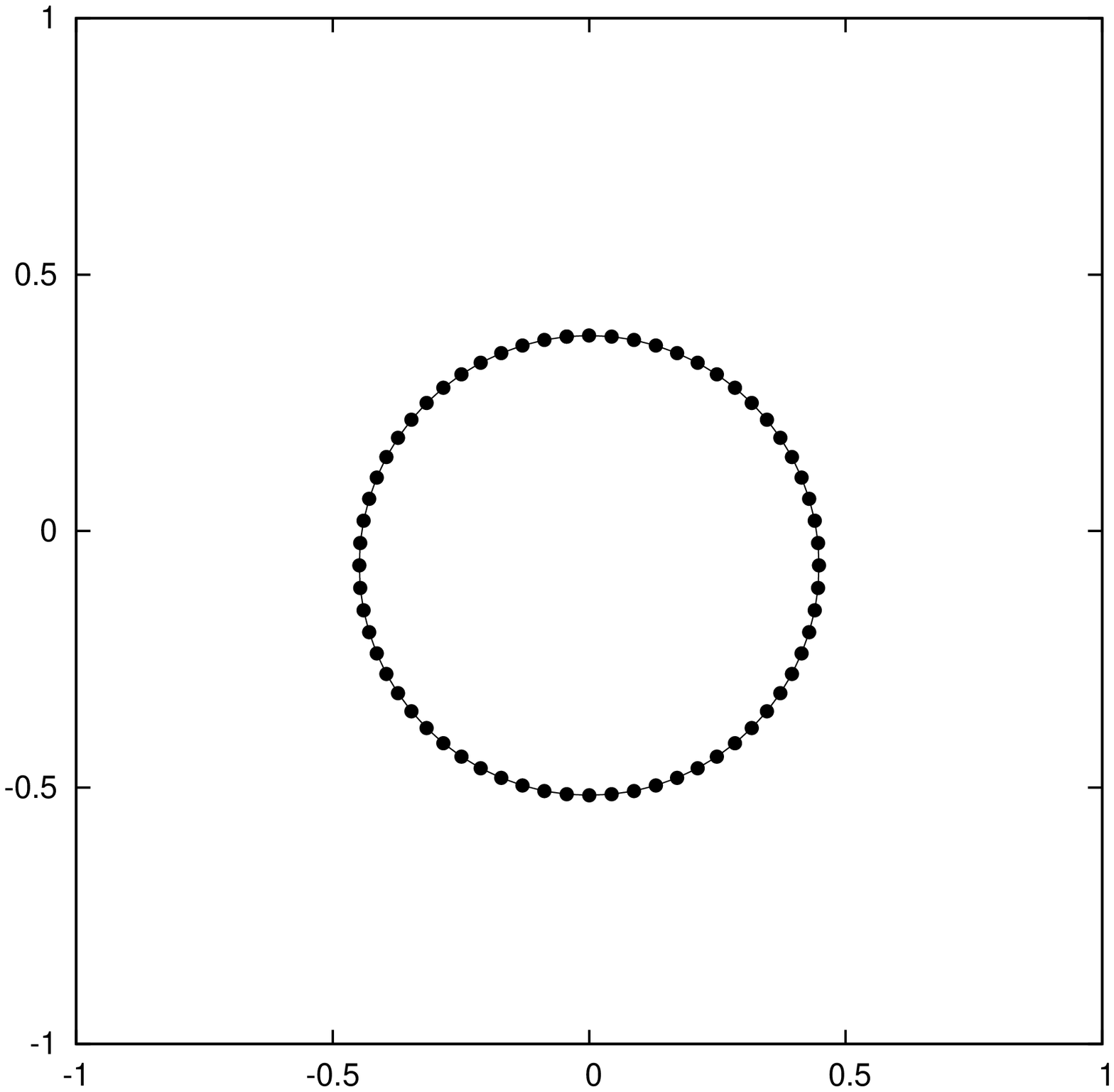}}
\caption{($\mu=\gamma=1$)
The discrete interface $\Gamma^m$ at times $t=0,\,1,\,5$ together with the adaptive bulk mesh (top), and details of the distribution of vertices on 
$\Gamma^m$ (bottom). Here we use the P2--P1 element with \XFEMGAMMA.}
\label{fig:usp}
\end{figure*}

\begin{figure}
\center
\includegraphics[angle=-90,width=0.4\textwidth]{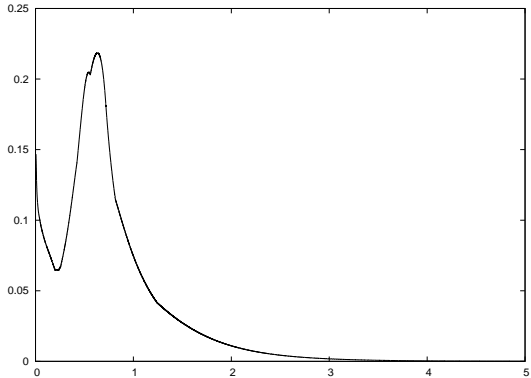}
\caption{($\mu=\gamma=1$)
Plot of $\|\vec U^m\|_{L^\infty(\Omega)}$
over the time interval $[0,5]$ for the P2--P1 element with 
\XFEMGAMMA.}
\label{fig:uspplot}
\end{figure}

We now demonstrate that this remarkable property is generic to our method, in
the sense that the circular, equidistributed numerical steady state solution 
is recovered by our method even if we choose very noncircular or very
nonuniform initial data $\Gamma^0$. Of course, this is the discrete analogue of
the fact that circles are the unique steady state solutions in the continuous
case, recall Section~\ref{sec:1}.
In particular, we choose $\Gamma^0$ to be a very nonuniform
approximation of $\Gamma(0)$, where we represent the upper half of the
circle by a single vertex, while the lower half is properly resolved to
resemble a semicircle. In total we use $K^0_\Gamma = 64$ vertices for 
$\Gamma^0$, and we use an adaptive bulk mesh with 
$h_c = 8\,h_f = 2^{-\frac12}$.
Choosing $\tau=10^{-4}$ we simulate the evolution with 
our scheme for the time interval $[0,5]$. In Figure~\ref{fig:usp} we show some
snapshots of the evolution, while in Figure~\ref{fig:uspplot} a plot of
$\|\vec U^m\|_{L^\infty(\Omega)}$ over time can be seen. 
Here %, and for the remainder of this paper, 
we use the P2--P1 element. As expected, the approximations $\Gamma^m$ converge
towards an equidistributed circle, while $\vec U^m$ converges to zero.

For our second set of convergence experiments we fix 
$\Omega = (-1,1)^2 \setminus [-\frac13,\frac13]^2$ and use the parameters
$$
\alpha = 0.15 \quad\text{and}\quad \mu = \gamma = 1
$$
for the true solution (\ref{eq:radialr2},b). 
With $T=1$ as before we obtain that $\Gamma(T)$ is a circle of radius
$r(1) = \sqrt{0.55} \approx 0.742$. Some errors for our approximation
(\ref{eq:HGa}--d), where we again use uniform bulk meshes with $h_c = h_f = h$
and $h^m_\Gamma \approx h/8$, 
are shown in Tables~\ref{tab:dn0p2p1}--\ref{tab:dn0p2p10XFEM}. We observe that
the computations without \XFEMGAMMA\ indicate a convergence for the pressure
error $\LerrorPp$ of $\mathcal{O}(h^\frac12)$. The simulations with
\XFEMGAMMA\ show a significant improvement in all errors compared to the
results in Tables~\ref{tab:dn0p2p1}--\ref{tab:dn0p2p10}.

\begin{table*}
\center
\begin{tabular}{lllll}
\hline
 $1/h$ & $\tau$ & $\errorXx$ & $\errorUu2$ & $\LerrorPp$ \\
 \hline % the following numbers are for N_f = N_c
   3  & $10^{-2}$ & 9.6884e-03 & 2.5200e-02 & 5.0519e-01 \\ % & 7.1419e-01 \\  
   6  & $10^{-3}$ & 3.4188e-03 & 1.6696e-02 & 3.5383e-01 \\ % & 4.7753e-01 \\
   12 & $10^{-4}$ & 1.7263e-03 & 7.7395e-03 & 2.4944e-01 \\ % & 6.5413e-01 \\
\hline
\end{tabular}
\caption{($\alpha = 0.15$, $\mu=\gamma=1$)
Expanding bubble problem on $(-1,1)^2 \setminus [-\frac13,\frac13]^2$ 
over the time interval $[0,1]$ 
for the P2--P1 element without \XFEMGAMMA.}
\label{tab:dn0p2p1}
\end{table*}%
\begin{table*}
\center
\begin{tabular}{lllll}
\hline
 $1/h$ & $\tau$ & $\errorXx$ & $\errorUu2$ & $\LerrorPp$ \\
 \hline % the following numbers are for N_f = N_c
   3  & $10^{-2}$ & 1.2003e-02 & 3.7696e-02 & 5.6246e-01 \\ % & 5.9399e-01 \\
   6  & $10^{-3}$ & 5.3537e-03 & 1.7035e-02 & 3.9392e-01 \\ % & 5.9340e-01 \\
   12 & $10^{-4}$ & 2.8129e-03 & 9.7768e-03 & 2.8218e-01 \\ % & 5.5011e-01 \\
\hline
\end{tabular}
\caption{($\alpha = 0.15$, $\mu=\gamma=1$)
Expanding bubble problem on $(-1,1)^2 \setminus [-\frac13,\frac13]^2$ 
over the time interval $[0,1]$ 
for the P2--P0 element without \XFEMGAMMA.}
\label{tab:dn0p2p0}
\end{table*}%
\begin{table*}
\center
\begin{tabular}{lllll}
\hline
 $1/h$ & $\tau$ & $\errorXx$ & $\errorUu2$ & $\LerrorPp$ \\
 \hline % the following numbers are for N_f = N_c
   3  & $10^{-2}$ & 5.0366e-03 & 1.5689e-02 & 4.3367e-01 \\ % & 7.4382e-01 \\
   6  & $10^{-3}$ & 7.4242e-04 & 6.2453e-03 & 2.9131e-01 \\ % & 4.9838e-01 \\
   12 & $10^{-4}$ & 3.8317e-04 & 3.3759e-03 & 2.0658e-01 \\ % & 7.7516e-01 \\
\hline
\end{tabular}
\caption{($\alpha = 0.15$, $\mu=\gamma=1$)
Expanding bubble problem on $(-1,1)^2 \setminus [-\frac13,\frac13]^2$ 
over the time interval $[0,1]$ 
for the P2--(P1+P0) element without \XFEMGAMMA.}
\label{tab:dn0p2p10}
\end{table*}%

\begin{table*}
\center
\begin{tabular}{llllll}
\hline
 $1/h$ & $\tau$ & $\errorXx$ & $\errorUu2$ & $\LerrorPpc$ 
 & $\errorLl$ \\ % & $\errorKk$ \\ 
 \hline % the following numbers are for N_f = N_c
   3  & $10^{-2}$ & 3.1597e-03 & 1.2035e-02 & 2.3486e-01 & 7.1367e-01 \\ 
   6  & $10^{-3}$ & 2.1580e-04 & 1.9889e-03 & 1.9056e-02 & 2.4283e-02 \\
   12 & $10^{-4}$ & 2.9774e-05 & 2.7515e-04 & 3.8804e-03 & 3.8788e-03 \\
\hline
\end{tabular}
\caption{($\alpha = 0.15$, $\mu=\gamma=1$)
Expanding bubble problem on $(-1,1)^2 \setminus [-\frac13,\frac13]^2$ 
over the time interval $[0,1]$ 
for the P2--P1 element with \XFEMGAMMA.}
\label{tab:dn0p2p1XFEM}
\end{table*}%
\begin{table*}
\center
\begin{tabular}{llllll}
\hline
 $1/h$ & $\tau$ & $\errorXx$ & $\errorUu2$ & $\LerrorPpc$ 
 & $\errorLl$ \\ % & $\errorKk$ \\ 
 \hline % the following numbers are for N_f = N_c
   3  & $10^{-2}$ & 2.2001e-03 & 6.0761e-03 & 1.0310e-01 & 2.2378e-01 \\ 
   6  & $10^{-3}$ & 1.6084e-04 & 9.7213e-04 & 5.6229e-03 & 3.0091e-02 \\
   12 & $10^{-4}$ & 2.6827e-05 & 1.2680e-04 & 9.1875e-04 & 3.1451e-03 \\
\hline
\end{tabular}
\caption{($\alpha = 0.15$, $\mu=\gamma=1$)
Expanding bubble problem on $(-1,1)^2 \setminus [-\frac13,\frac13]^2$ 
over the time interval $[0,1]$ 
for the P2--P0 element with \XFEMGAMMA.}
\label{tab:dn0p2p0XFEM}
\end{table*}%
\begin{table*}
\center
\begin{tabular}{llllll}
\hline
 $1/h$ & $\tau$ & $\errorXx$ & $\errorUu2$ & $\LerrorPpc$ 
 & $\errorLl$ \\ % & $\errorKk$ \\ 
 \hline % the following numbers are for N_f = N_c
   3  & $10^{-2}$ & 2.9540e-03 & 1.6556e-02 & 7.2662e-01 & 2.2105e-00 \\ 
   6  & $10^{-3}$ & 2.2475e-04 & 2.2898e-03 & 3.8525e-02 & 1.1152e-01 \\
   12 & $10^{-4}$ & 3.0103e-05 & 3.0178e-04 & 6.2078e-03 & 2.2342e-02 \\
\hline
\end{tabular}
\caption{($\alpha = 0.15$, $\mu=\gamma=1$)
Expanding bubble problem on $(-1,1)^2 \setminus [-\frac13,\frac13]^2$ 
over the time interval $[0,1]$ 
for the P2--(P1+P0) element with \XFEMGAMMA.}
\label{tab:dn0p2p10XFEM}
\end{table*}%

\begin{figure*}
\center
\mbox{
\includegraphics[angle=-0,width=0.33\textwidth]{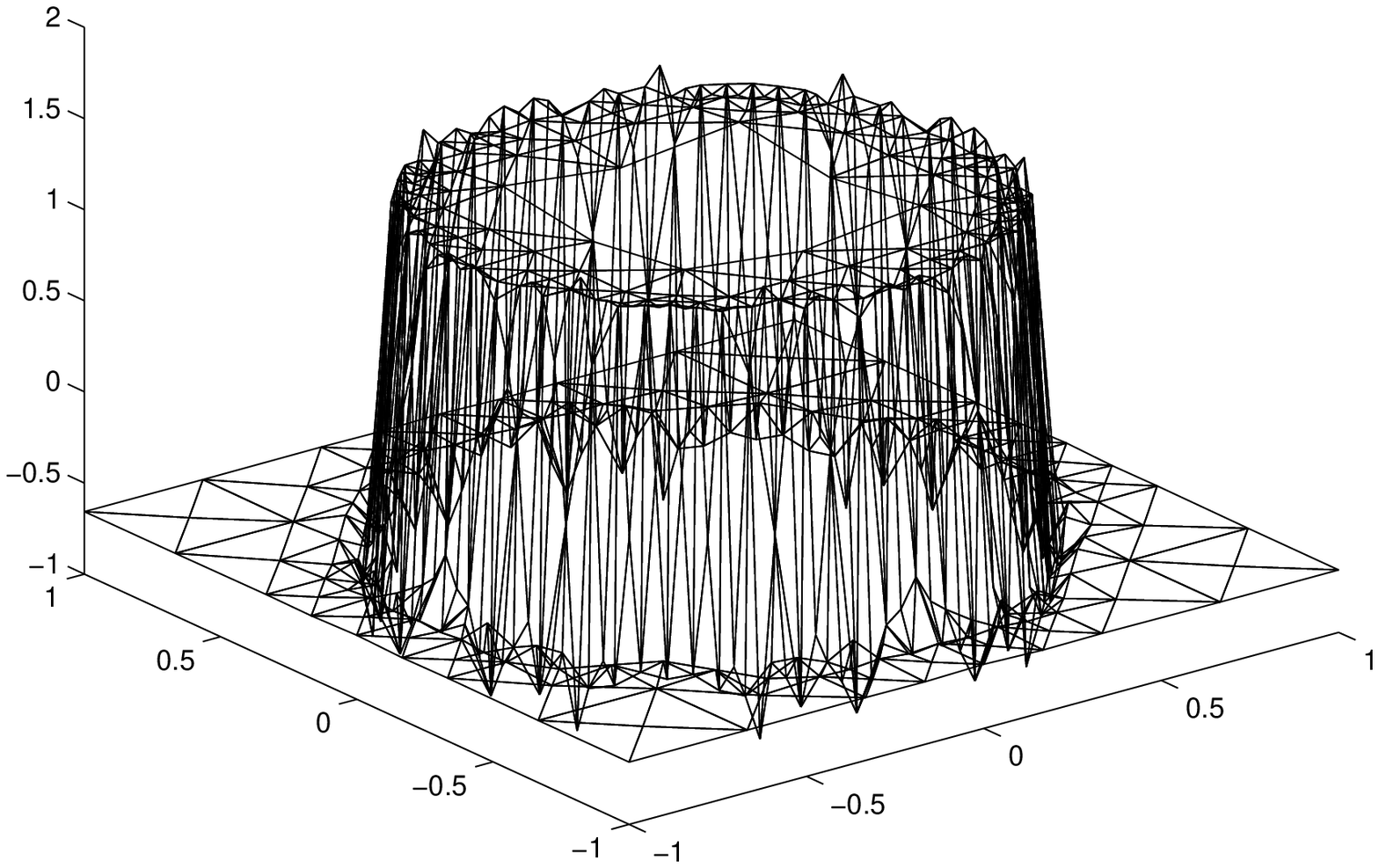}
\includegraphics[angle=-0,width=0.33\textwidth]{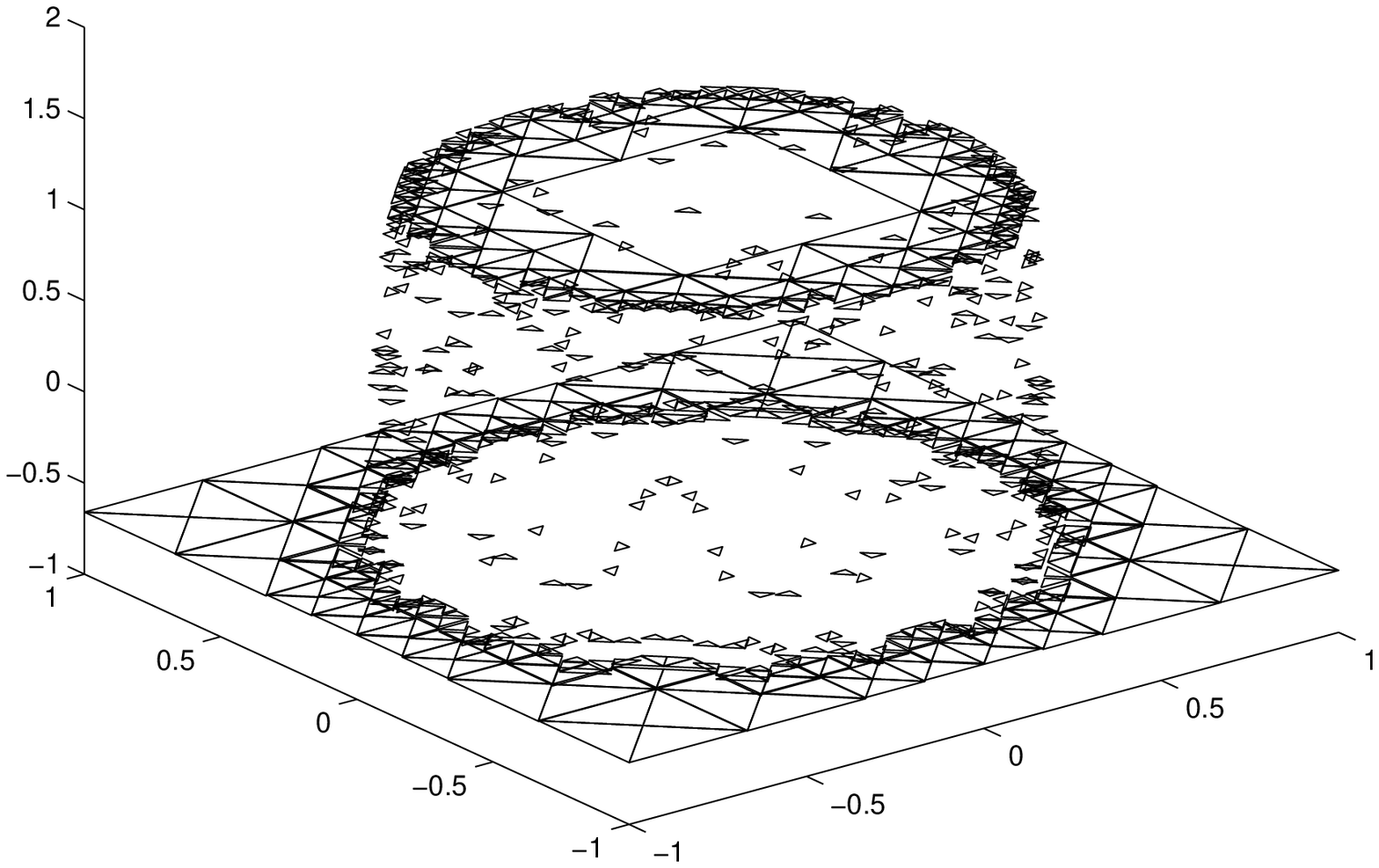}
\includegraphics[angle=-0,width=0.33\textwidth]{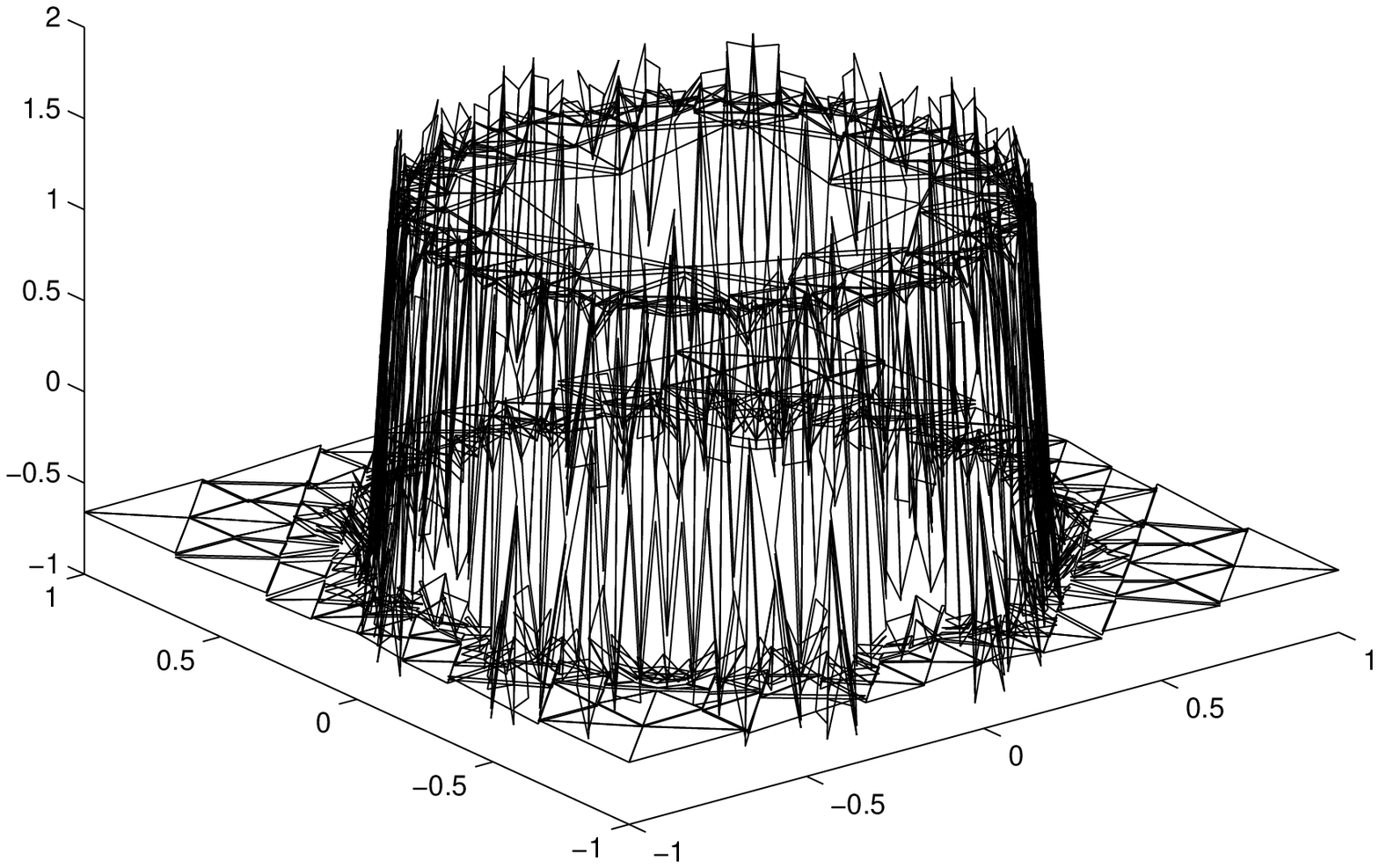}}
\mbox{
\includegraphics[angle=-0,width=0.33\textwidth]{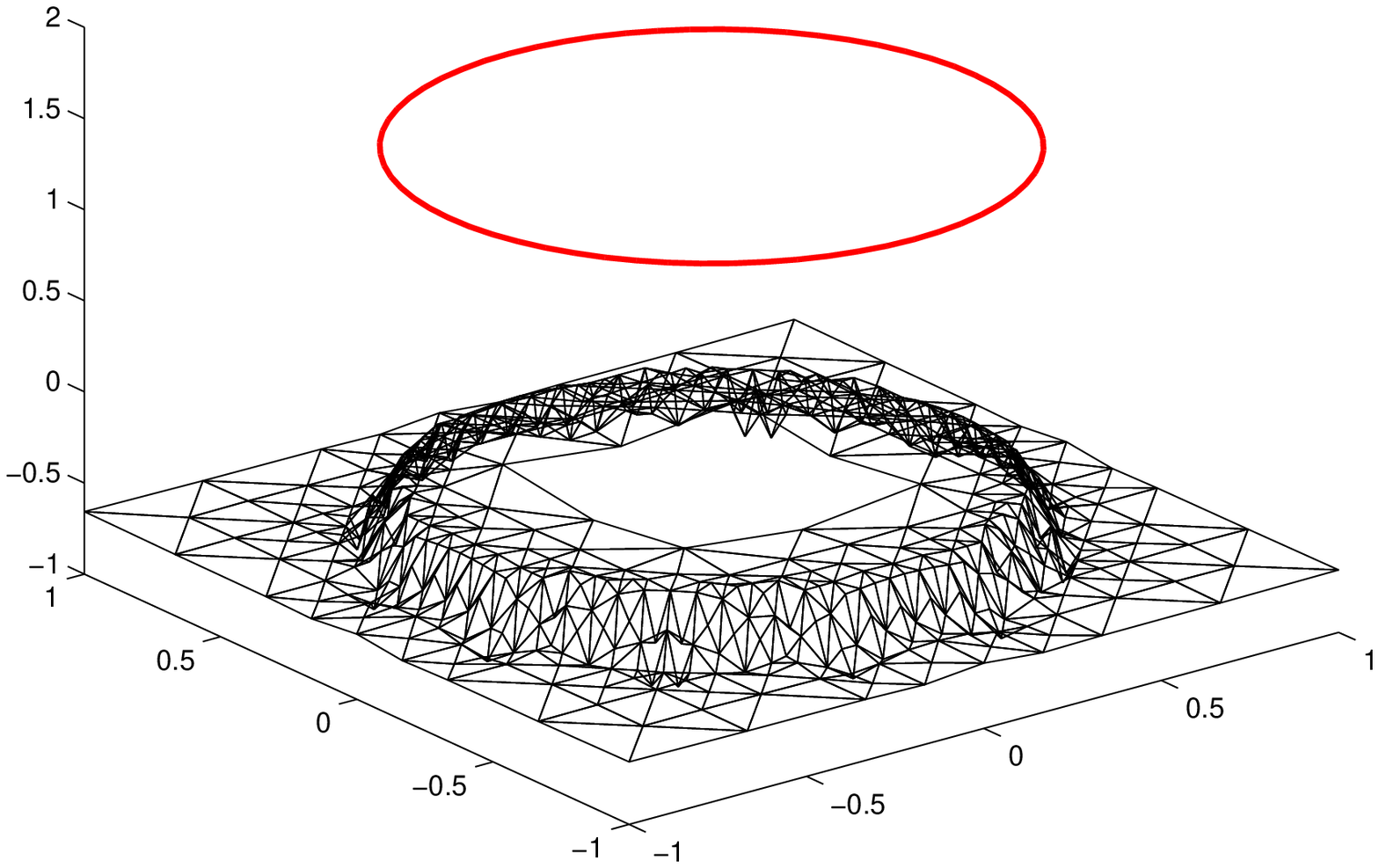}
\includegraphics[angle=-0,width=0.33\textwidth]{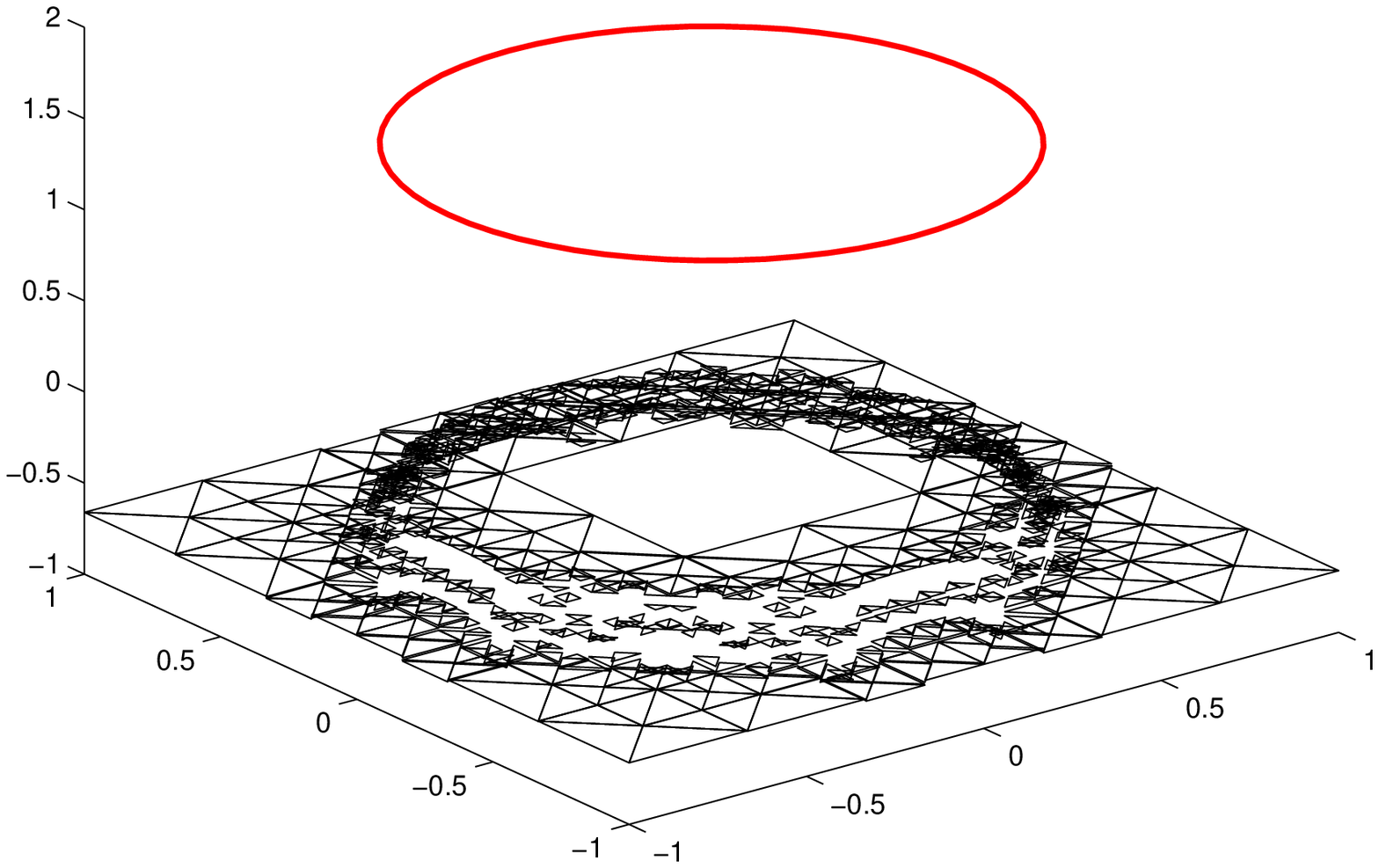}
\includegraphics[angle=-0,width=0.33\textwidth]{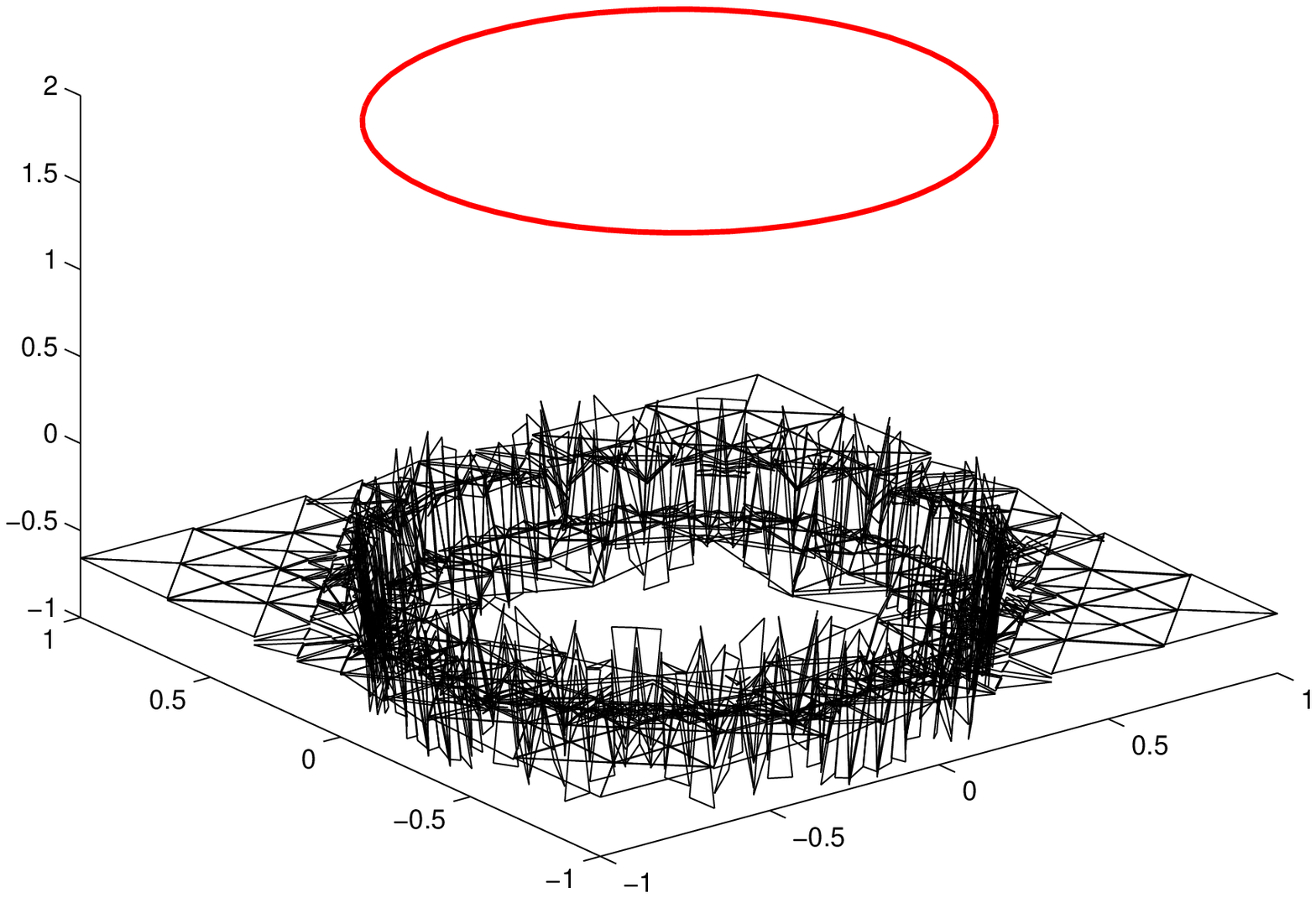}}
\caption{($\alpha = 0.15$, $\mu_+ = 10\,\mu_- = \gamma=1$)
Pressure plots at time $T=1$ for the expanding bubble problem. The
pressure spaces are P1, P0 and P1+P0 without and with \XFEMGAMMA.}
\label{fig:mudiff_pressures}
\end{figure*}

For our final set of convergence experiments we fix 
$\Omega = (-1,1)^2 \setminus [-\frac13,\frac13]^2$ and use the parameters
$$
\alpha = 0.15 \quad\text{and}\quad \mu_+ = 10\,\mu_- = \gamma = 1
$$
for the true solution (\ref{eq:radialr2},b). 
Some errors for our approximation (\ref{eq:HGa}--d) without \XFEMGAMMA\ 
can be seen in Tables~\ref{tab:mudn0p2p1}--\ref{tab:mudn0p2p10}.
Here we always choose the spatial discretization parameters such that 
$h_c = 8\,h_f$ and $h^m_\Gamma \approx h_f$.
\begin{table*}
\center
\begin{tabular}{lllll}
\hline
 $1/h_{f}$ & $\tau$ & $\errorXx$ & $\errorUu2$ & $\LerrorPp$ \\
 \hline % the following numbers are for N_f = 8 N_c
   24 & $10^{-2}$ & 2.1219e-03 & 2.7166e-02 & 3.1792e-01 \\ % & 1.7225e-01 \\  
   48 & $10^{-3}$ & 1.2865e-03 & 1.4738e-02 & 2.2452e-01 \\ % & 1.5813e-01 \\
   96 & $10^{-4}$ & 6.7906e-04 & 9.0482e-03 & 1.5907e-01 \\
\hline
\end{tabular}
\caption{($\alpha = 0.15$, $\mu_+ = 10\,\mu_- = \gamma=1$)
Expanding bubble problem on $(-1,1)^2 \setminus [-\frac13,\frac13]^2$ 
over the time interval $[0,1]$ 
for the P2--P1 element without \XFEMGAMMA.}
\label{tab:mudn0p2p1}
\end{table*}%
\begin{table*}
\center
\begin{tabular}{lllll}
\hline
 $1/h_{f}$ & $\tau$ & $\errorXx$ & $\errorUu2$ & $\LerrorPp$ \\
 \hline % the following numbers are for N_f = 8 N_c
   24 & $10^{-2}$ & 3.2325e-03 & 2.2733e-02 & 3.6175e-01 \\ % & 1.1928e-01 \\  
   48 & $10^{-3}$ & 1.7798e-03 & 1.1853e-02 & 2.5000e-01 \\ % & 1.3457e-01 \\
   96 & $10^{-4}$ & 9.3414e-04 & 7.3246e-03 & 1.7635e-01 \\
\hline
\end{tabular}
\caption{($\alpha = 0.15$, $\mu_+ = 10\,\mu_- = \gamma=1$)
Expanding bubble problem on $(-1,1)^2 \setminus [-\frac13,\frac13]^2$ 
over the time interval $[0,1]$ 
for the P2--P0 element without \XFEMGAMMA.}
\label{tab:mudn0p2p0}
\end{table*}%
\begin{table*}
\center
\begin{tabular}{lllll}
\hline
 $1/h_{f}$ & $\tau$ & $\errorXx$ & $\errorUu2$ & $\LerrorPp$ \\
 \hline % the following numbers are for N_f = 8 N_c
   24 & $10^{-2}$ & 3.2486e-04 & 1.7170e-02 & 2.8290e-01 \\ % & 8.3711e-02 \\  
   48 & $10^{-3}$ & 1.6914e-04 & 1.0384e-02 & 1.9551e-01 \\ % & 2.6046e-02 \\
   96 & $10^{-4}$ & 1.1920e-04 & 6.6529e-03 & 1.3816e-01 \\
\hline
\end{tabular}
\caption{($\alpha = 0.15$, $\mu_+ = 10\,\mu_- = \gamma=1$)
Expanding bubble problem on $(-1,1)^2 \setminus [-\frac13,\frac13]^2$ 
over the time interval $[0,1]$ 
for the P2--(P1+P0) element without \XFEMGAMMA.}
\label{tab:mudn0p2p10}
\end{table*}%
The same convergence experiments now for the pressure spaces enriched with
\XFEMGAMMA\ can be found in 
Tables~\ref{tab:mudn0p2p1XFEM}--\ref{tab:mudn0p2p10XFEM}.
\begin{table*}
\center
\begin{tabular}{lllll}
\hline
 $1/h_{f}$ & $\tau$ & $\errorXx$ & $\errorUu2$ & $\LerrorPp$ \\
 \hline % the following numbers are for N_f = 8 N_c
   24 & $10^{-2}$ & 7.7759e-04 & 1.8081e-02 & 1.3118e-01 \\ % & 6.1157e-02 \\  
   48 & $10^{-3}$ & 1.2812e-04 & 9.7040e-03 & 9.0830e-02 \\ % & 2.6215e-02 \\
   96 & $10^{-4}$ & 2.9108e-05 & 6.2708e-03 & 6.2309e-02 \\
\hline
\end{tabular}
\caption{($\alpha = 0.15$, $\mu_+ = 10\,\mu_- = \gamma=1$)
Expanding bubble problem on $(-1,1)^2 \setminus [-\frac13,\frac13]^2$ 
over the time interval $[0,1]$ 
for the P2--P1 element with \XFEMGAMMA.}
\label{tab:mudn0p2p1XFEM}
\end{table*}%
\begin{table*}
\center
\begin{tabular}{lllll}
\hline
 $1/h_{f}$ & $\tau$ & $\errorXx$ & $\errorUu2$ & $\LerrorPp$ \\
 \hline % the following numbers are for N_f = 8 N_c
   24 & $10^{-2}$ & 7.6544e-04 & 1.6390e-02 & 1.3925e-01 \\ % & 6.9919e-02 \\  
   48 & $10^{-3}$ & 1.3313e-04 & 1.0512e-02 & 9.6750e-02 \\ % & 4.0644e-02 \\
   96 & $10^{-4}$ & 3.0199e-05 & 6.8224e-03 & 6.6541e-02 \\
\hline
\end{tabular}
\caption{($\alpha = 0.15$, $\mu_+ = 10\,\mu_- = \gamma=1$)
Expanding bubble problem on $(-1,1)^2 \setminus [-\frac13,\frac13]^2$ 
over the time interval $[0,1]$ 
for the P2--P0 element with \XFEMGAMMA.}
\label{tab:mudn0p2p0XFEM}
\end{table*}%
\begin{table*}
\center
\begin{tabular}{lllll}
\hline
 $1/h_{f}$ & $\tau$ & $\errorXx$ & $\errorUu2$ & $\LerrorPp$ \\
 \hline % the following numbers are for N_f = 8 N_c
   24 & $10^{-2}$ & 8.4735e-04 & 1.5378e-02 & 1.5429e-01 \\ % & 8.3444e-02 \\  
   48 & $10^{-3}$ & 1.5376e-04 & 9.7116e-03 & 1.0642e-01 \\ % & 2.7408e-02 \\
   96 & $10^{-4}$ & 3.5443e-05 & 6.4129e-03 & 7.2764e-02 \\
\hline
\end{tabular}
\caption{($\alpha = 0.15$, $\mu_+ = 10\,\mu_- = \gamma=1$)
Expanding bubble problem on $(-1,1)^2 \setminus [-\frac13,\frac13]^2$ 
over the time interval $[0,1]$ 
for the P2--(P1+P0) element with \XFEMGAMMA.}
\label{tab:mudn0p2p10XFEM}
\end{table*}%
We visualize the final pressures for the coarsest runs in 
Tables~\ref{tab:mudn0p2p1}--\ref{tab:mudn0p2p10XFEM} in 
Figure~\ref{fig:mudiff_pressures}.
Here in the case of the enrichment \XFEMGAMMA\ being used, we plot the pressure
parts $P^M_c$ and $\lambda^M\,\charfcn{\Omega_-^{M-1}}$ separately.
We recall from (\ref{eq:radialup2}) that the jump $\lambda(T)$ in the pressure
is made up from two components. A jump of $0.55^{-\frac12} \approx 1.35$ due 
to the curvature, and a jump of $0.27 / 0.55 \approx 0.49$ due to the
difference in $\mu$. When the \XFEMGAMMA\ enrichment of the pressure space is
used, then the discretizations will treat these two jumps differently. As can
be seen from (\ref{eq:HGa}), the jump due to curvature can be absorbed by
$\charfcn{\Omega^m_-}$, while the jump in $\mu$, recall (\ref{eq:rhoma}), is
left to the standard bulk pressure space.
It appears from our plots in Figure~\ref{fig:mudiff_pressures}
that this is indeed what happens in practice.

\subsection{Numerical results in 3d}\label{numexpt3d}

Similarly to the simulation in Figures~\ref{fig:usp} and \ref{fig:uspplot}, we
show that our approximation (\ref{eq:HGa}--d) with \XFEMGAMMA\ naturally
eliminates spurious velocities also in three space dimensions. 
We recall that in our experiments in 2d it was necessary (and sufficient) for 
the polygonal curve
$\Gamma^m$ to be an equidistributed approximation to a circle in order to 
admit a constant discrete curvature $\kappa^{m+1} = \overline\kappa \in \R$. 
This leads to a
constant pressure jump across $\Gamma^m$, which can be picked up by our
extended finite element function $\charfcn{\Omega^m_-}$.

\begin{figure*}
\center
\mbox{
\includegraphics[angle=-90,width=0.26\textwidth]{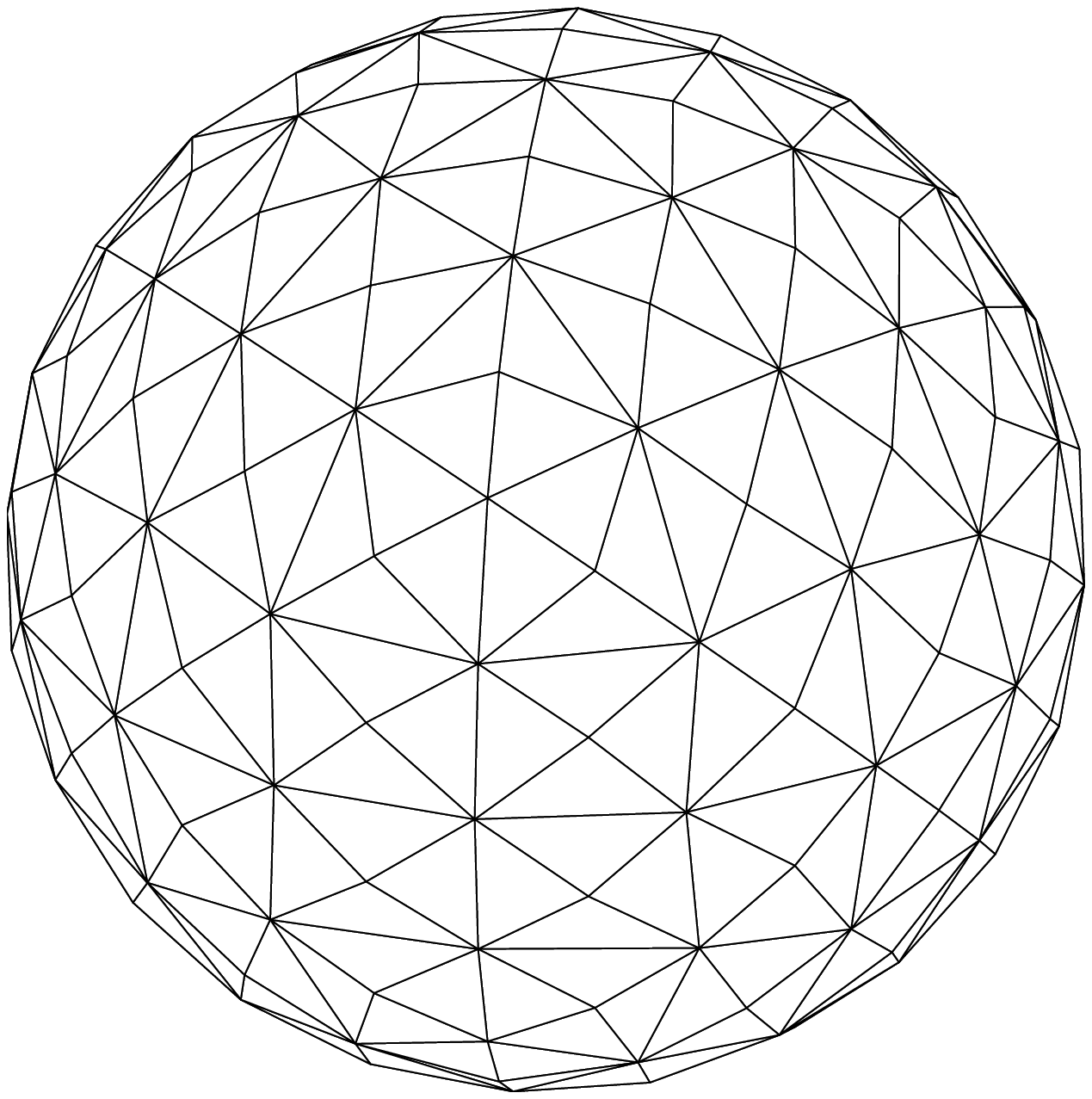} \quad
\includegraphics[angle=-90,width=0.26\textwidth]{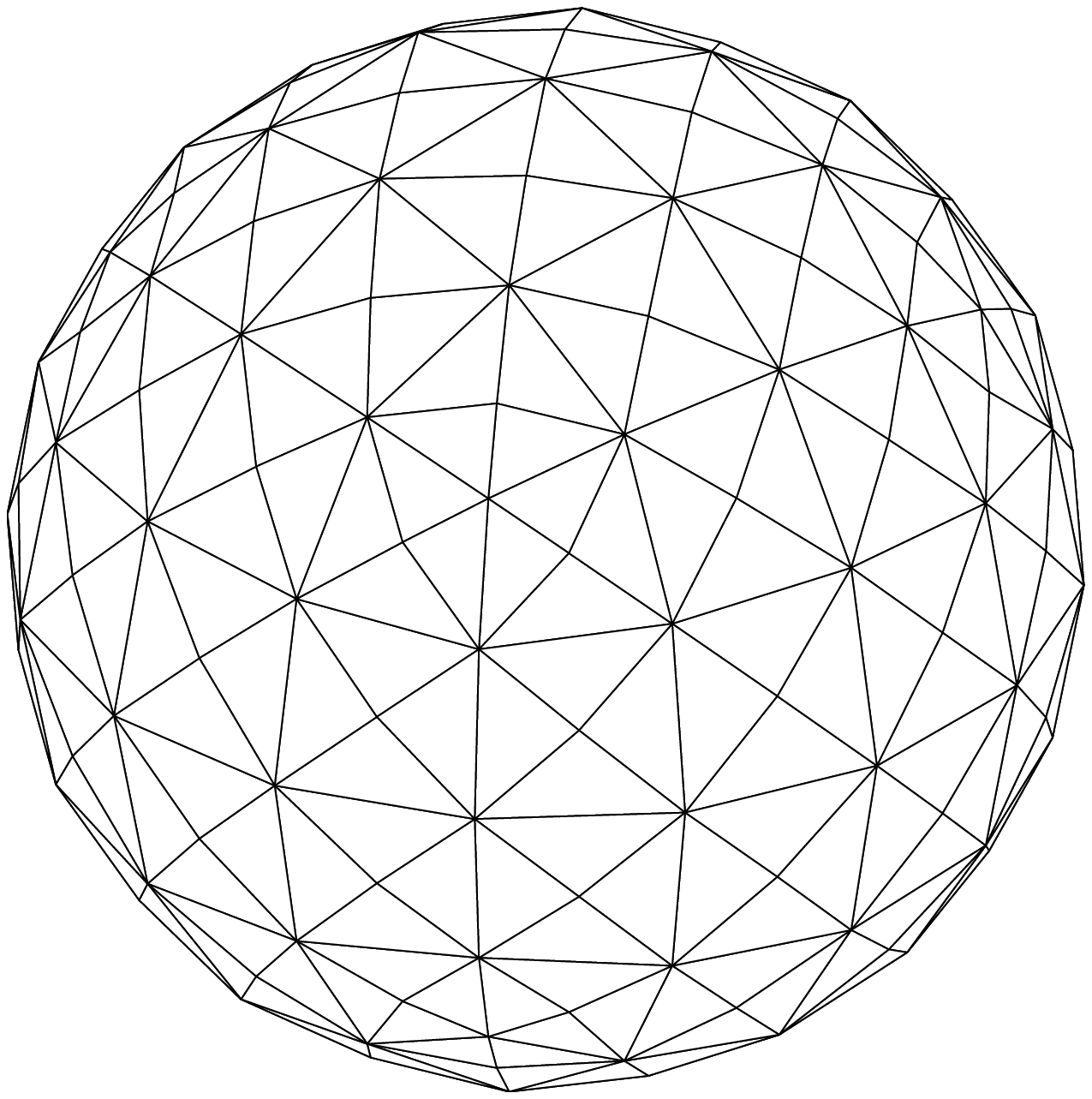}
\quad
\includegraphics[angle=-90,width=0.4\textwidth]{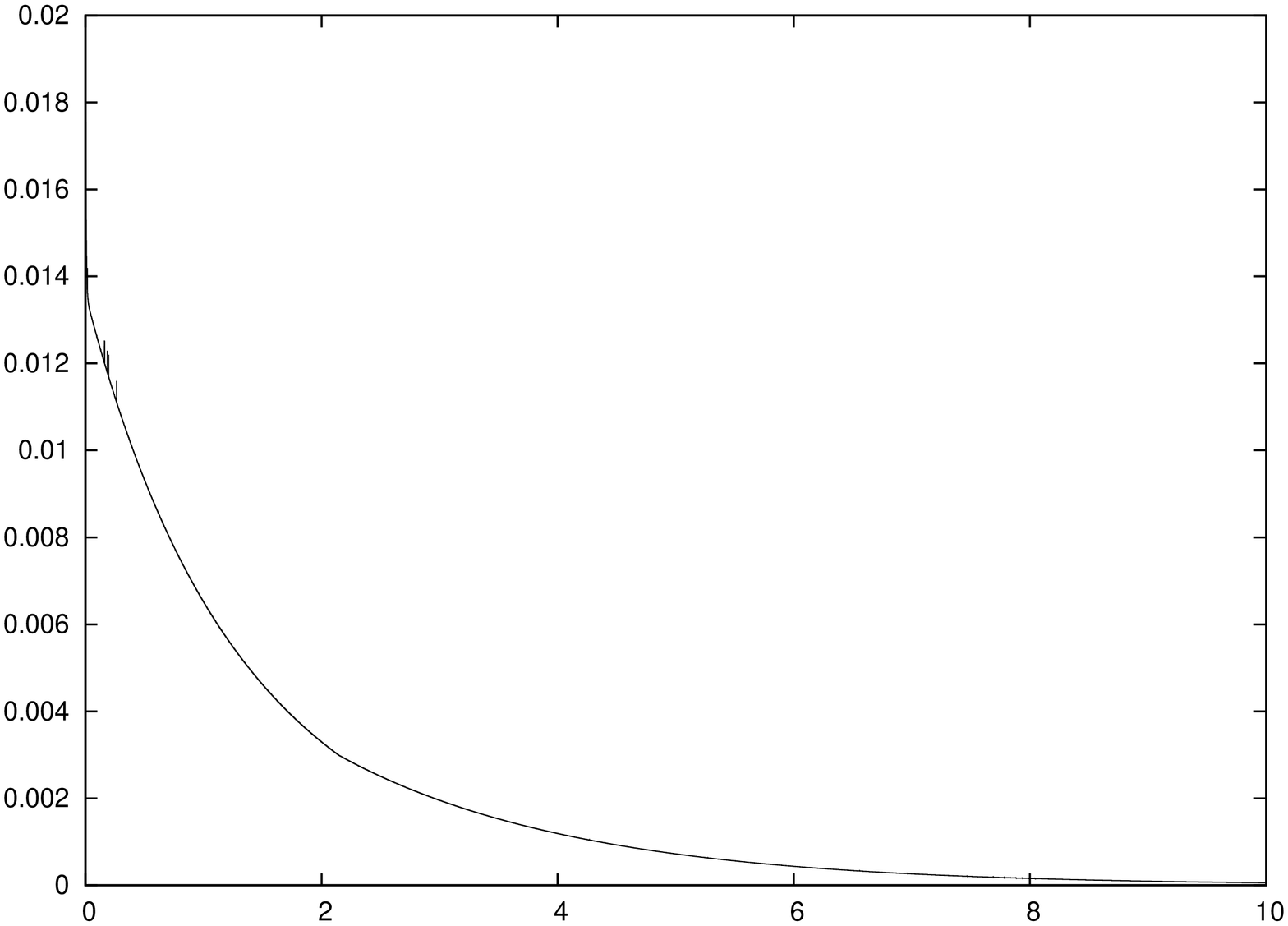}}
\caption{($\mu=\gamma=1$)
The discrete interface $\Gamma^m$ at times $t=0,\,5$. 
Here we use the P2--P1 element with \XFEMGAMMA.
On the right a plot of $\|\vec U^m\|_{L^\infty(\Omega)}$
over the time interval $[0,10]$ 
}
\label{fig:3dsphere}
\end{figure*}

In 3d a necessary condition for 
$\kappa^{m+1}$ to be constant requires $\Gamma^m$ to be a conformal polyhedral 
surface, recall Lemma~\ref{lem:stat2}. 
Here we note that the
tangential movement of vertices induced by (\ref{eq:HGd}) leads to stationary
solutions being conformal polyhedral surfaces, see \cite{gflows3d} for details.
We demonstrate this, and the fact that this property leads to the elimination
of spurious velocities, 
with a numerical simulation for (\ref{eq:HGa}--d) with \XFEMGAMMA\ on
$\Omega=(-1,1)^3$ with $\mu = \gamma = 1$. 
For the bulk mesh discretization we use an adaptive mesh with 
$h_c = 16\,h_f = 2^{\frac32}$.
We start with a standard
triangulation of the sphere $\Gamma(0)$, which is made up of $378$ elements,
so that $h^0_\Gamma = 0.2$, 
and compute the evolution over the time
interval $[0,10]$ with the time step size $\tau = 10^{-4}$. We visualize
$\Gamma^0$ and $\Gamma^M$ in Figure~\ref{fig:3dsphere}, with a plot of
$\|\vec U^m\|_{L^\infty(\Omega)}$ over time given on the right of the same 
figure.
We observe that the initial triangulation of the sphere evolves towards a
numerical steady state for our approximation \mbox{(\ref{eq:HGa}--d)}. Such a
numerical steady state will be given by a conformal polyhedral surface, 
recall (\ref{eq:conformal}). % see \cite{gflows3d}. 
Similarly to \cite[Fig.\ 11]{gflows3d} it can be observed
that the final triangulation in Figure~\ref{fig:3dsphere}
exhibits many groups of two, four and eight triangles that form 
``curved squares'', as well as groups of six and twelve triangles that form
``curved equilateral triangles''. These are typical for conformal polyhedral 
surfaces, and we conjecture that %this 
conformal polyhedral approximations of the sphere have
constant discrete mean curvature, i.e.\ satisfy (\ref{eq:constcurv}). 
In fact we note that for the simulation at hand 
the extremal values of $\kappa^1$ are given by
$-6.47$ and $-2.51$, while for $\kappa^M$ they are $-4.07$ and $-4.05$, 
i.e.\ $\kappa^M$ is close to being constant.

\begin{figure*}
\center
\mbox{\hspace*{-12mm}
\includegraphics[angle=-90,width=0.4\textwidth]{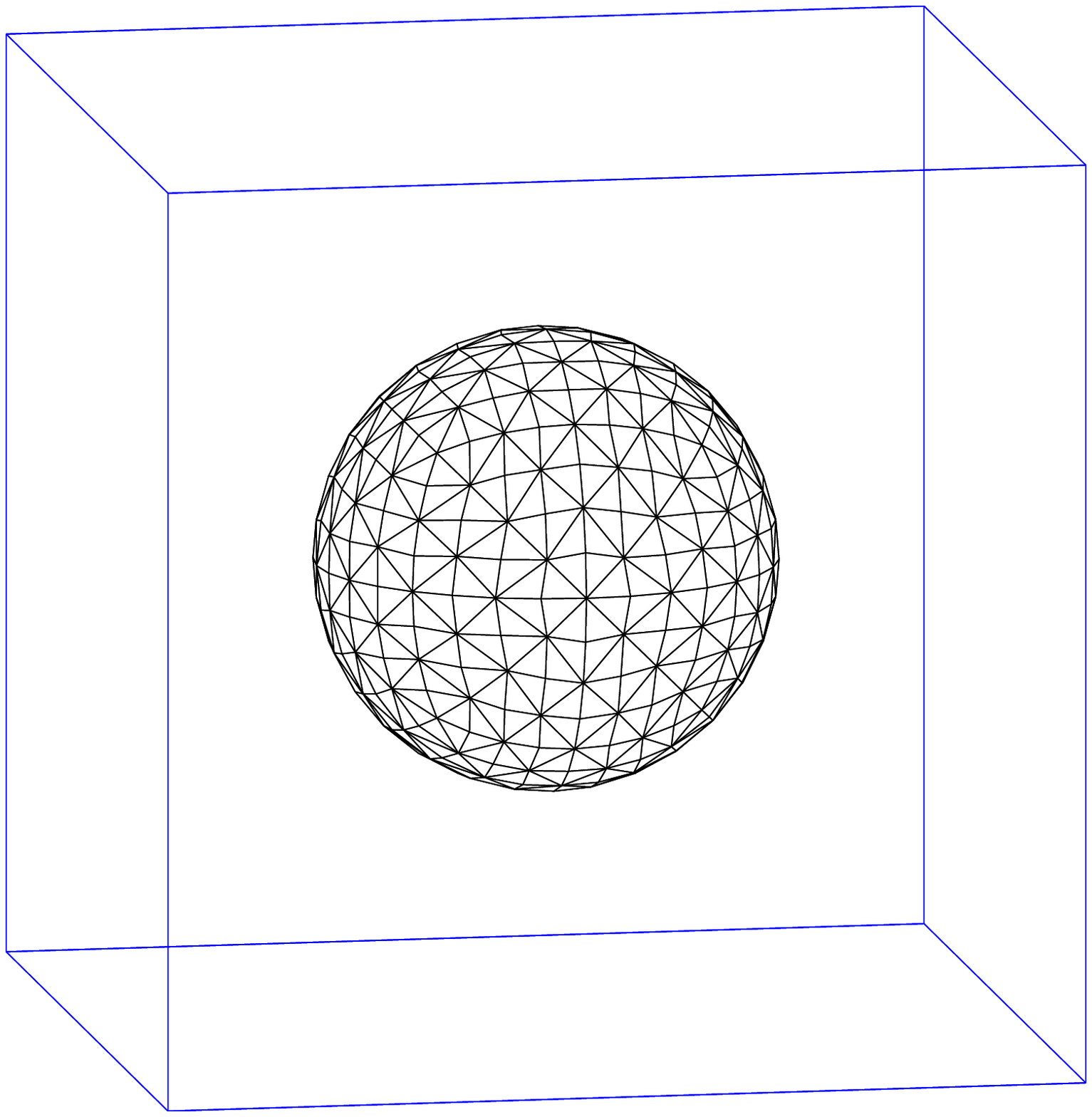} 
\includegraphics[angle=-90,width=0.4\textwidth]{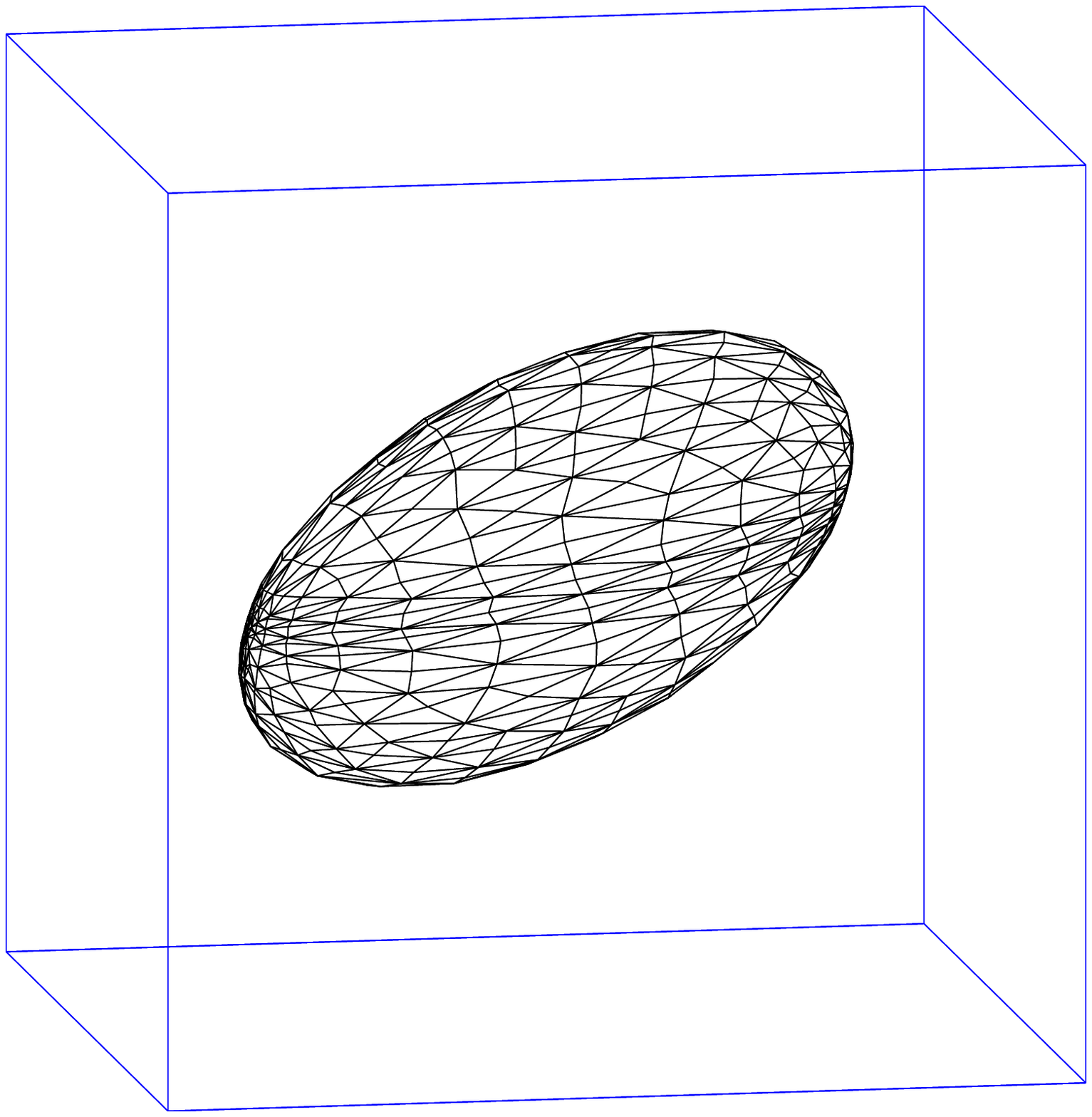}} 
\caption{($\mu=1$, $\gamma = 3$)
The discrete interface $\Gamma^m$ at times $t=0,\,4.2$.
Here we use the the alternative curvature scheme (\ref{eq:GDa}--d) with the
P2--P1 element and \XFEMGAMMA.}
\label{fig:shear3ddziuk}
\end{figure*}

We end this section with a shear flow experiment that is motivated by 
similar simulations in \cite{LiYLSJK13}. Here we take $\Omega = (-1,1)^3$ and
for $\vec u$
prescribe the inhomogeneous Dirichlet boundary condition
$\vec g(\vec z) = (z_3, 0, 0)^T$ on $\partial\Omega$. 
We also let $\mu=1$ and $\gamma = 3$. 
For the discretization parameters we fix $h_c = 8\,h_f = 2^{-\frac{1}2}$ 
and $\tau=10^{-2}$, and we use a triangulation of
$\Gamma^0$ that is made up of $768$ elements, and so 
$h^0_\Gamma = 0.14$.
At first we present the numerical results for the alternative formulation
(\ref{eq:GDa}--d), see Figure~\ref{fig:shear3ddziuk}, where we plot the
discrete interface $\Gamma^m$ at times $t=0$ and $t=4.2$. We note the elongated
elements and the high density of vertices in certain regions. In fact, the
simulation breaks down soon after because of the coalescence of mesh points.

\begin{figure*}
\center
\mbox{\hspace*{-12mm}
\includegraphics[angle=-90,width=0.4\textwidth]{figures/shearflow3d_t0} 
\includegraphics[angle=-90,width=0.4\textwidth]{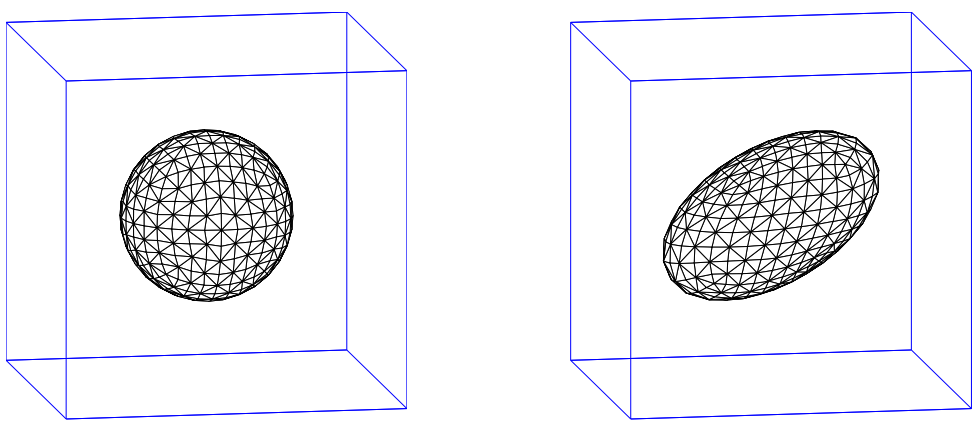}} 
\caption{($\mu=1$, $\gamma = 3$)
The discrete interface $\Gamma^m$ at times $t=0,\,5$. 
Here we use our scheme (\ref{eq:HGa}--d) with the P2--P1 element and
\XFEMGAMMA.}
\label{fig:shear3d}
\end{figure*}

The same simulation for our approximation (\ref{eq:HGa}--d) can be seen in
Figure~\ref{fig:shear3d}. Here 
the interface at the final time $T=5$ is close to being a numerical steady
state. The ``sphericity'' of $\Gamma^M$, see \cite{Wadell33}, is given
by $\pi^\frac13\,[6\,\mathcal{L}^3(\Omega_-^M)]^\frac23\,
[\mathcal{H}^{2}(\Gamma^M)]^{-1} = 0.957$.
We observe an excellent mesh quality throughout the evolution, in contrast to
the elongated elements that can be seen in Figure~\ref{fig:shear3ddziuk}
and in e.g.\ \cite[Fig.\ 14]{LiYLSJK13}.

\section*{Conclusions}
We have presented a novel front-tracking method for 
viscous incompressible two-phase flow which can be
shown to be stable. The numerical method couples a parametric finite element
approximation of the interface with a standard finite element approximation of
the Stokes equations in the bulk. Here the bulk mesh may be chosen to be either
fitted (or adapted) to the interface, or it can be totally independent of
the interface mesh. In the latter case, we introduce an XFEM approach to
guarantee that our scheme conserves the volumes of the two phases.

The stability of the proposed method implies that the well-known static bubble
test problem can be computed exactly. More generally, we can show that our
scheme admits time-independent discrete solutions, 
which all have the property that the velocity is zero.                     
This means, in particular, that it is possible to eliminate spurious
velocities for stationary solutions.

The second prominent feature of our numerical method is the excellent mesh
quality of the interface approximation. This is induced by an inherent discrete
tangential motion of the vertices that make up the discrete interface. In
particular, for a semidiscrete continuous-in-time variant of our scheme it can 
be shown that the discrete interfaces are equidistributed polygonal curves
($d=2$) and conformal polyhedral surfaces ($d=3$), respectively.

\def\soft#1{\leavevmode\setbox0=\hbox{h}\dimen7=\ht0\advance \dimen7
  by-1ex\relax\if t#1\relax\rlap{\raise.6\dimen7
  \hbox{\kern.3ex\char'47}}#1\relax\else\if T#1\relax
  \rlap{\raise.5\dimen7\hbox{\kern1.3ex\char'47}}#1\relax \else\if
  d#1\relax\rlap{\raise.5\dimen7\hbox{\kern.9ex \char'47}}#1\relax\else\if
  D#1\relax\rlap{\raise.5\dimen7 \hbox{\kern1.4ex\char'47}}#1\relax\else\if
  l#1\relax \rlap{\raise.5\dimen7\hbox{\kern.4ex\char'47}}#1\relax \else\if
  L#1\relax\rlap{\raise.5\dimen7\hbox{\kern.7ex
  \char'47}}#1\relax\else\message{accent \string\soft \space #1 not
  defined!}#1\relax\fi\fi\fi\fi\fi\fi}


\begin{thebibliography}{10}

\bibitem{AndersonMW98}
{\sc D.~M. Anderson, G.~B. McFadden, and A.~A. Wheeler}, {\em Diffuse-interface
  methods in fluid mechanics}, in Annual review of fluid mechanics, {V}ol. 30,
  Annual Reviews, Palo Alto, CA, 1998, pp.~139--165.

\bibitem{AusasBI12}
{\sc R.~F. Ausas, G.~C. Buscaglia, and S.~R. Idelsohn}, {\em A new enrichment
  space for the treatment of discontinuous pressures in multi-fluid flows},
  Internat. J. Numer. Methods Fluids, 70 (2012), pp.~829--850.

\bibitem{Bansch01}
{\sc E.~B{\"a}nsch}, {\em Finite element discretization of the
  {N}avier--{S}tokes equations with a free capillary surface}, Numer. Math., 88
  (2001), pp.~203--235.

\bibitem{triplejMC}
{\sc J.~W. Barrett, H.~Garcke, and R.~N\"urnberg}, {\em On the variational
  approximation of combined second and fourth order geometric evolution
  equations}, SIAM J. Sci. Comput., 29 (2007), pp.~1006--1041.

\bibitem{triplej}
\leavevmode\vrule height 2pt depth -1.6pt width 23pt, {\em A parametric finite
  element method for fourth order geometric evolution equations}, J. Comput.
  Phys., 222 (2007), pp.~441--462.

\bibitem{gflows3d}
\leavevmode\vrule height 2pt depth -1.6pt width 23pt, {\em On the parametric
  finite element approximation of evolving hypersurfaces in {${\mathbb R}^3$}},
  J. Comput. Phys., 227 (2008), pp.~4281--4307.

\bibitem{willmore}
\leavevmode\vrule height 2pt depth -1.6pt width 23pt, {\em Parametric
  approximation of {W}illmore flow and related geometric evolution equations},
  SIAM J. Sci. Comput., 31 (2008), pp.~225--253.

\bibitem{ejam3d}
\leavevmode\vrule height 2pt depth -1.6pt width 23pt, {\em Finite element
  approximation of coupled surface and grain boundary motion with applications
  to thermal grooving and sintering}, European J. Appl. Math., 21 (2010),
  pp.~519--556.

\bibitem{dendritic}
\leavevmode\vrule height 2pt depth -1.6pt width 23pt, {\em On stable parametric
  finite element methods for the {S}tefan problem and the {M}ullins--{S}ekerka
  problem with applications to dendritic growth}, J. Comput. Phys., 229 (2010),
  pp.~6270--6299.

\bibitem{jcg}
\leavevmode\vrule height 2pt depth -1.6pt width 23pt, {\em Numerical
  computations of faceted pattern formation in snow crystal growth}, Phys. Rev.
  E, 86 (2012), p.~011604.

\bibitem{crystal}
\leavevmode\vrule height 2pt depth -1.6pt width 23pt, {\em Finite element
  approximation of one-sided {S}tefan problems with anisotropic, approximately
  crystalline, {G}ibbs--{T}homson law}, Adv. Differential Equations, 18 (2013),
  pp.~383--432.

\bibitem{fluidfbp}
\leavevmode\vrule height 2pt depth -1.6pt width 23pt, {\em A stable parametric
  finite element discretization of two-phase {N}avier--{S}tokes flow}, 2013.
\newblock \url{http://arxiv.org/abs/1308.3335}.

\bibitem{BoffiCGG12}
{\sc D.~Boffi, N.~Cavallini, F.~Gardini, and L.~Gastaldi}, {\em Local mass
  conservation of {S}tokes finite elements}, J. Sci. Comput., 52 (2012),
  pp.~383--400.

\bibitem{DeckelnickDE05}
{\sc K.~Deckelnick, G.~Dziuk, and C.~M. Elliott}, {\em Computation of geometric
  partial differential equations and mean curvature flow}, Acta Numer., 14
  (2005), pp.~139--232.

\bibitem{Dziuk91}
{\sc G.~Dziuk}, {\em An algorithm for evolutionary surfaces}, Numer. Math., 58
  (1991), pp.~603--611.

\bibitem{Feng06}
{\sc X.~Feng}, {\em Fully discrete finite element approximations of the
  {N}avier--{S}tokes--{C}ahn--{H}illiard diffuse interface model for two-phase
  fluid flows}, SIAM J. Numer. Anal., 44 (2006), pp.~1049--1072.

\bibitem{FrancoisCDKSW06}
{\sc M.~M. Francois, S.~J. Cummins, E.~D. Dendy, D.~B. Kothe, J.~M. Sicilian,
  and M.~W. Williams}, {\em A balanced-force algorithm for continuous and sharp
  interfacial surface tension models within a volume tracking framework}, J.
  Comput. Phys., 213 (2006), pp.~141--173.

\bibitem{GanesanMT07}
{\sc S.~Ganesan, G.~Matthies, and L.~Tobiska}, {\em On spurious velocities in
  incompressible flow problems with interfaces}, Comput. Methods Appl. Mech.
  Engrg., 196 (2007), pp.~1193--1202.

\bibitem{GanesanT08}
{\sc S.~Ganesan and L.~Tobiska}, {\em An accurate finite element scheme with
  moving meshes for computing 3{D}-axisymmetric interface flows}, Internat. J.
  Numer. Methods Fluids, 57 (2008), pp.~119--138.

\bibitem{GerbeauLB97}
{\sc J.-F. Gerbeau, C.~Le~Bris, and M.~Bercovier}, {\em Spurious velocities in
  the steady flow of an incompressible fluid subjected to external forces},
  Internat. J. Numer. Methods Fluids, 25 (1997), pp.~679--695.

\bibitem{GiraultR86}
{\sc V.~Girault and P.-A. Raviart}, {\em Finite Element Methods for
  {N}avier--{S}tokes}, Springer-Verlag, Berlin, 1986.

\bibitem{GrossR07}
{\sc S.~Gro{\ss} and A.~Reusken}, {\em An extended pressure finite element
  space for two-phase incompressible flows with surface tension}, J. Comput.
  Phys., 224 (2007), pp.~40--58.

\bibitem{GrossR11}
\leavevmode\vrule height 2pt depth -1.6pt width 23pt, {\em Numerical methods
  for two-phase incompressible flows}, vol.~40 of Springer Series in
  Computational Mathematics, Springer-Verlag, Berlin, 2011.

\bibitem{HirtN81}
{\sc C.~W. Hirt and B.~D. Nichols}, {\em Volume of fluid {(VOF)} method for the
  dynamics of free boundaries}, J. Comput. Phys., 39 (1981), pp.~201--225.

\bibitem{Jacqmin99}
{\sc D.~Jacqmin}, {\em Calculation of two-phase {N}avier--{S}tokes flows using
  phase-field modeling}, J. Comput. Phys., 155 (1999), pp.~96--127.

\bibitem{JametTB02}
{\sc D.~Jamet, D.~Torres, and J.~U. Brackbill}, {\em On the theory and
  computation of surface tension: the elimination of parasitic currents through
  energy conservation in the second-gradient method}, J. Comput. Phys., 182
  (2002), pp.~262--276.

\bibitem{LeVequeL97}
{\sc R.~J. LeVeque and Z.~Li}, {\em Immersed interface methods for {S}tokes
  flow with elastic boundaries or surface tension}, SIAM J. Sci. Comput., 18
  (1997), pp.~709--735.

\bibitem{LiYLSJK13}
{\sc Y.~Li, A.~Yun, D.~Lee, J.~Shin, D.~Jeong, and J.~Kim}, {\em
  Three-dimensional volume-conserving immersed boundary model for two-phase
  fluid flows}, Comput. Methods Appl. Mech. Engrg., 257 (2013), pp.~36--46.

\bibitem{OsherF03}
{\sc S.~Osher and R.~Fedkiw}, {\em Level Set Methods and Dynamic Implicit
  Surfaces}, vol.~153 of Applied Mathematical Sciences, Springer-Verlag, New
  York, 2003.

\bibitem{Peskin02}
{\sc C.~S. Peskin}, {\em The immersed boundary method}, Acta Numer., 11 (2002),
  pp.~479--517.

\bibitem{Popinet09}
{\sc S.~Popinet}, {\em An accurate adaptive solver for surface-tension-driven
  interfacial flows}, J. Comput. Phys., 228 (2009), pp.~5838--5866.

\bibitem{PopinetZ99}
{\sc S.~Popinet and S.~Zaleski}, {\em A front-tracking algorithm for the
  accurate representation of surface tension}, Internat. J. Numer. Methods
  Fluids, 30 (1999), pp.~775--793.

\bibitem{RenardyR02}
{\sc Y.~Renardy and M.~Renardy}, {\em P{ROST}: a parabolic reconstruction of
  surface tension for the volume-of-fluid method}, J. Comput. Phys., 183
  (2002), pp.~400--421.

\bibitem{SauerlandF12}
{\sc H.~Sauerland and T.-P. Fries}, {\em The stable {XFEM} for two-phase
  flows}, Comput. \& Fluids,  (2012).
\newblock DOI: 10.1016/j.compfluid.2012.10.017.

\bibitem{Alberta}
{\sc A.~Schmidt and K.~G. Siebert}, {\em Design of Adaptive Finite Element
  Software: The Finite Element Toolbox {ALBERTA}}, vol.~42 of Lecture Notes in
  Computational Science and Engineering, Springer-Verlag, Berlin, 2005.

\bibitem{Sethian}
{\sc J.~A. Sethian}, {\em Level Set Methods and Fast Marching Methods},
  Cambridge University Press, Cambridge, 1999.

\bibitem{SussmanSO94}
{\sc M.~Sussman, P.~Semereka, and S.~Osher}, {\em A level set approach for
  computing solutions to incompressible two-phase flow}, J. Comput. Phys., 114
  (1994), pp.~146--159.

\bibitem{TongW07}
{\sc A.~Y. Tong and Z.~Wang}, {\em A numerical method for capillarity-dominant
  free surface flows}, J. Comput. Phys., 221 (2007), pp.~506--523.

\bibitem{UnverdiT92}
{\sc S.~O. Unverdi and G.~Tryggvason}, {\em A front-tracking method for
  viscous, incompressible multi-fluid flows}, J. Comput. Phys., 100 (1992),
  pp.~25--37.

\bibitem{Wadell33}
{\sc H.~Wadell}, {\em Sphericity and roundness of rock particles}, J. Geol., 41
  (1933), pp.~310--331.

\bibitem{ZahediKK12}
{\sc S.~Zahedi, M.~Kronbichler, and G.~Kreiss}, {\em Spurious currents in
  finite element based level set methods for two-phase flow}, Internat. J.
  Numer. Methods Fluids, 69 (2012), pp.~1433--1456.

\end{thebibliography}
\end{document}